\documentclass[11pt,a4paper]{article}
\usepackage[francais]{babel}
\usepackage[latin1]{inputenc}

\usepackage{amsmath}
\usepackage{amsfonts}
\usepackage{amssymb}
\usepackage{amsthm}
\usepackage{graphicx,epsfig}
\usepackage{amscd}
\usepackage{array}

\newtheorem{Prop}{Proposition}[section]
\newtheorem{Def}[Prop]{D\'efinition}
\newtheorem{Lem}[Prop]{Lemme}
\newtheorem{Thm}[Prop]{Th\'eor\`eme}
\newtheorem{Cor}[Prop]{Corollaire}
\newtheorem*{nThm}{Th\'eor\`eme}

\newcommand{\N}{\mathbb{N}}
\newcommand{\n}{\nabla}
\newcommand{\na}{\nabla ^{*}}
\newcommand{\tr}{\mathrm{tr\,}}
\newcommand{\trs}{\mathrm{tr}_\Sigma\,}
\newcommand{\im}{\mathrm{Im\,}}
\newcommand{\<}{\langle}
\renewcommand{\>}{\rangle}
\renewcommand{\d}{\partial}
\newcommand{\RO}{\mathring{R}}

\renewcommand{\cosh}{\mathrm{ch}}
\renewcommand{\sinh}{\mathrm{sh}}
\renewcommand{\tanh}{\mathrm{th}}

\title{D\'eformations Einstein infinit\'esimales de c\^one-vari\'et\'es hyperboliques \\
\Large Infinitesimal Einstein deformations of hyperbolic cone-manifolds
}
\author{Gr\'egoire Montcouquiol \\ \small{Equipe de Topologie et Dynamique, Laboratoire de 
Math\'ematiques (UMR 8628),
  Universit\'e Paris-Sud XI}}

\begin{document}

\maketitle

\begin{minipage}{16cm}
\small

\noindent{\bf Abstract}

\vskip 0.5\baselineskip
\noindent Starting with a compact hyperbolic cone-manifold of dimension $n\geq 3$, we study the deformations of the metric in order to get Einstein cone-manifolds. If the singular locus is a closed codimension $2$ submanifold and all cone angles are smaller than $2\pi$, we show that there is no non-trivial infinitesimal Einstein deformations preserving the cone angles. This result can be interpreted as a higher-dimensional case of the celebrated Hodgson and Kerckhoff's theorem on deformations of hyperbolic $3$-cone-manifolds.\\
If all cone angles are smaller than $\pi$, we then give a construction which associates to any variation of the angles a corresponding infinitesimal Einstein deformation. We also show that these deformations are smooth on the singular locus.

\vskip 0.5\baselineskip
\noindent{\bf R\'esum\'e}

\vskip 0.5\baselineskip
\noindent Partant d'une c\^one-vari\'et\'e hyperbolique compacte de dimension $n\geq 3$, on \'etudie les d\'eformations de la m\'etrique dans le but d'obtenir des c\^ones-vari\'et\'es Einstein. Dans le cas o\`u le lieu singulier est une sous-vari\'et\'e ferm\'ee de codimension $2$ et que tous les angles coniques sont plus petits que $2\pi$, on montre qu'il n'existe pas de d\'eformations Einstein infinit\'esimales non triviales pr\'eservant les angles coniques. Ce r\'esultat peut s'interpr\'eter comme une g\'en\'eralisation en dimension sup\'erieure du c\'el\`ebre th\'eor\`eme de Hodgson et Kerckhoff sur les d\'eformations des c\^ones-vari\'et\'es hyperboliques de dimension $3$.\\
Si tous les angles coniques sont inf\'erieurs \`a $\pi$, on donne aussi une construction qui \`a chaque variation donn\'ee des angles associe une d\'eformation Einstein infinit\'esimale correspondante. On montre ensuite que ces d\'eformations sont lisses sur le lieu singulier.
\end{minipage}

\normalsize

\vskip 2cm

L'\'etude des vari\'et\'es Einstein, un domaine de recherche actif depuis plusieurs dizaines 
d'an\-n\'ees, est r\'ecemment revenue au coeur de l'actualit\'e math\'ematique, gr\^ace notamment aux 
travaux de G. Perelman \cite{Perelman} sur la conjecture de g\'eom\'etrisation de Thurston 
via le flot de Ricci. Les exemples de vari\'et\'es admettant des m\'etriques Einstein 
sont plus en plus nombreux, mais restent souvent cantonn\'es \`a des familles bien particuli\`eres.
Ainsi, on conna\^{\i}t de nombreux exemples de vari\'et\'es Einstein \`a courbure n\'egative, mais 
tr\`es peu sont non homog\`enes.
Le probl\`eme de trouver de telles vari\'et\'es est rendu plus difficile par le fait 
qu'elles sont rigides dans le cas compact: on ne peut pas les 
d\'eformer pour obtenir d'autres vari\'et\'es Einstein. 
Ces r\'esultats de rigidit\'e sont connus depuis longtemps pour les vari\'et\'es hyperboliques et pour 
les espaces sym\'etriques en g\'en\'eral \cite{Mostow}. Mais la situation n'est plus la m\^eme d\`es que 
l'on quitte les vari\'et\'es ferm\'ees~: la rigidit\'e est alors en g\'en\'eral subordonn\'ee \`a d'autres 
param\`etres, comme par exemple la structure conforme du bord \`a l'infini pour les vari\'et\'es 
hyperboliques.

Dans leur c\'el\`ebre article \cite{HK}, Hodgson et Kerckhoff montrent que contrairement au 
cas compact, il est possible de d\'eformer des vari\'et\'es hyperboliques \`a singularit\'es coniques.
Plus pr\'ecis\'ement, pour
une large classe de c\^ones-vari\'et\'es hyperboliques de dimension $3$, l'espace
des structures coniques hyperboliques au voisinage d'une c\^one-vari\'et\'e donn\'ee
est param\'etr\'e par les angles coniques. Si l'on se donne une petite variation
des angles, il existe donc une unique structure de c\^ones-vari\'et\'es hyperboliques proche de la 
structure de d\'epart et r\'ealisant la variation donn\'ee des angles coniques.
Leur r\'esultat principal est le th\'eor\`eme
de rigidit\'e infinit\'esimale suivant : si $M$ est une c\^one-vari\'et\'e hyperbolique
de dimension $3$ de volume fini, dont le lieu singulier forme un entrelacs et
dont tous les angles coniques sont inf\'erieurs \`a $2\pi$, alors il est
impossible de la d\'eformer sans modifier ses angles.
Cet article, compl\'et\'e par des travaux plus r\'ecents
(voir notamment \cite{HK2}, \cite{Kojima} et \cite{Weiss}), a \'et\'e le point de
d\'epart de nombreux d\'eveloppements dans l'\'etude de la g\'eom\'etrie des vari\'et\'es
hyperboliques de dimension $3$, tels que la g\'eom\'etrisation des petites
orbifolds ou l'\'etude des groupes kleiniens (\cite{BLP}, \cite{BB}).

Le principe de la d\'emonstration du th\'eor\`eme de rigidit\'e infinit\'esimale de
Hodgson et Kerckhoff est de
r\'eussir \`a appliquer la m\'ethode de Calabi-Weil (cf \cite{Calabi}, \cite{Garland}, \cite{Weil}) aux
c\^ones-vari\'et\'es : on montre que la repr\'esentation d'holonomie
n'admet pas de d\'eformations non triviales de la forme voulue. Cela n\'ecessite
d'\'etablir des formules d'int\'egration par parties ainsi qu'un r\'esultat du type
th\'eor\`eme de Hodge. Ce genre de difficult\'es est inh\'erent \`a l'\'etude
des c\^ones-vari\'et\'es; on verra dans cet article comment les aborder.

Dans le cas des vari\'et\'es ferm\'ees, Koiso \cite{Koiso} a donn\'e un analogue de la
m\'ethode de Calabi-Weil, qui n'utilise plus la repr\'esentation d'holonomie mais
\'etudie directement les d\'eformations de la m\'etrique (cf aussi \cite{Besse}, \S
12.H). Cette deuxi\`eme m\'ethode pr\'esente l'avantage de s'appliquer, en dimension sup\'erieure, \`a une 
classe de vari\'et\'es plus vaste, et en particulier aux vari\'et\'es Einstein (v\'erifiant de bonnes 
conditions de courbure).

Il est int\'eressant de regarder si ces techniques s'appliquent aux vari\'et\'es \`a singularit\'es 
coniques, et permettent d'obtenir une g\'en\'eralisation du th\'eor\`eme de Hodgson et Kerckhoff.
Il devrait \^etre alors
possible de construire, \`a partir d'une vari\'et\'e hyperbolique \`a cusps (que l'on peut 
consid\'erer comme une c\^one-vari\'et\'e d'angles coniques nuls), des c\^ones-vari\'et\'es Einstein proches, 
dont les angles coniques (suffisamment petits) sont donn\'es. On peut choisir 
ces angles de la forme $2\pi /n$ ; en prenant ensuite un rev\^etement appropri\'e, on obtient une 
vari\'et\'e compacte non singuli\`ere, dont la m\'etrique a priori non homog\`ene est Einstein, \`a courbure 
sectionnelle n\'egative. Comme il a \'et\'e mentionn\'e, on conna\^{\i}t actuellement tr\`es peu d'exemples de telles vari\'et\'es riemanniennes ; la construction ci-dessus en donnerait toute une famille.

Le but de cet article est d'utiliser ces techniques pour montrer qu'infinit\'esimalement, la 
situation en dimension sup\'erieure \`a $3$ est la m\^eme qu'en dimension $3$.
On adapte pour cela la m\'ethode de Koiso, ce qui ne se fait pas sans difficult\'e. En effet, pour \'eliminer  les d\'eformations g\'eom\'etriquement triviales, il est indispensable de normaliser les d\'eformations infinit\'esimales (l'analogue dans le cadre des repr\'esentations d'holonomie consiste \`a trouver un repr\'esentant harmonique d'une classe de cohomologie donn\'ee). Le probl\`eme, d\^u au caract\`ere singulier de la m\'etrique, est que la solution n'est pas unique et pr\'esente en g\'en\'eral une perte de r\'egularit\'e. Il est alors n\'ecessaire d'\'etudier plus en d\'etail l'\'equation aux d\'eriv\'ees partielles intervenant dans la normalisation et le laplacien associ\'e, qui est un ``op\'erateur d'ar\^ete'' elliptique selon la terminologie de \cite{Mazzeo}. On observe que le comportement des solutions de l'\'equation de normalisation est explicitement contr\^ol\'e par la valeur des angles coniques, ce qui permet de donner des bons domaines de r\'esolution quand les angles sont suffisamment petits.
On peut alors d\'emontrer que,
sous des hypoth\`eses voisines de celles du th\'eor\`eme de Hodgson et Kerckhoff, il est impossible de  
d\'eformer une c\^one-vari\'et\'e hyperbolique en des c\^ones-vari\'et\'es Einstein sans en modifier les
angles coniques. En particulier, on red\'emontre dans le cas de la dimension
trois le th\'eor\`eme de rigidit\'e infinit\'esimale ci-dessus.
On donne aussi une construction qui, \`a toute variation donn\'ee du $p$-uplet des angles coniques, associe une d\'eformation Einstein 
infinit\'esimale r\'ealisant cette variation au premier ordre.

Ces deux r\'esultats montrent que dans un certain sens, l'espace tangent \`a une c\^one-vari\'et\'e hyperbolique parmi les structures de c\^ones-vari\'et\'es Einstein est de dimension finie, param\'etr\'e par les variations du $p$-uplet des angles coniques. La question naturelle est alors de savoir s'il est possible d'int\'egrer les d\'eformations Einstein infinit\'esimales, ce qui est un probl\`eme plus difficile. Il existe cependant des raisons de penser que la r\'eponse est oui (au moins si les angles coniques sont suffisamment petits), qui viennent du fait que les d\'eformations Einstein infinit\'esimales modifiant les angles coniques ont un comportement relativement r\'egulier. En particulier, on montre dans la derni\`ere section de cet article que ces d\'eformations sont lisses sur le lieu singulier. Cette situation contraste fortement avec le cas des m\'etriques asymptotiquement hyperboliques, o\`u les d\'eformations sont en g\'en\'eral beaucoup moins r\'eguli\`eres sur le bord qu'\`a l'int\'erieur.

\section{Pr\'esentation des r\'esultats}

Les principaux r\'esultats de cet article sont les deux th\'eor\`emes suivants :

\begin{nThm}[\ref{thmppal}]
Soit $M$ une c\^one-vari\'et\'e hyperbolique compacte, dont le lieu singulier
  forme une sous-vari\'et\'e ferm\'ee de codimension $2$, et dont tous les
  angles coniques sont strictement inf\'erieurs \`a $2\pi$.
Alors toute d\'eformation Einstein infinit\'esimale ne
  modifiant pas les angles coniques est triviale.
\end{nThm}

\medskip

\begin{nThm}[\ref{consdef}]
Soit $M$ une c\^one-vari\'et\'e hyperbolique compacte, dont le lieu singulier
  forme une sous-vari\'et\'e ferm\'ee de codimension $2$, et dont les angles coniques $\alpha_1,\ldots \alpha_p$ 
sont tous strictement inf\'erieurs \`a $\pi$. Soit $\dot{\alpha} = (\dot{\alpha_1}, \ldots 
\dot{\alpha_p})$ une variation donn\'ee du $p$-uplet des angles coniques. Alors il existe une 
d\'eformation Einstein infinit\'esimale, unique \`a d\'eformartions triviales pr\`es, 
induisant la variation des angles coniques donn\'ee.
\end{nThm} 

\bigskip

Les restrictions impos\'ees \`a la g\'eom\'etrie des c\^ones-vari\'et\'es sont
essentiellement les m\^emes que dans l'article de Hodgson et
Kerckhoff \cite{HK} (on aurait pu remplacer l'hypoth\`ese ``$M$ compacte''
par l'hypoth\`ese ``$M$ de volume fini'', mais les choses sont quand m\^eme plus
simples dans le cas compact). La condition sur la g\'eom\'etrie du
lieu singulier est plus cruciale : c'est elle qui permet d'avoir
un bon mod\`ele local pour faire les calculs, car de
mani\`ere g\'en\'erale, le lieu singulier d'une c\^one-vari\'et\'e peut \^etre
beaucoup plus compliqu\'e. Enfin, la condition sur les angles
coniques est une hypoth\`ese technique qui para\^{\i}t de prime abord
assez myst\'erieuse. En fait, on verra dans la section \ref{vois3}
que les angles coniques r\'egissent en partie la croissance au voisinage du lieu
singulier des solutions d'un laplacien ; plus les angles sont petits, plus on
contr\^ole ces solutions.

L'outil principal dans la d\'emonstration de la rigidit\'e infinit\'esimale est
connu sous le nom
de {\em technique de Bochner}. En partant d'une \'equation du type $Pu=0$ o\`u $P$
est un op\'erateur diff\'erentiel du deuxi\`eme ordre de type laplacien, on
exprime $P$ comme somme d'un op\'erateur auto-adjoint positif $Q^*Q$
de degr\'e $2$ et d'un op\'erateur $R$ de degr\'e $0$ faisant intervenir la
courbure. Une telle d\'ecomposition
$$P = Q^*Q + R $$
 s'appelle une {\em formule de Weitzenb\"ock}; on en rencontrera \`a de nombreuses
reprises dans cet article. Ensuite, si elle est valide, une int\'egration par parties donne
$$0 = \<Pu,u\> = ||Qu||^2 + \<Ru,u\>.$$
Si l'op\'erateur $R$ est tel que $\<Ru,u\> \geq c ||u||^2$ avec $c>0$, on
trouve alors $u=0$. Le lecteur int\'eress\'e par le sujet pourra se r\'ef\'erer \`a
\cite{Besse}, \S 1.I. Et si on se place sur les bons espaces, le fait d'avoir montr\'e par cette technique l'injectivit\'e de l'op\'erateur $P$ suffit pour obtenir son inversibilit\'e, ce qui permettra de construire les d\'eformations Einstein infinit\'esimales. Mais avant d'en arriver l\`a, il faut d'abord montrer qu'il existe des formules de Stokes pour les int\'egrations par partie propos\'ees, et il faut aussi v\'erifier que les objets que l'on consid\`ere rentrent dans le cadre de ces formules.

\bigskip

Dans les premi\`eres parties de cet article, on met en place le cadre et les diff\'erents outils utilis\'es par la suite.
En particulier, on donne la d\'efinition pr\'ecise des c\^ones-vari\'et\'es 
envisag\'ees ainsi que les restrictions impos\'es \`a leur topologie, 
d'o\`u l'on d\'eduit que les
d\'eformations infinit\'esimales d'une telle structure peuvent
toujours se mettre sous une forme standard (i.e. appartenant \`a une
certaine famille de d\'eformations) au voisinage du lieu singulier. 
En particulier, une d\'eformation ne
modifiant pas les angles coniques a la propri\'et\'e d'\^etre $L^2$, \`a d\'eriv\'ee
covariante $L^2$~; c'est entre autres pour cette raison que l'on sera amen\'e
ensuite \`a travailler principalement dans le cadre $L^2$. 
On rappelle aussi la d\'efinition
des m\'etriques Einstein, des d\'eformations infinit\'esimales Einstein et des d\'eformations triviales. 

La section suivante commence par quelques r\'esultats de la th\'eorie des op\'erateurs non born\'es
d'un espace de Hilbert, qui nous seront utiles pour r\'esoudre les \'equations aux d\'eriv\'ees partielles apparaissant naturellement dans ce type de probl\`eme d'analyse g\'eom\'etrique.
On passe ensuite \`a la th\'eorie de Hodge $L^2$ et aux diff\'erentes formules de Stokes dont on aura besoin pour faire fonctionner des techniques de Bochner.
On d\'emontre en particulier le th\'eor\`eme suivant :

\begin{nThm}[\ref{ipp}] Soient $u \in C^\infty(T^{(r,s)}M)$, $v \in
  C^\infty(T^{(r+1,s)}M)$ tels que $u$, $\n u$, $v$, $\na v$ soient dans
  $L^2$. Alors $\<u,\na v\> = \<\n u,v\>.$
\end{nThm}

Les r\'esultats de cette section \ref{secIPP} seront d'usage constant dans la suite de l'article. Ici encore, on verra qu'il est naturel de travailler avec des objets appartenant \`a des
espaces $L^2$. En plus des th\'eor\`emes d'int\'egrations par partie,
on donne leur interpr\'etation en termes d'op\'erateurs non born\'es ainsi que d'autres r\'esultats 
utilisant les m\^emes techniques.

\bigskip

Pour \'eliminer les d\'eformations triviales, on cherche dans la section \ref{Bianchi} \`a imposer la condition de jauge de Bianchi. Un premier r\'esultat dans ce sens est le suivant :

\begin{nThm}[\ref{normgen}] Soit $M$ une c\^one-vari\'et\'e Einstein \`a courbure de Ricci n\'egative. On a 
la d\'e\-com\-po\-si\-tion en somme directe
$$L^2(S^2M) = \ker \beta_{max} \oplus \im \delta^*_{min}$$
\end{nThm}

Ce r\'esultat signifie que toute d\'eformation infinit\'esimale $L^2$ peut \^etre normalis\'ee selon la jauge de Bianchi, ce qui justifie a posteriori son choix. Cependant pour pouvoir adapter la m\'ethode de Koiso on a besoin de r\'esultats plus forts, garantissant une certaine r\'egularit\'e \`a la d\'eformation normalis\'ee.
On est alors amen\'e \`a \'etudier de plus pr\`es l'\'equation de normalisation et
l'op\'erateur correspondant
$$L = \na \n + (n-1)Id = \Delta + 2(n-1)Id$$
agissant sur les 1-formes. Il est facile de trouver des solutions dans l'espace de Sobolev $L^{1,2}$, et le but est de comprendre le comportement de ces solutions. Pour ce faire, et apr\`es
avoir pr\'ealablement exhib\'e une d\'ecomposition adapt\'ee en s\'eries de Fourier
g\'en\'eralis\'ees (\S \ref{vois2}), on \'etudie les
solutions de l'\'equation homog\`ene au voisinage de la singularit\'e. On
montre que leur comportement est \'etroitement li\'e aux angles coniques; par
exemple, la norme ponctuelle d'une solution donn\'ee au voisinage d'une
composante connexe du lieu singulier d'angle conique $\alpha$ est en $r^k$
avec $k \in \{\pm 1 \pm 2p\pi\alpha^{-1}, \pm 2p\pi\alpha^{-1}\ |\ p \in
\mathbb{Z} \}$. Les restrictions impos\'ees sur les angles coniques permettent
de contr\^oler suffisamment les solutions de l'\'equation homog\`ene, et finalement
les solutions de l'\'equation de normalisation tout court ; et plus on restreint les angles, plus on contr\^ole les solutions. On aboutit ainsi aux
th\'eor\`eme suivant :

\begin{nThm}[\ref{nablad}] Soit $M$ une c\^one-vari\'et\'e hyperbolique dont tous les
  angles coniques sont strictement inf\'erieurs \`a $2\pi$.
Soit $\phi$ une forme lisse appartenant \`a $L^2(T^*M)$. Alors il existe une
unique forme $\eta \in
C^\infty(T^*M)$ solution de l'\'equation $L\eta=\phi$ telle que $\eta$, $\n
\eta$, $d\delta \eta$, et $\n d\eta$ soient dans $L^2$.
\end{nThm}
 
\begin{nThm}[\ref{nnu}]
Soit $M$ une c\^one-vari\'et\'e hyperbolique dont tous les
angles coniques sont strictement inf\'erieurs \`a $\pi$. Alors l'op\'erateur 
$$L = \na \n + (n-1)Id\ :\ L^{2,2}(T^*M) \to L^2(T^*M)$$
est un isomorphisme.
\end{nThm}

Une fois ces r\'esultats \'etablis, il est relativement facile de faire fonctionner
la m\'ethode de Koiso pour d\'emontrer les th\'eor\`emes \ref{thmppal} et \ref{consdef} ; c'est
l'objet de la section \ref{demo}. Pour la rigidit\'e infinit\'esimale, le principe est le suivant. Partant d'une
d\'eformation infinit\'esimale Einstein $h_0$ pr\'eservant les angles
(donc \`a d\'eriv\'ee covariante $L^2$) d'une c\^one-vari\'et\'e hyperbolique,
dont tous les angles coniques sont inf\'erieurs \`a $2\pi$,
la d\'emonstration de sa trivialit\'e se fait en deux temps. On a
d'abord besoin de se d\'ebarrasser des d\'eformations triviales, on utilise donc
le r\'esultat \ref{nablad} mentionn\'e ci-dessus pour r\'esoudre l'\'equation de normalisation.
On applique ensuite une technique de Bochner \`a la d\'eformation normalis\'ee
$h=h_0-\delta^*\eta$. En utilisant la formule de Weitzenb\"ock idoine et le
premier r\'esultat d'int\'egration par parties, on obtient
$${\delta^\n d^\n h + (n-2) h =0}.$$
Une deuxi\`eme int\'egration par parties, un peu plus compliqu\'ee,
permet de conclure que $h_0 = \delta^*\eta$, et donc que l'on a bien rigidit\'e
infinit\'esimale relativement aux angles coniques au sein des c\^ones-vari\'et\'es
Einstein.
 
Pour le th\'eor\`eme \ref{consdef}, la m\'ethode de construction 
est la suivante : on part d'une d\'eformation 
infinit\'esimale $h_0$, Einstein au voisinage du lieu singulier, et induisant la variation voulue 
des angles. On cherche 
alors \`a lui rajouter une d\'eformation $L^{1,2}$ (donc ne modifiant pas les angles) de telle sorte 
que la somme v\'erifie l'\'equation d'Einstein lin\'earis\'ee. Cela revient \`a r\'esoudre une \'equation de 
la forme $$(\na \n -2 \RO) h = \phi,$$
o\`u $\phi=E'(h_0)$ est un $2$-tenseur v\'erifiant la condition de jauge de Bianchi, et \`a 
s'assurer que la solution trouv\'ee v\'erifie aussi cette condition, ce qui se fait assez ais\'ement \`a partir des th\'eor\`emes \ref{normgen} et \ref{nnu}. La d\'eformation $h-h_0$ est 
alors Einstein et a les propri\'et\'es voulues.

\bigskip

La derni\`ere section de cet article est consacr\'e aux d\'eformations Einstein infinit\'esimales donn\'ees par le th\'eor\`eme \ref{consdef}, et en particulier \`a la r\'egularit\'e des d\'eformations induites de la m\'etrique du lieu singulier. Il est pour cela n\'ecessaire d'\'etudier en d\'etails l'op\'erateur $\na \n -2 \RO$ agissant sur les $2$-tenseurs sym\'etriques. La m\'ethode est en grande partie la m\^eme que pour l'\'etude de l'op\'erateur $\na \n + (n-1)Id$, r\'ealis\'ee dans la section \ref{secLap}. Ici encore, le comportement des solutions de l'\'equation homog\`ene est directement li\'e \`a la valeur des angles coniques. On aboutit alors au th\'eor\`eme suivant, conclusion de cet article :

\begin{nThm}[\ref{regdef}] Soit $M$ une c\^one-vari\'et\'e hyperbolique dont tous les angles coniques $\alpha_1,\ldots \alpha_p$ sont strictement inf\'erieurs \`a $\pi$. Soit $\dot{\alpha} = (\dot{\alpha_1}, \ldots 
\dot{\alpha_p})$ une variation donn\'ee du $p$-uplet des angles coniques, et soit $h_{\dot{\alpha}}$ la 
d\'eformation Einstein infinit\'esimale normalis\'ee correspondante. Alors la d\'eformation infinit\'esimale $h_\Sigma$ de la m\'etrique du lieu singulier, induite par $h_{\dot{\alpha}}$, est $C^\infty$.
\end{nThm}

\section{Pr\'eliminaires}

\subsection{C\^ones-vari\'et\'es}\label{defnot}

On va maintenant pr\'eciser le cadre dans lequel on se place.
La notion de c\^one-vari\'et\'e, plus g\'en\'erale que celle d'orbifold, a
\'et\'e introduite par Thurston \cite{ThurstonGeom} pour l'\'etude des
d\'eformations des vari\'et\'es hyperboliques \`a cusps en dimension $3$.
Le cas le plus fr\'equemment rencontr\'e est celui des
c\^ones-vari\'et\'es \`a courbure constante. Celles-ci sont relativement
simples \`a d\'efinir, soit g\'eom\'etriquement comme un recollement (local) de
simplexes g\'eod\'esiques, soit en explicitant la m\'etrique en
coordonn\'ees; c'est cette derni\`ere approche qui sera utilis\'ee ici.
Le lecteur int\'eress\'e pourra se reporter \`a \cite{ThurstonShapes}
pour une d\'efinition par r\'ecurrence des c\^ones-vari\'et\'es model\'ees sur
une g\'eom\'etrie.

La g\'eom\'etrie du lieu singulier d'une c\^one-vari\'et\'e arbitraire peut \^etre tr\`es
compliqu\'ee. Dans le cadre de notre \'etude on se limitera au cas o\`u il
forme une sous-vari\'et\'e de codimension deux, ce qui permet de parler d'angle
conique le long de chaque composante connexe du lieu singulier et d'avoir des
bons mod\`eles locaux pour mener \`a bien les calculs.

Enfin, comme notre but est de s'int\'eresser \`a des vari\'et\'es Einstein, on
s'autorise une classe assez large de m\'etriques \`a singularit\'es : on demande
juste que la m\'etrique conique ressemble asymptotiquement au produit de la
m\'etrique du lieu singulier avec la m\'etrique d'un c\^one (de dimension deux).

\bigskip

Soit $M$ une vari\'et\'e compacte de dimension $n\geq 3$, et $\Sigma =
\coprod_{i=1}^p \Sigma_i$ une sous-vari\'et\'e ferm\'ee plong\'ee de
codimension $2$, dont les $\Sigma_i$ sont les composantes
connexes. Dans la suite de ce texte on emploiera souvent la
notation $M$ pour d\'esigner improprement $M \setminus \Sigma$.

\begin{Def}\label{defcv} Soient $\alpha_1,\ldots,\alpha_p$ des r\'eels
  positifs. La vari\'et\'e $M$ est munie d'une structure de {\em
  c\^{o}ne-vari\'et\'e}, de lieu singulier $\Sigma = \coprod_{i=1}^p \Sigma_i$ et
  d'angles coniques les $\alpha_i$, si :
\begin{itemize}

\item $M\setminus \Sigma$ est munie d'une m\'etrique riemannienne $g$, non
  compl\`ete;

\item pour tout $i$, $\Sigma_i$ est munie d'une m\'etrique riemannienne $g_i$;

\item pour tout $i$, tout point $x$ de $\Sigma_i$ a un voisinage $V$ dans
  $M$ diff\'eomorphe \`a $D^2\times U$, avec $U=V\cap \Sigma_i$ un voisinage de
  $x$ dans $\Sigma_i$,
dans lequel $g$ s'exprime en coordonn\'ees cylindriques locales sous
la forme $$g=dr^2 + (\frac{\alpha_i}{2\pi})^2 r^2d\theta^2+ g_i +
q,$$ o\`u $q$ est un 2-tenseur sym\'etrique v\'erifiant $g(q,q)=o(r^2)$
et $g(\n q,\n q)=o(1)$.
\end{itemize}
\end{Def}

\medskip

Dans la suite on exprimera souvent la m\'etrique $g$ sous la forme
l\'eg\`erement diff\'erente $$g = dr^2 + r^2d\theta^2+ g_i + q,$$
 o\`u la coordonn\'ee d'angle $\theta$ est d\'efinie non plus modulo
 $2\pi$ mais modulo l'angle conique $\alpha_i$.

Une c\^one-vari\'et\'e hyperbolique est alors une c\^one-vari\'et\'e telle
que les m\'etriques $g$ et $g_i$ sont hyperboliques. On a dans ce
cas, en reprenant les notations de la d\'efinition, $$q =
(\frac{\alpha_i}{2\pi})^2(\sinh(r)^2 - r^2) d\theta^2 +
(\cosh(r)^2 - 1) g_i.$$

Pour d\'emontrer un certain nombre de r\'esultats, on aura besoin d'un contr\^ole sur les angles
coniques ; par exemple, la preuve de la rigidit\'e infinit\'esimale \`a partir de la fin de la section 
\ref{secLap} n\'ecessite que tous les angles soient inf\'erieurs \`a $2\pi$. Ces conditions seront 
explicit\'ees quand elles appara\^itront.

\bigskip

Le caract\`ere singulier des c\^ones-vari\'et\'es pose probl\`eme pour adapter
 certaines m\'ethodes d'analyse, comme une technique de Bochner. Il faut
toujours v\'erifier si les choses marchent de la m\^eme mani\`ere que dans le cas
compact.

La premi\`ere difficult\'e va venir des int\'egrations par parties. Premi\`erement,
pour garantir que les expressions manipul\'ees ont un sens, on sera oblig\'e
de travailler avec des objets $L^2$. Deuxi\`emement, il va falloir d\'emontrer
qu'on peut effectivement appliquer des formules de type Stokes : ce sera
l'objet de la partie \ref{secIPP}. Au final on sera en mesure d'effectuer
des int\'egrations par parties pour les op\'erateurs $d$ et $\delta$, et $\n$ et
$\na$. Mais un tel r\'esultat n'existe pas (\`a la connaissance de l'auteur) pour les
op\'erateurs $d^\n$ et $\delta^\n$; on devra donc contourner cette
difficult\'e quand on en aura besoin (section \ref{rigidite}).

La plus grande difficult\'e va venir de l'\'equation correspondant \`a l'op\'erateur d'Einstein 
lin\'earis\'e et de l'\'equation de normalisation,
\'etudi\'ees dans les sections \ref{secLap} et \ref{secLap2}. Bien qu'en pr\'esence de
sympathiques op\'erateurs elliptiques de la forme $\na \n$ plus un terme born\'e d'ordre $0$, on ne 
peut pas appliquer la th\'eorie classique sur une c\^one-vari\'et\'e, dont la m\'etrique est
singuli\`ere. Les \'equations admettront encore des solutions, mais
celles-ci ne seront plus uniques, et on aura des probl\`emes de 
perte de r\'egularit\'e. Cependant, en imposant que les angles
coniques soient assez petits, on arrivera \`a avoir suffisamment de contr\^ole 
sur les solutions et la norme de certaines combinaisons lin\'eaires de
leurs d\'eriv\'ees pour faire fonctionner les d\'emonstrations.

\subsection{D\'eformations des c\^ones-vari\'et\'es}

Soit $(M,g)$ une c\^one-vari\'et\'e au sens ci-dessus, de lieu singulier
$\Sigma$. Soit maintenant $g_t$ une famille de m\'etriques singuli\`eres,
d\'erivable, telle que $g_0=g$ et que pour tout $t$,
$(M,g_t)$ soit une c\^one-vari\'et\'e de lieu singulier $\Sigma$.

Si $x$ est un point de $\Sigma$, pour tout $t$ il existe par d\'efinition un
voisinage de $x$ dans $M$ dans lequel on a l'expression ci-dessus pour la
m\'etrique en coordonn\'ees cylindriques. Quitte \`a les restreindre, ces voisinages
sont tous diff\'eomorphes, et on peut donc faire agir une famille $\phi_t$ de
diff\'eomorphismes de telles fa\c{c}ons que les coordonn\'ees cylindriques locales
pour l'expression de $\phi_t^*g_t$ soient les m\^emes pour tout $t$.

Dit d'une autre mani\`ere, il existe un voisinage $V$ de $x$ dans $M$,
diff\'eomorphe \`a $D^2 \times U$ o\`u $U = V \cap \Sigma$ est un voisinage de $x$
dans $\Sigma$, dans lequel on peut trouver des coordonn\'ees cylindriques telles
que pour tout $t$, on ait  :
$$\phi_t^*g_t = dr^2 + (\frac{\alpha_t}{2\pi})^2r^2d\theta^2+ h_t + q_t.$$
Dans cette expression, $h_t$ d\'esigne une m\'etrique sur $U$ et $q_t$ est un
2-tenseur sym\'etrique qui v\'erifie les conditions de la d\'efinition \ref{defcv}.

Finalement, quitte \`a modifier la famille $g_t$ par des diff\'eomorphismes, ce
qui revient \`a modifier la d\'eformation infinit\'esimale par une d\'eformation
g\'eom\'etriquement triviale, on peut montrer que $h = \frac{dg_t}{dt}|_{t=0}$ est
au voisinage du lieu singulier combinaison lin\'eaire des quatre types de
d\'eformations suivants, modifiant :
\begin{list}{-}{}
\item l'angle,
\item la m\'etrique du lieu singulier,
\item le reste,
\item et enfin, la fa\c{c}on de ``recoller'' la variable d'angle quand on passe
  d'un syst\`eme de coordonn\'ees \`a un autre.
\end{list}
Au voisinage du lieu singulier, une d\'eformation $h$ modifiant le reste v\'erifie 
$|h| = o(r)$ et $|\n h| = o(1)$. Les autres d\'eformations sont de la forme (\`a une 
d\'eformation modifiant le reste pr\`es) $\lambda r^2 d\theta^2$ pour celle modifiant l'angle,
$h_\Sigma$, extension d'un $2$-tenseur sym\'etrique d\'efini sur $\Sigma$ pour celle modifiant la 
m\'etrique du lieu singulier, et $r^2 d\theta.\omega$, o\`u $\omega$ est l'extension d'une 
$1$-forme d\'efinie sur $\Sigma$, pour celle modifiant la variable d'angle.

Il est important de noter que les toutes ces d\'eformations infinit\'esimales sont
$L^2$. Par contre seules les trois derni\`eres ont leur d\'eriv\'ee covariante dans
$L^2$. En effet, on peut constater que sur les expressions donn\'ee ci-dessus, $|\n h|$ est
en $o(1)$, $|\n h_\Sigma|$ et $|\n r^2 d\theta.\omega|$ sont en $O(1)$, alors que $|\n r^2 
d\theta^2|$ est en $r^{-1}$, et n'est donc pas $L^2$.
Ainsi, c'est au niveau du caract\`ere $L^2$ ou non de la d\'eriv\'ee
covariante de la d\'eformation que l'on voit si celle-ci pr\'eserve ou non les
angles coniques.

\subsection{D\'eformations Einstein infinit\'esimale et \'equa\-tion
de normalisation}\label{secEin}

Par d\'efinition, une {\em m\'etrique Einstein} est une m\'etrique
riemannienne $g$ v\'erifiant l'\'equation $$ric(g) = c g,$$ o\`u le
terme de gauche est le tenseur de courbure de Ricci et o\`u $c$ est
une constante. Notons que si on remplace $g$ par $\lambda g$, avec $\lambda$
une constante strictement positive, alors la nouvelle m\'etrique v\'erifie
l'\'equation ci-dessus en rempla\c{c}ant $c$ par $\lambda^{-1}c$;
donc en fait c'est principalement le signe et non la valeur exacte
de la constante $c$ qui compte. On peut ainsi distinguer trois
grandes classes de m\'etriques Einstein suivant que $c$ est n\'egatif,
positif ou nul.

Les m\'etriques \`a courbure sectionnelle constante sont toujours
Einstein; en dimension $3$ ce sont les seules. Par contre d\`es la
dimension $4$ il y a beaucoup plus de m\'etriques Einstein que de
m\'etriques \`a courbure sectionnelle constante; on peut donc
consid\'erer la condition Einstein comme un affaiblissement ou une
g\'en\'eralisation de la condition de courbure sectionnelle constante.

Puisque l'on s'int\'eresse principalement aux c\^ones-vari\'et\'es
hyperboliques, on ne consid\`erera que des m\'etriques Einstein
v\'erifiant $E(g)=0$, avec $$E(g) = ric(g) + (n-1) g.$$ La constante
$(n-1)$ est choisie de telle sorte que les m\'etriques hyperboliques
v\'erifient cette \'equation.

\bigskip

Soit $g_t$ une famille lisse de m\'etriques Einstein (c'est-\`a-dire
v\'erifiant $E(g_t)= 0$) sur une vari\'et\'e donn\'ee $M$, avec $g_0=g$.
Le $2$-tenseur sym\'etrique $h=\frac{d}{dt}g_t|_{t=0}$ v\'erifie alors
l'\'equation d'Einstein lin\'earis\'ee $$E'_g(h)=0.$$ Le calcul de
$E'_g$ est classique, voir par exemple \cite{Besse} \S 1.K :
\begin{eqnarray}E'_g(h) = \na_g\n_g h - 2 \mathring{R}_gh - \delta^*_g(2\delta_g
h + d\tr_g h).\label{0019}\end{eqnarray}

Les op\'erateurs utilis\'es ici n\'ecessitent un peu d'explication. La
notation $\n_g$, ou $\n$ pour simplifier, d\'esigne la d\'eriv\'ee covariante ou
connexion de
Levi-Civit\`a associ\'ee \`a la m\'etrique riemannienne $g$. Elle admet un
adjoint formel not\'e $\na_g$ : si $(e_i)_{i=1\ldots n}$ est une base
orthonorm\'ee, on a $$\na_g \eta (X_1,\ldots,X_p) = - \sum_{i=1}^n
(\n_{e_i} \eta)(e_i,X_1,\ldots,X_p).$$

Pour les tenseurs sym\'etriques, on d\'efinit $\delta^*_g :
\mathcal{S}^p M \to \mathcal{S}^{p+1} M$ comme \'etant la compos\'ee
de la d\'eriv\'ee covariante et de la sym\'etrisation. En particulier,
si $\eta \in \Omega^1 M = \mathcal{S}^1 M$, alors
\begin{eqnarray*}
\delta^*_g \eta (x,y) & = &\frac{1}{2}((\n_x \eta)(y) + (\n_y
\eta)(x)\\
& = &\frac{1}{2}(g(\n_x
\eta^\sharp,y) + g(\n_y \eta^\sharp,x))\\
& = &\frac{1}{2}L_{\eta^\sharp}g(x,y),
\end{eqnarray*}
o\`u $L_{\eta^\sharp}$ d\'esigne la d\'eriv\'ee de Lie le long du champ
de vecteur $\eta^\sharp$ dual (pour la m\'etrique $g$) \`a la forme
$\eta$. L'adjoint
formel de l'op\'erateur $\delta^*_g$ se note $\delta_g$ ; c'est juste la
restriction de $\na_g$ \`a $S^{p+1}M$.

Ensuite, $\mathring{R}_g$ d\'esigne l'action du tenseur de courbure
$R_g$ sur les $2$-tenseurs sym\'etriques : si $h$ est une section de
$S^2M$, on pose
$$\mathring{R}_g h(x,y) = \sum_{i=1}^n h(R_g(x,e_i)y,e_i),$$
o\`u $(e_i)$ est une base orthonormale pour $TM$;
c'est encore un $2$-tenseur sym\'etrique. Si $g$ est hyperbolique, on a alors
\begin{eqnarray} \mathring{R}_g h = h - (\tr_g h)g. \label{0020} \end{eqnarray}
L'op\'erateur $\mathring{R}_g$ appara\^it souvent dans les 
probl\`emes de d\'eformations de m\'etriques ; la propri\'et\'e suivante (\cite{Besse}, \S 1.132) met 
l'accent sur son lien avec les m\'etriques Einstein.

\begin{Prop}
Une m\'etrique riemannienne $g$ est Einstein si et seulement si l'op\'erateur $\mathring{R}_g$ 
envoie l'espace $S^2_0$ des $2$-tenseurs sym\'etriques sans trace dans lui-m\^eme.
\end{Prop}

Enfin, la notation $\tr_g$ d\'esigne juste la trace par rapport \`a
$g$ : si $h$ est un $2$-tenseur,
$$\tr_g h = \sum_{i=1}^n h(e_i,e_i).$$
Dans la suite et pour all\'eger les notations, on
omettra le plus fr\'equemment l'indice $g$.

\bigskip

Par d\'efinition, une {\em d\'eformation Einstein infinit\'esimale} de
la vari\'et\'e Einstein $(M,g)$ est un $2$-tenseur sym\'etrique $h$
v\'erifiant l'\'equation $E'_g(h)=0$.

Maintenant, si $g$ est Einstein et si $\phi$ est un
diff\'eomorphisme de $M$, alors la m\'etrique tir\'ee en arri\`ere $\phi^*
g$ est aussi Einstein. Par cons\'equent, si $\phi_t$ est une famille
lisse de diff\'eomorphismes telle que $\phi_0$ soit l'identit\'e,
alors la d\'eformation infinit\'esimale associ\'ee
$\frac{d}{dt}\phi_t^*g|_{t=0}$ est naturellement Einstein. Une
telle d\'eformation est qualifi\'ee de {\em triviale}. Soit $X$ le
champ de vecteurs sur $M$ d\'efini par
$X(x)=\frac{d}{dt}(\phi_t(x))|_{t=0}$, et soit $\eta = X^\flat$
la $1$-forme duale, c'est-\`a-dire v\'erifiant $\eta(Y) = g(X,Y)$
pour tout vecteur $Y$. On a les relations $$\frac{d}{dt}(\phi_t^*
g)|_{t=0} = L_X g = 2 \delta^*_g \eta;$$ l'espace des
d\'eformations infinit\'esimales triviales est donc \'egal \`a $\im
\delta^*_g$.

\bigskip

La fa\c{c}on habituelle de se d\'ebarrasser des d\'eformations triviales
est d'imposer une condition de jauge, c'est-\`a-dire de ne
consid\'erer que des d\'eformations infinit\'esimales v\'erifiant une
certaine \'equation. On en trouve plusieurs dans la litt\'erature, on
utilisera ici la jauge de Bianchi (voir \cite{Biquard} \S I.1.C,
\cite{Anderson} \S 2.3). On
veut donc que nos d\'eformations infinit\'esimales $h$ v\'erifient
$$\beta_g(h)=0,$$ o\`u $\beta_g : \mathcal{S}^2M \to \Omega^1M$ est
l'op\'erateur de Bianchi (associ\'e \`a la m\'etrique $g$) d\'efini par
$$\beta_g(h)=\delta_g h +\frac{1}{2}d\tr_g h.$$

D'un point de vue g\'eom\'etrique, d'autres conditions de normalisation sont plus naturelle ;
par exemple la condition $\delta_g h=0$, qui correspond \`a regarder des d\'eformations 
$L^2$-orthogonales aux d\'eformations triviales, ou la condition d'\^etre infinit\'esimalement 
harmonique, voir \cite{Besse} \S 12.C, cf aussi \cite{HK2}. Mais en g\'en\'eral ces conditions 
co\"{\i}ncident sur les d\'eformations infinit\'esimales Einstein. La condition de jauge de Bianchi est 
plus naturelle d'un point de vue analytique, pour rendre les op\'erateurs elliptiques ; cf par 
exemple \cite{DeTurck}.

Ainsi, \'etant donn\'ee une d\'eformation infinit\'esimale $h_0$, on veut
pouvoir la modifier par une d\'e\-for\-ma\-tion triviale, de fa\c{c}on
essentiellement unique, de telle sorte que le r\'esultat v\'erifie la
condition de jauge. Dit plus pr\'ecis\'ement, on veut trouver une
$1$-forme $\eta$ telle que la d\'eformation normalis\'ee $h=h_0 -
\delta^*\eta$ satisfasse $\beta(h)=0$; de fa\c{c}on \'equivalente, on
cherche \`a r\'esoudre {\em l'\'equation de normalisation} (on omet les
indices) 
\begin{eqnarray}\beta\circ \delta^* \eta = \beta(h_0). \label{0021} \end{eqnarray}

L'\'etude de cette \'equation et de l'op\'erateur $\beta\circ\delta^*$ est l'objet
de la section \ref{secLap}. 

On se placera entre autre dans le cadre de la
th\'eorie des op\'erateurs non born\'es entre espace de Hilbert, dont les r\'esultats
principaux sont cit\'es dans la section suivante.

\section{Analyse $L^2$ sur les c\^ones-vari\'et\'es}\label{secIPP}

\subsection{Quelques rappels sur les op\'erateurs non
born\'es}\label{rapop}

On va maintenant annoncer un certain nombre de d\'efinitions et propri\'et\'es concernant
les op\'e\-ra\-teurs non born\'es; le lecteur int\'eress\'e pourra consulter \cite{Riesz},
chapitre 8, ou \cite{Rudin}, chapitre 13.

Soient $E$ et $F$ deux espaces de Hilbert. Un {\em op\'erateur non born\'e} est une
application lin\'eaire
$$A : D(A) \to F$$
o\`u $D(A)$ (le domaine de A) est un
sous-espace vectoriel de $E$. En particulier, toute application lin\'eaire
(continue ou non) de $E$ dans $F$ est un op\'erateur non born\'e.

Soit $A$ et $B$ deux op\'erateurs non born\'es. On dit que $B$ est un {\em
  prolongement} de $A$, not\'e $A\subset B$, si $D(A) \subset D(B)$ et
  $B_{|D(A)} = A$.

Un op\'erateur non born\'e $A$ est {\em ferm\'e} si son graphe $G(A)= \{(u,A(u)) |
u\in D(A) \}$ est ferm\'e dans $E\times F$, ce qui revient \`a dire que pour toute
suite $(u_n)$ de $D(A)$ telle que $u_n \to u \in E$ et $A(u_n) \to v \in F$,
on a $u \in D(A)$ et $v=A(u)$.

Si $A$ est \`a domaine dense dans $E$, on peut d\'efinir son {\em adjoint}
$A^* : D(A^*) \subset F \to E$ de la
fa\c{c}on suivante :
$$v \in D(A^*) \Longleftrightarrow \exists w\in
E {\rm \ tel\ que\ } \forall u\in D(A),\ \<u,w\>_E = \<A(u),v\>_F.$$
Comme $D(A)$ est
dense dans $E$, l'\'el\'ement $w$ (si il existe) est unique ; on pose
$w=A^*(v)$. Remarquons que l'adjoint d'un op\'erateur est toujours ferm\'e.
On a aussi la propri\'et\'e \'evidente (si les op\'erateurs consid\'er\'es sont \`a domaine
dense) $A \subset B \implies B^* \subset A^*$.

Pour d\'efinir $A^{**}$, il faut v\'erifier que $A^*$ est \`a domaine dense, ce qui
n'est pas toujours le cas. Mais on a la propri\'et\'e suivante (voir \cite{Riesz} \S
117):

\begin{Prop} Soit $A$ un op\'erateur non born\'e de $E$ dans $F$, \`a domaine
  dense. Alors $A^*$ est \`a domaine dense si et seulement si $A$ admet un
  prolongement ferm\'e. Dans ce cas, $A^{**}$ est le plus petit prolongement
  ferm\'e de $A$, i.e. si on a $A \subset B$ avec $B$ ferm\'e, alors $A^{**}
  \subset B$.
\end{Prop}

On remarque aussi que le graphe de $A^{**}$ n'est autre que l'adh\'erence dans
$E\times F$ du graphe de $A$. D'autre part, si $A$ est ferm\'e, on a $A^{**} =
A$. En particulier, d\`es que cela a un sens, on a toujours $A^{***} = A^*$
(notons au passage que l'on a bien $(A^*)^{**} = (A^{**})^*$).

Si $A$ est injectif, on peut d\'efinir son inverse $A^{-1}$ : son domaine n'est
autre que l'image de $A$.

Pour pouvoir d\'efinir la composition de deux op\'erateurs non born\'es $A : D(A)
\subset E \to F$ et  $B : D(B) \subset F \to G$, on pose, par d\'efinition,
$$D(B\circ A) = \left\{ x \in D(A) | A(x) \in D(B) \right\}.$$
 De m\^eme, la somme se d\'efinit naturellement sur le domaine
$$D(A+A') = D(A) \cap D(A').$$
Il se peut \'evidemment que ces domaines soient r\'eduits \`a l'\'el\'ement
nul. Cependant, on a le th\'eor\`eme relativement surprenant suivant
(\cite{Riesz}, §118, ou \cite{Rudin}, th\'eor\`eme 13.13), et son corollaire :

\begin{Thm}\label{1+*} Si l'op\'erateur non born\'e $A : E \to F$ est ferm\'e et de
domaine dense, alors les op\'erateurs
$$ B = (A^*\circ A + Id)^{-1},\ C=A\circ (A^*\circ A + Id)^{-1}$$
sont des applications lin\'eaires {\em continues} de $E$ dans $E$ et
de $E$ dans $F$; de plus $||B||\leq 1$, $||C||\leq 1$, et $B$ est
auto-adjointe positive.
\end{Thm}

\medskip

\begin{Cor}\label{1+*b} Si l'op\'erateur non born\'e $A : E \to F$ est ferm\'e et de
domaine dense, et si $B : E \to E$ est une application lin\'eaire continue et auto-adjointe, 
alors l'op\'erateur $A^* \circ A + B$ est auto-adjoint.
\end{Cor}

\bigskip

Maintenant, soit $M$ une vari\'et\'e riemannienne, et soient $E$ et $F$
deux fibr\'es vectoriels sur $M$, munis de m\'etriques riemanniennes
$(.,.)_E$ et $(.,.)_F.$ On note $C^\infty_0(E)$ (resp.
$C^\infty(E),$ resp. $L^2(E)$) l'espace des sections de E qui sont
$C^\infty$ \`a support compact (resp. $C^\infty$, resp. $L^2$); de
m\^eme pour $F$. La m\'etrique sur $E$ et la forme volume sur $M$ font
de $L^2(E)$ un espace de Hilbert (pour le produit scalaire
${\<f,g\>_E = \int_M (f,g)_E dvol_M}$) dont $C^\infty_0(E)$ est un
sous-espace dense; de m\^eme pour $F$.

Soit $A$ un op\'erateur diff\'erentiel agissant sur les sections de
$E$. On le consid\`ere comme un op\'erateur non born\'e de domaine les
sections $C^\infty$ \`a support compact, i.e. $$A : C^\infty_0(E)
\to C^\infty_0(F) \subset L^2(F),$$
 et on suppose que $A$ admet un
{\em adjoint formel} $A^t : C^\infty_0(F) \to C^\infty_0(E)$, i.e.
tel que $$\<Au,v\>_F = \<u,A^tv\>_E\ \forall u \in C^\infty_0(E)
{\rm \ et\ } \forall v \in C^\infty_0(F).$$
 On a clairement $A^t \subset A^*$ donc $A^*$ est \`a domaine dense.

On pose alors $A_{min} = A^{**}$, c'est, on l'a vu,
le plus petit prolongement ferm\'e de $A$. Le graphe de $A^{**}$ est
l'adh\'erence du graphe de $A$, donc (et on peut prendre \c{c}a comme d\'efinition)
$$u \in D(A_{min}) \Leftrightarrow  \exists (u_n) \in C^\infty_0(E) {\rm \
  telle\ que\ } \lim_{n\to \infty} u_n = u {\rm \ et\ que\ la\ suite\ } (Au_n)
  {\rm \ converge\ dans\ } L^2,$$
$A_{min} u$ est alors la valeur de cette limite.

On pose aussi $A_{max} = (A^t)^*$ ; comme $A^t\subset A^*$, on a
$A^{**}\subset A_{max}$ et donc $A \subset A_{max}$. De plus $A^t \subset
(A^t)^{**}=(A_{max})^*$, et, vu la propri\'et\'e de minimalit\'e de $^{**}$, on en
d\'eduit que $A_{max}$ est le plus grand prolongement
de $A$ dont l'adjoint prolonge aussi $A^t$. Plus pr\'ecis\'ement,
$$u \in D(A_{max}) \Longleftrightarrow \exists v \in
L^2(F) {\rm \ tel\ que\ } \forall \phi \in C^\infty_0(F),\ \<u,A^t\phi\>_E =
\<v,\phi\>_F,$$
ce qui signifie exactement que $v=Au$ ``au sens des
distributions''. En utilisant des techniques standards d'analyse
(convolution), on montre qu'on peut approcher
$u \in D(A_{max})$ par des sections lisses, i.e. (et on peut prendre \c{c}a comme
d\'efinition)
$$u \in D(A_{max}) \Leftrightarrow  \exists (u_n) \in
C^\infty(E) {\rm \ telle\ que\ } \lim_{n \to \infty} u_n = u {\rm \ et\ que\
  la\ suite\ } (Au_n)  {\rm \ converge\ dans\ } L^2$$ ($A_{max} u$ est alors
la valeur de cette limite).

\subsection{Formules de Stokes}

Pour faire fonctionner la technique de Bochner on a besoin de proc\'eder \`a
des int\'e\-gra\-tions par parties. Les deux r\'esultats suivants ainsi que leur
interpr\'etation en termes d'op\'erateurs non born\'es sont \`a notre disposition.
Le premier th\'eor\`eme d'int\'egration par parties sur une c\^one-vari\'et\'e est le
suivant, d\^u \`a Cheeger \cite{Cheeger} :

\begin{Thm}\label{Cheeg} Soient $\eta \in \Omega^p M$ et $\sigma \in
  \Omega^{p+1}M$ deux formes $C^\infty$ sur $M$ telles que $\eta$,
  $d\eta$, $\sigma$, et $\delta \sigma$ soient dans $L^2$. Alors
$$\< \eta, \delta \sigma \> = \< d\eta, \sigma\>.$$
\end{Thm}

En fait il faut adapter un tout petit peu la d\'emonstration, ou combiner deux
r\'esultats de l'article cit\'e (cf aussi \cite{HK}, appendice).

\bigskip

Le deuxi\`eme r\'esultat concerne les tenseurs et non plus les formes
diff\'erentielles :

\begin{Thm} \label{ipp} Soient $u \in C^\infty(T^{(r,s)}M)$, $v \in
  C^\infty(T^{(r+1,s)}M)$ tels que $u$, $\n u$, $v$, $\na v$ soient dans
  $L^2$. Alors $$\<u,\na v\> = \<\n u,v\>.$$
\end{Thm}

\begin{proof}
On va d\'emontrer ce r\'esultat en utilisant une m\'ethode similaire \`a
  celle de Cheeger \cite{Cheeger}. Pour simplifier, on supposera que la
  m\'etrique au voisinage de $\Sigma$ est exactement, en coordonn\'ees locales, de
  la forme $g=dr^2+r^2d\theta^2+g_{|\Sigma_i}$, o\`u $\theta$ est d\'efinie modulo
  l'angle conique $\alpha$. Le cas g\'en\'eral se traite exactement de la m\^eme
  fa\c{c}on, les expressions sont juste un peu plus compliqu\'ees.

Soit $a$ un r\'eel positif suffisamment petit pour que le
$a$-voisinage ferm\'e de $\Sigma$ dans $M$ soit tubulaire. Pour $t\leq a$, on note $U_t$ le 
$t$-voisinage de $\Sigma$ dans $M$, $\Sigma_t = \partial U_t$,
et $M_t = M\setminus U_t$. Le vecteur $\frac{\partial}{\partial r} = e_r$
est une normale unitaire en tout point de $\Sigma_t$. Avec ces notations,
une int\'egration par parties ( c'est-\`a-dire la formule de Stokes ) nous donne :
\begin{eqnarray} \int_{M_t} (g(u,\na v)-g(\n u,v)) = \int_{\Sigma_t} g_{|\Sigma_t}(u,i_{e_r}
v) \label{0016} \end{eqnarray}
o\`u $i_{e_r} v = v(e_r,.)$. Le terme de gauche converge vers $\<u,\na v\> - \<\n
u,v\>$ quand $t$ tend vers $0$.
Notons $I_t$ le terme de droite de l'\'egalit\'e, qui
  correspond au terme de bord. On a alors les in\'egalit\'es suivantes (la
  notation $|.|$ d\'esigne la valeur absolue aussi bien que la norme ponctuelle
  pour la m\'etrique $g$) :
\begin{eqnarray*}
|I_t| & \leq & \int_{\Sigma_t} |u| |i_{e_r} v| \\
      & \leq & \left(\int_{\Sigma_t} |u|^2\right)^{1/2} \left(\int_{\Sigma_t}
      |i_{e_r} v|^2\right)^{1/2}
\end{eqnarray*}

On va montrer que le fait que $\n u$ soit $L^2$ permet d'avoir une
bonne majoration de $\int_{\Sigma_t} |u|^2$. Et comme $i_{e_r}v$ est $L^2$
(car $v$ l'est aussi), $\int_{\Sigma_t} |i_{e_r}v|^2$ ne
peut pas cro\^itre trop vite quand $t$ tend vers $0$. Ce deux r\'esultats nous
permettront d'affirmer que $I_{t_n}$ tend vers $0$ pour une suite $t_n$
tendant vers $0$.

En tout point o\`u $|u| \neq 0$, la fonction $|u|$ est d\'erivable, et
$d|u|(x) = g(\n _x u, \frac{u}{|u|})$. On pose 
$$\frac{\partial
|u|}{\partial r} = g(\n_{e_r} u, \frac{u}{|u|})$$ 
si $|u| \neq 0$, et $\frac{\partial |u|}{\partial r} =0$ sinon. Il s'agit de la
d\'eriv\'ee partielle distributionnelle de $|u|$, au sens o\`u l'on a,
si les coordonn\'ees autres que $r$ restent fix\'ees, 
$$|u(t)| - |u(a)| = \int^t_a \frac{\partial |u|}{\partial r}(r) dr.$$ 
Or $|\frac{\partial |u|}{\partial r}| \leq |\n _{e_r} u|$, 
et donc, si $t\leq a$, 
$$|u(t)| \leq |u(a)| + \int_t^a |\n_{e_r} u| dr$$
 et 
 $$|u(t)|^2 \leq 2|u(a)|^2 + 2(\int_t^a |\n_{e_r} u| dr)^2.$$

En appliquant l'in\'egalit\'e de Cauchy-Schwarz il vient
\begin{eqnarray*}
\left(\int_t^a |\n_{e_r} u|dr \right)^2 &\leq &\int_t^a \frac{dr}{r}
\int_t^a r|\n_{e_r} u|^2 dr\\ & \leq & |\ln(\frac{t}{a})| \int_t^a r|\n
u|^2 dr
\end{eqnarray*}

En int\'egrant sur $\Sigma_t$ on trouve
\begin{eqnarray}
\int_{\Sigma_t} |u|^2 & \leq & \int_{\Sigma_t} 2|u(a)|^2 +
\int_{\Sigma_t}\left( 2 \ln(\frac{t}{a}) \int_t^a r|\n_{e_r} u|^2
dr\right) \nonumber \\ 
& \leq & 2\int_{\Sigma_t} |u(a)|^2 + 2
|\ln(\frac{t}{a})| \int_{\Sigma_t}\int_t^a r|\n_{e_r} u|^2 dr \nonumber \\ 
& \leq & 2\int_{\Sigma}\int_{\theta=0}^\alpha |u(a)|^2td\theta
dvol_\Sigma + 2 |\ln(\frac{t}{a})| \int_{\Sigma}
\int_{\theta=0}^\alpha \left(\int_t^a r|\n_{e_r} u|^2
  dr\right) td\theta dvol_\Sigma \nonumber \\
& \leq &  2\frac{t}{a}\int_{\Sigma_a} |u|^2 +
2t|\ln(\frac{t}{a})|\int_{\Sigma} \int_{\theta=0}^\alpha \int_t^a
|\n_{e_r} u|^2 rdrd\theta dvol_\Sigma \nonumber \\ 
& \leq & 2\frac{t}{a}\int_{\Sigma_a} |u|^2 + 2t|\ln(\frac{t}{a})|
\int_{U_a} |\n_{e_r} u|^2 \nonumber \\ 
& \leq & 2\frac{t}{a}\int_{\Sigma_a} |u|^2 + 2t|\ln(\frac{t}{a})|
\int_{U_a} |\n u|^2 \nonumber \\ 
& = & O(t|\ln(t)|) \label{0017}
\end{eqnarray}

Il reste \`a contr\^oler le terme $\int_{\Sigma_t} |i_{e_r}v|^2 \leq
\int_{\Sigma_t} |v|^2$. Comme $v$ est $L^2$, $$\int_0^a
(\int_{\Sigma_t} |v|^2)dt = \int_{U_a} |v|^2 < +\infty,$$ et donc
la fonction $t \mapsto \int_{\Sigma_t} |v|^2$ est int\'egrable sur
$]0,a]$. Or la fonction $(t\ln(t))^{-1}$ n'est pas int\'egrable en
$0$. On en d\'eduit donc qu'il existe une suite $t_n$ tendant vers
$0$ pour laquelle 
$$\int_{\Sigma_{t_n}} |v|^2 = o((t_n\ln(t_n))^{-1}).$$ 
En combinant avec la majoration \eqref{0017}
pour $\int_{\Sigma_t} |u|^2$, on en d\'eduit imm\'ediatement que
$${\lim_{n\to +\infty} I_{t_n} = 0}.$$ 
Or $I_t \to \<u,\na v\> - \<\n u,v\>$ quand $t \to 0$; on a donc bien  $\<u,\na v\> = \<\n
u,v\>.$
\end{proof}

\bigskip

{\em Remarque :} dans les th\'eor\`emes ci-dessus, la condition $L^2$ para\^it
naturelle, ne serait-ce que pour s'assurer de l'existence des termes du type
$\<\n u,v\>$. Cependant il est int\'eressant de noter que la d\'emonstration
donn\'ee du th\'eor\`eme \ref{ipp} ne fonctionne pas avec des hypoth\`eses du type
$u$, $\n u \in L^p$ et $v$, $\na v \in L^q$, avec $p$ et $q$ des exposants
conjugu\'es. La condition $L^2$ est donc plus importante qu'il n'y para\^it.

\bigskip

Ces deux formules de Stokes ont une interpr\'etation en terme d'op\'erateurs non born\'es. Au vue de la d\'efinition de l'extension maximale d'un op\'erateur, il est clair qu'en passant \`a la limite dans les formules de Stokes \ref{Cheeg} et \ref{ipp}, on obtient
$$\<\eta, \delta_{max} \sigma \> = \< d_{max}\eta, \sigma\> {\rm \ et}$$ 
$$\<u, \na_{max} v \> = \<\n_{max} u,v\> $$
quel que soient $\eta \in D(d_{max})$, $\sigma \in D(\delta_{max})$, $u \in D(\n_{max})$ et $v \in D(\na_{max}$.
On en d\'eduit imm\'ediatement (cf aussi \cite{Gaffney}) que :

\begin{Cor}\label{ipp1} Les op\'erateurs $d_{max}$ et $\delta_{max}$ sont
  adjoints l'un de
l'autre; on a $d_{max} = d_{min}$ et $\delta_{max}= \delta_{min}$. \\
 Les op\'erateurs $\n_{max}$ et $\na_{max}$ sont
  adjoints l'un de
l'autre; on a $\n_{max} = \n_{min}$ et $\na_{max}= \na_{min}$.
\end{Cor}

\medskip

Dans la suite les notations $\n$ et $\na$ d\'esigneront donc (sauf exception)
les op\'erateurs $\n_{max} (=\n_{min})$ et $\na_{max}(=\na_{min}).$
Puisqu'il n'y a pas de risque d'ambigu\"{\i}t\'e, on utilisera \'epi\-so\-di\-que\-ment la notation 
$L^{1,2}$ (resp. $L^{2,2}$) \`a la place de $D(\n)$ (resp. $D(\n \circ \n)$).

On emploiera fr\'equemment le corollaire suivant, simple reformulation du
pr\'ec\'edent :

\begin{Cor}\label{approx} Soit $u$ appartenant \`a $D(\n)$, c'est-\`a-dire tel que
  $u$ et $\n u$
  sont $L^2$. Alors il existe une suite $(u_n)$, $C^\infty$ \`a support
  compact, telle que $u_n \to u$ et $\n u_n \to \n u$ dans $L^2$ quand $n \to
  \infty$.
\end{Cor}

\begin{proof} C'est juste la d\'efinition de $u \in D(\n_{min})$.
\end{proof}

\bigskip

Les techniques d\'evelopp\'ees dans la d\'emonstration du th\'eor\`eme \ref{ipp}, et en particulier la 
majoration \eqref{0017}, vont nous permettre de montrer deux autres r\'esultats dans la m\^eme 
lign\'ee. Commen\c{c}ons par la proposition suivante, qui simplifiera plusieurs preuves par la suite :

\begin{Prop}\label{derl2}
Soit $u$ une section de $T^{(r,s)}M$ telle que $u$, $\n_{e_r} u$ et $\na \n u$ soient dans 
$L^2$. Alors $\n u$ appartient \`a $L^2(T^{(r+1,s)}M)$.
\end{Prop}

\begin{proof} La preuve est juste une relecture de la d\'emonstration du th\'eor\`eme \ref{ipp}. Plus 
pr\'eci\-s\'e\-ment, on constate que pour d\'emontrer l'existence d'une suite $t_n$ tendant vers $0$ telle 
que 
$$\lim_{n \to \infty} \int_{\Sigma_{t_n}} g_{|\Sigma_{t_n}}(u,i_{e_r}(v)) = 0,$$
il suffit que $\n_{e_r} u$ et $i_{e_r}(v)$ soient $L^2$. Quand $v=\n u$ comme c'est le cas ici, 
ces deux conditions reviennent \`a une seule, \`a savoir $\n_{e_r} u \in L^2$.

On applique cela \`a l'\'egalit\'e \eqref{0016}
$$\int_{M_{t_n}} g(\n u, \n u) = \int_{M_{t_n}} g(u,\na \n u)- \int_{\Sigma_{t_n}} 
g_{|\Sigma_{t_n}}(u,\n_{e_r} u).$$
Comme $u$ et $\na \n u$ sont $L^2$, le terme $\int_{M_{t_n}} g(u,\na \n u)$ converge vers 
$\<u,\na \n u\>$ quand $t_n$ tend vers $0$, et comme le deuxi\`eme terme de droite tend vers $0$, 
on a
$$\lim_{n \to \infty} \int_{M_{t_n}} g(\n u, \n u) = \<u,\na \n u\>.$$
L'existence de cette limite, plus le fait que $g(\n u, \n u)$ est partout positif, montre que 
$\n u$ appartient \`a $L^2$ et que $||\n u||^2 = \<u,\na \n u\>.$
\end{proof}

\bigskip

Le r\'esultat suivant, moins anecdotique, nous renseigne encore sur la d\'eriv\'ee covariante (je 
remercie chaleureusement Gilles Carron pour m'avoir donn\'e l'id\'ee de la preuve) :

\begin{Thm}\label{imferme}
L'image de l'op\'erateur $\n : L^{1,2}(T^{(r,s)}M)  \to L^2(T^{(r+1,s)})$ est ferm\'ee.
\end{Thm}

\begin{proof}
On va proc\'eder en deux temps : on va regarder ce qu'il se passe en dehors du lieu 
singulier, puis au voisinage du lieu singulier.
On reprend les notations de la d\'emonstration du th\'eor\`eme \ref{ipp}. Soit $a$ un r\'eel 
positif suffisamment petit pour que le $a$-voisinage ferm\'e du lieu singulier $\Sigma$ dans 
$M$ soit tubulaire, et on suppose donn\'e un r\'eel $b$ tel que $0<b<a$ (on verra ensuite comment 
choisir $b$). Pour $t\leq a$, on note $U_t$ le $t$-voisinage de $\Sigma$ dans $M$, $\Sigma_t = 
\partial U_t$, et $M_t = M \setminus U_t$. On pose aussi $\Omega = M_b = M \setminus U_b$.

\bigskip

Commen\c{c}ons par regarder ce qui se passe sur $\Omega$. C'est un domaine de $M$, \`a bord lisse. On 
va montrer que l'image de $\n$ (ou plut\^ot de son extension maximale) y est ferm\'ee.

Sur $\Omega$, les extensions minimales et maximales de la d\'eriv\'ee covariante sont distinctes. 
On consid\`ere alors l'op\'erateur $L_N = \na_{min|\Omega} \circ \n_{max|\Omega}$ ; il s'agit de 
l'extension de Neumann du laplacien de connexion $\na \n$ sur $\Omega$. Il est 
auto-adjoint (corollaire \ref{1+*b}), positif, de spectre discret. Pour tout $v$ appartenant \`a 
$D(L_N)$, on a
$$\int_\Omega g(L_N v,v) = \int_\Omega g(\n_{max|\Omega}\,v, \n_{max|\Omega}\,v) ;$$
le noyau de $L_N$ co\"{\i}ncide donc avec celui de $\n_{max|\Omega}$, c'est exactement l'ensemble des 
sections parall\`eles sur $\Omega$. Notons $\lambda$ la premi\`ere valeur propre non nulle de $L_N$. 
D'apr\`es la classique formule du minimax, 
$$ \lambda = \inf \left\{ \frac{\int_\Omega g(\n_{max|\Omega}\,v, \n_{max|\Omega}\,v)}{\int_\Omega g(v,v)}\ |\ v \in D(\n_{max|\Omega}) \cap (\ker \n_{max|\Omega})^\perp,\ v\neq 0 
\right\}.$$
En particulier, si $v$ appartient \`a $D(\n_{max|\Omega}) \cap (\ker \n_{max|\Omega})^\perp$, 
alors 
\begin{eqnarray}||v|| \leq (\lambda)^{1/2} ||\n_{max|\Omega}\,v|| \label{0023} \end{eqnarray} 
(il s'agit des normes $L^2$ sur 
$\Omega$). Cette in\'egalit\'e suffit \`a montrer que l'image de $\n_{max|\Omega}$ est ferm\'ee. En 
effet, soit $(v_n)$ une suite de $D(\n_{max|\Omega})$ telle que la suite 
$(\n_{max|\Omega}\,v_n)$ 
converge en norme $L^2$ vers une limite $l$. On note $p_n$ le projet\'e orthogonal de $v_n$ sur 
$(\ker \n_{max|\Omega})^\perp$. Alors $\n_{max|\Omega}\,p_n = \n_{max|\Omega}\,v_n$, et la suite 
des $\n_{max|\Omega}\,p_n$ est donc de Cauchy. En appliquant l'in\'egalit\'e \eqref{0023} \`a $p_{n+k} 
- p_n$, on montre directement que la suite $(p_n)$ est de Cauchy, donc converge vers une limite 
$p$. Maintenant, comme $\n_{max|\Omega}$ est un op\'erateur ferm\'e, la limite $p$ appartient \`a 
$D(\n_{max|\Omega})$ et $\n_{max|\Omega}\,p = l$, ce qui termine de montrer que l'image de 
$\n_{max|\Omega}$ est ferm\'ee.

\bigskip

On va maintenant montrer que sur $U_a$, l'image de l'extension minimale de la d\'eriv\'ee covariante 
est ferm\'ee. On utilise pour cela la 
m\^eme majoration \eqref{0017} que pour la d\'emonstration du th\'eor\`eme \ref{ipp}, qui donne, pour 
toute section $u$ lisse, $L^2$, \`a d\'eriv\'ee covariante $L^2$, 
$$\int_{\Sigma_t} |u(t)|^2 \leq 2\frac{t}{a}\int_{\Sigma_a} |u|^2 + 2t|\ln(\frac{t}{a})| 
\int_{U_a} |\n u|^2.$$
En fait cette in\'egalit\'e est pour une m\'etrique plate sur le c\^one. Dans le cas g\'en\'eral la m\^eme 
majoration est vraie, \`a un facteur multiplicatif $c^2$ pr\`es. 
Si $u$ est \`a support compact dans $U_a$, alors $u$ est nul sur $\Sigma_a$, et
en int\'egrant ce qu'il reste entre $0$ et $a$, on obtient
$||u|| \leq \frac{a c}{\sqrt{2}} ||\n u||,$
o\`u $||.||$ d\'esigne la norme $L^2$ sur $U_a$. Par passage \`a la limite, on a
\begin{eqnarray} ||u|| \leq \frac{a c}{\sqrt{2}} ||\n_{min|U_a}\,u|| \label{0024} \end{eqnarray}
pour tout $u$ appartenant \`a $D(\n_{min|U_a})$.
Cette in\'egalit\'e (une sorte de lemme de Poincar\'e) implique, comme pr\'ec\'edemment, que l'image de 
$\n_{min|U_a}$ est ferm\'ee. En effet, soit $(v_n)$ une suite de $D(\n_{min|U_a})$ telle que la 
suite $(\n_{min|U_a}\,v_n)$ converge en norme $L^2$ vers une limite $l$. La suite 
$(\n_{min|U_a}\,v_n)$ est donc de Cauchy, et en appliquant l'in\'egalit\'e \eqref{0024} \`a $v_{n+k} - 
v_n$, on montre 
directement que la suite $(v_n)$ est de Cauchy, donc converge vers une limite $v$, qui v\'erifie 
(l'op\'erateur $\n_{min|U_a}$) \'etant ferm\'e) $v \in D(\n_{min|U_a})$ et $\n_{min|U_a}\,v =l$, ce 
qui termine de montrer que l'image de $\n_{min|U_a}$ est ferm\'ee.

\bigskip

Voyons maintenant comment en d\'eduire un r\'esultat sur tout $M$. Soit $w$ un point de l'adh\'erence 
de l'image de $\n$. Il existe alors une suite $(u_n)$ d'\'el\'ements de $D(\n)$, que l'on peut tous 
choisir $C^\infty$, telle que la suite 
$\n u_n$ converge vers $w$. On regarde alors la restriction de $u_n$ \`a $\Omega$ : $u_{n|\Omega}$ 
appartient \`a $D(\n_{max|\Omega})$, et $\n_{max|\Omega}\,u_{n|\Omega}$ converge vers 
$w_{|\Omega}$. Alors on a vu que la suite des projections orthogonales des $u_{n|\Omega}$ sur 
$(\ker \n_{max|\Omega})^\perp$ \'etait convergente, c'est-\`a-dire qu'il existe une suite $h_n$ de 
sections parall\`eles sur $\Omega$ telle que la suite $(u_{n|\Omega} - h_n)$ converge en norme 
$L^2$ sur $\Omega$.

Si tous les $h_n$ sont nuls (par exemple si le fibr\'e n'admet pas de sections parall\`eles 
locales), ou si tous les $h_n$ se prolongent en des sections parall\`eles sur tout $M$, on peut 
terminer la preuve facilement. A la place de la suite $(u_n)$, on s'int\'eresse \`a la suite 
des $v_n = u_n - h_n$, et on vient de voir que cette suite converge en norme $L^2$ sur $\Omega$.

On choisit ensuite une fonction de coupure $\rho$, telle que $\rho$ soit lisse, \`a support dans 
$U_a$, et identiquement \'egale \`a $1$ sur $U_b$. Alors pour tout $n$, $\rho v_n$ appartient \`a 
$D(\n_{min|U_a})$, et
$$\n_{min|U_a}\,\rho v_n = \n \rho v_n = \rho \n v_n + d\rho \otimes v_n = \rho \n u_n + d\rho 
\otimes v_n.$$
Comme $\n u_n$ converge vers $w$ et que $\rho$ est born\'ee, la suite $(\rho \n u_n)$ converge en 
norme $L^2$ sur $U_a$ vers $\rho w$. D'autre part $d\rho$ est born\'ee, \`a support dans $U_a 
\setminus U_b$ qui est inclus dans $\Omega$ ; et comme la suite $(v_n)$ est convergente sur 
$\Omega$, donc sur $U_a \setminus U_b$, la suite $(d\rho \otimes v_n)$ est aussi convergente
en norme $L^2$ sur $U_a$. La suite $(\n_{min|U_a}\,\rho v_n)$ est donc convergente sur $U_a$, ce 
qui implique directement \`a l'aide de l'in\'egalit\'e \eqref{0024} la convergence de la suite $\rho 
v_n$ sur $U_a$. Comme $\rho$ est identiquement \'egale \`a $1$ sur $U_b$, on en d\'eduit la 
convergence de la suite $(v_{n|U_b})$ en norme $L^2$ sur $U_b$. Or on a vu qu'en plus 
$(v_{n|\Omega})$ converge sur $\Omega = M\setminus U_b$. La suite $(v_n)$ 
converge donc en norme $L^2$ sur $M$ entier vers une limite $v$. Et comme l'op\'erateur $\n$ est 
ferm\'e, $v$ appartient \`a $D(\n)$ et v\'erifie $\n v = w$, ce qui termine de montrer dans ce cas que 
$\im \n$ est ferm\'ee.

\medskip

Pour conclure, il suffit de montrer que l'on peut toujours choisir $b$ tel
que toute section parall\`ele sur $\Omega = M\setminus U_b$ se prolonge en section parall\`ele sur 
tout $M$. Pour cela, on utilise les deux faits suivants. D'une part, l'espace des sections 
parall\`eles sur $M_t = M\setminus U_t$ est un espace vectoriel de dimension finie. En effet, une 
section parall\`ele sur (une composante connexe d') un ouvert est enti\`erement d\'etermin\'ee par sa 
valeur en un point. D'autre part, si $0<t<t'<a$, alors $M_{t'} \subset M_t$, et l'espace des 
sections parall\`eles sur $M_t$ est inclus dans l'espace des sections parall\`eles sur $M_{t'}$ 
(plus formellement, s'injecte via l'application de restriction dans...). Ainsi quand $t$ tend en 
d\'ecroissant vers $0$, l'espace des sections parall\`eles sur $M_t$, qui est de dimension finie, 
d\'ecro\^it aussi, donc finit par \^etre constant. C'est-\`a-dire qu'il existe un r\'eel $t_0$, $0<t_0<a$, 
tel que pour tout $t$ inf\'erieur \`a $t_0$, l'espace des sections parall\`eles sur $M_t$ et sur 
$M_{t_0}$ co\"{\i}ncide : toute section parall\`ele sur $M_{t_0}$ se prolonge en section parall\`ele sur 
$M_t$. Et comme cette derni\`ere propri\'et\'e est vraie pour tout $t \in ]0,t_0]$, on en d\'eduit que 
toute section parall\`ele sur $M_{t_0}$ se prolonge en une section parall\`ele de $M$. On 
choisit alors $b = t_0$ pour terminer la d\'emonstration.
\end{proof}

\bigskip

On en d\'eduit imm\'ediatement le corollaire suivant :

\begin{Cor}\label{corimferme}
On consid\`ere l'op\'erateur $\n : L^{1,2}(T^{(r,s)}M)  \to L^2(T^{(r+1,s)})$. Il existe une 
constante $c$ telle que pour tout \'el\'ement $\xi \in \im \n$, il existe une section $\zeta \in 
L^{1,2}(T^{(r,s)}M)$ v\'erifiant $\n \zeta = \xi$ et $||\zeta|| \leq c ||\xi||$.
\end{Cor}

\begin{proof}
L'espace $L^{1,2}(T^{(r,s)}M)$ est un Hilbert pour le produit scalaire $\<u,v\>_2 = \<u,v\> + 
\<\n u, \n v\>$. L'application $\n$ est alors continue pour la norme $||.||_2$ associ\'ee. Sa 
restriction \`a $(\ker \n)^\perp$ (o\`u l'orthogonalit\'e est au sens du produit scalaire sus-cit\'e) 
est une application lin\'eaire, continue, bijective ; c'est donc un isomorphisme d'apr\`es le 
th\'eor\`eme de l'image ouverte. Il existe alors une constante $c'$ telle que $||\zeta||_2 \leq c' 
||\n \zeta||$ pour tout $\zeta \in (\ker \n)^\perp$. Comme $||\zeta||^2_2 = ||\zeta||^2 + ||\n 
\zeta||^2$, on en d\'eduit imm\'ediatement que $||\zeta|| \leq c ||\n \zeta||$ pour tout $u \in 
(\ker \n)^\perp$, avec $c^2 = {c'}^2 -1$. Pour conclure, \'etant donn\'e $\xi \in \im \n$, il suffit 
de prendre pour $\zeta$ son unique ant\'ec\'edent dans $(\ker \n)^\perp$.
\end{proof}

\section{Normalisation par la jauge de Bianchi}\label{Bianchi}

Comme on l'a vu dans la section \ref{secEin}, pour se d\'ebarrasser des
d\'eformations triviales on cherche \`a imposer la condition de jauge de
Bianchi. Montrer qu'une d\'eformation infinit\'esimale $h_0$ peut se mettre sous une
forme normalis\'ee v\'erifiant la condition de jauge revient \`a trouver une $1$-forme $\eta$ telle que 
$$h_0 - \delta^*\eta \in \ker \beta,$$
 ce qui \'equivaut \`a r\'esoudre l'\'equation
de normalisation \eqref{0021}:
$$\beta \circ \delta^* \eta = \beta h_0.$$

Cette \'equation peut se mettre sous une forme plus lisible. Pour cela, on
utilise le fait que 
$$\n \eta = \delta^* \eta + \frac{1}{2}d \eta$$ 
(il s'agit juste de la d\'ecomposition du $2$-tenseur $\n \eta$ en partie
sym\'etrique et anti-sym\'etrique), que $\delta$ est toujours la restriction de
$\na$ au sous-fibr\'e correspondant, et donc que 
$$\delta \eta = \na \eta = - \tr \n \eta = - \tr \delta^* \eta,$$
la trace de $d \eta$ \'etant nulle puisque $d\eta$ est anti-sym\'etrique. On
 obtient alors
\begin{eqnarray*}
2 \beta (\delta^* \eta) & = & 2 \delta \delta^* \eta + d \tr \delta^*
\eta \\
& = & 2 \na (\n \eta - \frac{1}{2}d \eta) - d \delta \eta \\
& = & 2 \na \n \eta - \delta d \eta - d \delta \eta \\
& = & 2 \na \n \eta - \Delta \eta.
\end{eqnarray*}
 Ici $\Delta = d\delta + \delta d$ est l'op\'erateur de Laplace-Beltrami sur les
 1-formes. Or $\Delta$ et $\na \n$ (parfois nomm\'e {\em laplacien de connexion}) sont
 reli\'es par la classique formule de Weitzenb\"ock 
 \begin{eqnarray}\Delta\eta = \na \n \eta + ric(\eta), \label{0015} \end{eqnarray}
voir \cite{Besse} \S 1.155. En utilisant cette formule, on trouve
\begin{eqnarray*}
2 \beta (\delta^* \eta) & = & \Delta \eta - 2 ric(\eta) \\
& = & \na \n \eta - ric(\eta).
\end{eqnarray*}
Pour une m\'etrique Einstein, $ric(\eta) = c \eta$, et l'expression ci-dessus se 
simplifie encore. Dans le cas qui nous int\'eresse, la constante $c$ vaut $-(n-1)$, et
\begin{eqnarray} 2 \beta (\delta^* \eta) = \na \n \eta + (n-1) \eta. \label{0022}
\end{eqnarray}

On sera donc amen\'e \`a \'etudier l'op\'erateur $L : \eta  \mapsto \na \n \eta +
(n-1) \eta$.

\subsection{Un premier r\'esultat}\label{premres}

Le th\'eor\`eme suivant montre que pour toute d\'eformation infinit\'esimale $h\in L^2(S^2M)$, 
il existe une (unique) $1$-forme $\eta \in D(\delta^*_{min})$ telle que $h - \delta^* \eta$ v\'erifie 
(distributionnellement) la condition de jauge de Bianchi $\beta(h - \delta^* \eta) = 0$. On peut 
donc normaliser en un certain sens toute d\'eformation infinit\'esimale $L^2$, de fa\c{c}on unique.
Ce r\'esultat est assez g\'en\'eral puisqu'on suppose juste que la m\'etrique conique est Einstein \`a 
courbure de Ricci n\'egative; il n'y a pas de limitations \`a la valeur des angles coniques.

\begin{Thm}\label{normgen} Soit $M$ une c\^one-vari\'et\'e Einstein \`a courbure de Ricci n\'egative. On a 
la d\'ecomposition suivante 
$$L^2(S^2M) = \ker \beta_{max} \oplus \im \delta^*_{min}$$
\end{Thm}

\medskip

Quitte \`a multiplier la m\'etrique $g$ par une constante (voir section \ref{secEin}), on va 
supposer dans la suite qu'elle v\'erifie $E(g)= ric(g) +(n-1)g = 0$.
La d\'emonstration de ce th\'eor\`eme n\'ecessite plusieurs r\'esultats interm\'ediaires, regroup\'es dans la 
proposition suivante.

\begin{Prop}\label{prop} \mbox{}
\begin{enumerate}
\item Le sous-espace $\im \delta^*_{min}$ est ferm\'e dans $L^2(S^2M)$.
\item Les domaines $D(\delta^*_{min})$, $D(\n)$ et $D(\beta^t_{min})$ sont \'egaux.
\item Les deux sous-espaces $\im \delta^*_{min}$ et $\ker \beta_{max}$ sont en somme directe 
(i.e. $\im \delta^*_{min} \cap \ker \beta_{max} = \{0\}$), et la projection canonique de $\ker 
\beta_{max} \oplus \im \delta^*_{min}$ sur le deuxi\`eme facteur est une application
lin\'eaire continue.

\end{enumerate}
\end{Prop}

\begin{proof}[D\'emonstration du 1]
Si $\eta$ appartient \`a $C^\infty_0$, alors
\begin{eqnarray*} ||\delta^* \eta||^2 & = & \<\delta^* \eta, \delta^* \eta\>\\
 & = & \< \delta \delta^* \eta,\eta\> \\
 & = & \< \na (\n \eta - \frac{1}{2}d\eta),\eta\> \\
 & = & \< \na \n \eta  - \frac{1}{2}\delta d \eta,\eta \>.
\end{eqnarray*}
On utilise la formule de Weitzenb\"ock \eqref{0015} $\Delta \eta = d\delta \eta +
\delta d \eta = \na \n \eta - (n-1) \eta$ :
\begin{eqnarray}
||\delta^* \eta||^2 & = & \< \na \n \eta  - \frac{1}{2}\delta d \eta,\eta \> \nonumber \\
& = & \< d\delta \eta + \delta d \eta + (n-1) \eta - \frac{1}{2}\delta d
\eta,\eta \> \nonumber \\
 & = & \< d\delta \eta + \frac{1}{2}\delta d \eta  + (n-1) \eta,\eta \> \nonumber \\
 & = & ||\delta \eta||^2 + \frac{1}{2}||d \eta ||^2 + (n-1) ||\eta||^2 \label{0018}
\end{eqnarray}
et donc $$ ||\delta^* \eta||^2 \geq (n-1) ||\eta||^2.$$ Cette in\'egalit\'e
est encore vraie si $\eta \in D(\delta^*_{min})$; il suffit de consid\'erer
une suite $\eta_n$ d'\'el\'ements $C^\infty_0$ avec $\lim_{n \to \infty}
\eta_n = \eta$ et $\lim_{n \to \infty} \delta^* \eta_n = \delta^*_{min} \eta$.

Cette in\'egalit\'e implique imm\'ediatement que $\im \delta^*_{min}$ est
ferm\'ee dans $L^2$. En effet, si $x_n \in \im \delta^*_{min}$ converge
vers $x$ dans $L^2$, alors il existe une suite $\eta_n$ dans
$D(\delta^*_{min})$ telle que $x_n = \delta^*_{min} \eta_n$; la suite $(x_n)$
est de Cauchy, et comme pour tout $n$, $p$ on a
 $$||x_p - x_n|| = ||\delta^*_{min} (\eta_p - \eta_n)|| \geq \sqrt{n-1} ||\eta_p -
\eta_n||,$$
 la suite $(\eta_n)$ est aussi de Cauchy, donc converge vers un
\'el\'ement $\eta$ de $L^2$. On a alors $\lim_{n \to \infty} \eta_n = \eta$ et
$\lim_{n \to \infty} \delta^*_{min} \eta_n = x$; comme l'op\'erateur
$\delta^*_{min}$ est ferm\'e, on en d\'eduit que $\eta \in D(\delta^*_{min})$ et
que $x= \delta^*_{min} \eta$, donc que $x$ appartient \`a $\im \delta^*_{min}$, ce qui montre 
que $\im \delta^*_{min}$ est ferm\'ee.

\bigskip

\noindent {\em D\'emonstration du 2.}
Revenons au calcul \eqref{0018} : si on utilise diff\'eremment
la formule de Weitzenb\"ock, on trouve
\begin{eqnarray*}
||\delta^* \eta||^2 & = & \< \na \n \eta  - \frac{1}{2}\delta d \eta,\eta \>
\\
& = & \frac{1}{2} \< \na \n \eta  + \Delta \eta + (n-1) \eta - \delta d \eta,\eta
\> \\
 & = & \frac{1}{2} \< \na \n \eta  + (n-1) \eta + d \delta \eta,\eta \> \\
 & = & \frac{1}{2}\left( ||\n \eta ||^2 + (n-1) ||\eta||^2 + ||\delta \eta
 ||^2\right),
\end{eqnarray*}
et donc $||\delta^* \eta||^2 \geq \frac{1}{2} ||\n \eta ||^2$, ceci
\'etant valable pour $\eta \in C^\infty_0$. Si $\eta \in D(\delta^*_{min})$,
on prend une suite $\eta_n$ d'\'el\'ements $C^\infty_0$ avec $\lim_{n \to
\infty} \eta_n = \eta$ et $\lim_{n \to \infty} \delta^* \eta_n = \delta^*_{min}
\eta$. La suite $(\delta^* \eta_n)$ est donc de Cauchy, et l'in\'egalit\'e
ci-dessus implique que la suite $(\n \eta _n)$ est aussi de Cauchy,
donc converge vers un \'el\'ement $x$ dans $L^2$. Comme $\eta_n \in
C^\infty_0 \subset D(\n)$ et que l'op\'erateur $\n$ est ferm\'e, on en
d\'eduit que $\eta = \lim_{n \to \infty} \eta_n$ appartient \`a $D(\n)$, et
donc que $D(\n) \subset D(\delta^*_{min})$. On en d\'eduit aussi que
l'in\'egalit\'e $||\delta^* \eta||^2 \geq \frac{1}{2} ||\n \eta ||^2$ est
vraie pour tout $\eta \in D(\delta^*_{min})$.

Maintenant, $\delta^* \eta$ \'etant la partie sym\'etrique de $\n \eta$, on
a aussi $||\delta^* \eta|| \leq ||\n \eta ||$, et le m\^eme argument
montre que $D(\delta^*_{min}) \subset D(\n)$, et donc $D(\delta^*_{min}) =
D(\n)$. Au passage on a aussi 
$$\frac{1}{2}||\n \eta ||^2 \leq ||\delta^*_{min} \eta||^2 \leq ||\n \eta ||^2$$ 
pour tout $\eta \in D(\n)=D(\delta^*_{min})$.

\medskip

Montrons ensuite que $D(\delta^*_{min}) = D(\beta^t_{min})$. 
On rappelle que $\beta^t$ est l'adjoint formel de
$\beta$, c'est-\`a-dire que pour $\eta \in C^\infty_0$, $\beta^t \eta =
\delta^* \eta + \frac{1}{2}(\delta \eta) g$. Maintenant si $\eta \in
D(\delta^*_{min})$, on prend une suite $\eta_n$ d'\'el\'ements $C^\infty_0$
avec $\lim_{n \to \infty} \eta_n = \eta$ et $\lim_{n \to \infty}
\delta^* \eta_n = \delta^*_{min} \eta$. Comme $(\delta \eta_n)g = - (\tr
(\delta^* \eta_n) ) g$, la suite $((\delta \eta_n)g)$ est aussi
convergente, donc $\beta^t \eta_n = \delta^* \eta_n + \frac{1}{2}(\delta
\eta_n) g$ converge dans $L^2$. Ceci implique que $\eta \in
D(\beta^t_{min})$ par d\'efinition de l'extension minimale $\beta^t_{min}$,
et donc $D(\delta^*_{min}) \subset D(\beta^t_{min})$.

R\'eciproquement, si $\eta \in D(\beta^t_{min})$, on prend une suite $\eta_n$
d'\'el\'ements $C^\infty_0$ avec $\lim_{n \to \infty} \eta_n = \eta$ et
$\lim_{n \to \infty} \beta^t \eta_n = \beta^t_{min} \eta$. Comme $\tr
\beta^t \eta_n = \frac{n-2}{2}\delta \eta_n$ et que $n>2$, on en d\'eduit
que la suite $(\delta \eta_n)$ est aussi convergente, ainsi donc que
$\delta^* \eta_n = \beta^t \eta_n - \frac{1}{2}(\delta \eta_n) g$. Par
suite, l'op\'erateur $\delta^*_{min}$ \'etant ferm\'e, $\eta$ appartient bien \`a
$D(\delta^*_{min})$, et $D(\beta^t_{min}) \subset D(\delta^*_{min})$.

\medskip

{\em Remarque :} \label{rem} On peut d\'emontrer exactement de la m\^eme fa\c{c}on que 
$D(\delta^*_{max}) = D(\beta^t_{max})$ ; il suffit juste de remplacer $C^\infty_0$ par 
$C^\infty$ dans les deux derniers paragraphes.

\bigskip

\noindent {\em D\'emonstration du 3.}
Soit $\eta\in D(\delta^*_{min})$ tel que $\delta^*_{min} \eta \in \ker
\beta_{max}$. Alors on a (au sens des distributions au moins) 
$$0 = \beta \delta^* \eta = \frac{1}{2}(\na \n \eta +(n-1) \eta).$$ 
Comme $\eta$ est
$L^2$, $\na \n \eta$ est aussi $L^2$; et comme $\eta \in D(\delta^*_{min}) =
D(\n)$, $\n \eta$ est aussi $L^2$. On peut alors faire une
int\'egration par partie contre $\eta$ pour trouver $$||\n \eta||^2 +
(n-1) ||\eta||^2 = 0$$ et donc $\eta = 0$ et $\delta^*_{min} \eta =0$.
Cela montre que $\im \delta^*_{min}$ et $\ker \beta_{max}$ sont en somme directe.

\medskip

Pour d\'emontrer la deuxi\`eme partie de ce point on aura besoin du lemme technique suivant :

\begin{Lem} Il existe une constante $c>0$ telle que quel que soit
$\eta$ appartenant \`a $D(\delta^*_{min})$, on a
 $$\<\delta^*_{min} \eta, \beta^t_{min} \eta\> \geq c ||\delta^*_{min} \eta||.||\beta^t_{min}
 \eta||.$$
\end{Lem}

\begin{proof}[D\'emonstration du lemme] On commence par une majoration : si $\eta$ est $C^\infty_0$, 
alors
\begin{eqnarray*}||\beta^t \eta|| & = & ||\delta^* \eta + \frac{1}{2}(\delta \eta)
g|| \\
 & \leq & ||\delta^* \eta|| + \frac{n}{2} ||\frac{1}{n}(\delta \eta)g||.
\end{eqnarray*}
Comme $\delta^* \eta$ et $\frac{1}{n}(\delta \eta)g$ sont respectivement
la partie sym\'etrique et la partie en trace de $\n \eta$, on a
$||\delta^* \eta|| \leq ||\n \eta||$ et $||\frac{1}{n}(\delta \eta)g|| \leq
||\n \eta||$, donc
 $$||\beta^t \eta||\leq \frac{n+2}{2}||\n \eta ||.$$

Ensuite, toujours pour $\eta \in C^\infty_0$, on a
\begin{eqnarray*}
\<\delta^* \eta, \beta^t \eta\> & = & \<\beta \circ \delta^* \eta, \eta \> \\
 & = & \frac{1}{2} \<\na \n \eta + (n-1) \eta,\eta\>\\
 & = & \frac{1}{2} ||\n \eta||^2 + \frac{n-1}{2}||\eta||^2\\
 & \geq & \frac{1}{2} ||\n \eta||^2.
\end{eqnarray*}
Comme $||\n \eta|| \geq ||\delta^* \eta||$ et $||\n \eta|| \geq
\frac{2}{n+2} ||\beta^t \eta||$, on a
 $$\<\delta^* \eta, \beta^t \eta\> \geq \frac{1}{n+2}||\delta^* \eta||.||\beta^t
 \eta||.$$

Maintenant si $\eta$ appartient \`a $D(\delta^*_{min})= D(\n)$, on prend
une suite $\eta_n \in C^\infty_0$ telle que $\eta_n$ tend vers $\eta$ et
$\n \eta_n$ tend vers $\n \eta$; alors $\delta^* \eta_n$ et $\beta^t \eta_n$
convergent respectivement vers $\delta^*_{min} \eta$ et $\beta^t_{min} \eta$, et
en passant \`a la limite dans l'in\'egalit\'e ci-dessus on trouve le
r\'esultat voulu avec $c = \frac{1}{n+2}$.
\end{proof}

\medskip

Revenons \`a la d\'emonstration de la continuit\'e de la projection canonique
de $\ker \beta_{max} \oplus \im \delta^*_{min}$ sur le deuxi\`eme facteur.
Ce que le lemme pr\'ec\'edent nous montre, c'est que pour tout \'el\'ement de $\im \delta^*_{min}$, il 
existe un \'el\'ement de $\im \beta^t_{min} \subset (\ker \beta_{max})^\perp$ tel que l'angle entre 
les deux reste \'eloign\'e de $\pi/2$. Cette propri\'et\'e va permettre d'estimer la norme du 
projet\'e sur $\im \delta^*_{min}$ \`a partir de la norme du projet\'e orthogonal sur $(\ker 
\beta_{max})^\perp$. C'est ce qui va assurer la continuit\'e.

\medskip

Soit $x = h + \delta^*_{min} \eta \in \ker \beta_{max} \oplus \im \delta^*_{min}$
(avec $h \in \ker \beta_{max}$), on veut montrer qu'il existe une
constante $C \geq 0$ telle que $||\delta^*_{min} a|| \leq C ||x||$.
Notons $p$ le projet\'e orthogonal de $\delta^*_{min} \eta$ sur $(\ker
\beta_{max})^\perp$ : on a $\delta^*_{min} \eta = p + k$ avec $k \in
\ker(\beta_{max})$. Par d\'efinition de la projection orthogonale,
$||k||^2 = \inf \{||\delta^*_{min} \eta - y||^2\ |\ y
\in\ker(\beta_{max})^\perp \}.$

Pour majorer $||k||^2$ on va choisir un bon $y$. On sait que $\ker(\beta_{max})^\perp = 
\overline{\im \beta^t_{min}}$ ; en
particulier $\beta^t_{min} \eta \in \ker(\beta_{max})^\perp$. Si $\beta^t_{min} \eta
= 0$, alors en prenant la trace on trouve $\delta \eta =0$ puis
$\delta^*_{min} \eta = 0$, et dans ce cas on a bien $||\delta^*_{min} \eta|| \leq C ||x||$.

Si $\beta^t_{min} \eta \neq 0$, on note $p'$
le projet\'e orthogonal de $\delta^*_{min} \eta$ sur ${\rm Vect}(\beta^t_{min}
\eta)$ ; on a
$$p' = \<\delta^*_{min} \eta ,\frac{\beta^t_{min} \eta}{||\beta^t_{min}
\eta||}\>\,\frac{\beta^t_{min} \eta}{||\beta^t_{min} \eta||},$$
et $||\delta^*_{min} \eta||^2 = ||p'||^2 + ||\delta^*_{min} \eta - p'||^2$.
Comme $p' \in (\ker\beta_{max})^\perp$, d'apr\`es la d\'efinition de la projection orthogonale on a
\begin{eqnarray*}
||k||^2 & \leq & ||\delta^*_{min} \eta - p'||^2 \\
& \leq & ||\delta^*_{min} \eta||^2 - ||p'||^2 \\
& \leq & ||\delta^*_{min} \eta||^2 - \<\delta^*_{min} \eta ,\frac{\beta^t_{min} \eta}{||\beta^t_{min}
\eta||}\>^2.
\end{eqnarray*}
En utilisant le lemme pr\'ec\'edent, on trouve
\begin{eqnarray*}
||k||^2 & \leq & ||\delta^*_{min} \eta||^2 - \frac{1}{(n+2)^2} ||\delta^*_{min} \eta||^2\\
& \leq & \left(1-\frac{1}{(n+2)^2}\right) ||\delta^*_{min} \eta||^2,
\end{eqnarray*}
 d'o\`u
 \begin{eqnarray*}
||p||^2 & = & ||\delta^*_{min} \eta||^2 - ||k||^2 \\
& \geq & \left(1 - (1-\frac{1}{(n+2)^2}) \right) ||\delta^*_{min} \eta||^2\\
& \geq & \frac{1}{(n+2)^2} ||\delta^*_{min} \eta||^2.
\end{eqnarray*}

 Ensuite, comme $x = h + \delta^*_{min} \eta = h + k + p$, avec $h$ et
$k$ dans $\ker\beta_{max}$ et $p$ dans $(\ker\beta_{max})^\perp$, on a
\begin{eqnarray*}
||x||^2 & = & ||h+k||^2 + ||p||^2 \\
 & \geq & ||p||^2 \\
 & \geq & \frac{1}{(n+2)^2} ||\delta^*_{min} \eta||^2,
\end{eqnarray*}
soit $||\delta^*_{min} \eta|| \leq (n+2) ||x||$, et on a bien montr\'e ce
qu'on voulait, \`a savoir la continuit\'e de la projection canonique
de $\ker \beta_{max} \oplus \im \delta^*_{min}$ sur $\im \delta^*_{min}$.
\end{proof}

\bigskip

\begin{proof}[D\'emonstration du th\'eor\`eme \ref{normgen}] 
La d\'emonstration se fait en deux \'etapes. On montre
d'abord que $C^\infty_0$ est un sous-espace de $\ker \beta_{max} \oplus \im
\delta^*_{min}$, puis que $\ker \beta_{max} \oplus \im \delta^*_{min}$ est un
sous-espace vectoriel ferm\'e de $L^2$. La densit\'e de $C^\infty_0$ dans $L^2$ permet ensuite de  
conclure.

\medskip

Montrons que $C^\infty_0 \subset \ker \beta_{max} \oplus \im
\delta^*_{min}$ : \\
 Soit $\phi \in C^\infty_0$. On cherche \`a \'ecrire
$\phi = k + \delta^*_{min} \eta$ avec $k \in \ker \beta_{max}$. Or comme $2
\beta \circ \delta^* = \na \n + (n-1) Id$, on peut trouver une
solution de l'\'equation $\beta(\phi) = \beta(\delta^* \eta)$ avec $\eta$,
$\n \eta$, $\na \n \eta$ dans $L^2$. Mais si $\eta$ et $\n \eta$ sont $L^2$,
alors $\eta \in D(\n) = D(\delta^*_{min})$, et on a alors la
d\'ecomposition voulue en \'ecrivant $\phi = (\phi - \delta^*_{min} \eta) +
\delta^*_{min} \eta$.

\medskip

Avec tout le travail pr\'eparatoire qui a \'et\'e fait dans la proposition \ref{prop}, il est 
maintenant facile de montrer que $\ker \beta_{max} \oplus \im
\delta^*_{min}$ est un sous-espace vectoriel ferm\'e de $L^2$. En effet,
si $x_n = h_n + \delta^*_{min} \eta_n$ est une suite de $\ker \beta_{max}
\oplus \im \delta^*_{min}$ convergeant vers $x \in L^2$, alors la
suite $\delta^*_{min} \eta_n$ converge aussi par continuit\'e, ainsi par
cons\'equent que la suit $h_n = x_n - \delta^*_{min} \eta_n$. Or $\im
\delta^*_{min}$ et $\ker \beta_{max}$ sont des sous-espaces vectoriels
ferm\'es de $L^2$ (pour $\im \delta^*_{min}$, c'est le {\it 1.} de la proposition \ref{prop}
ci-dessus; et pour $\ker \beta_{max}$, c'est parce que par d\'efinition
$\ker \beta_{max} = (\im \beta^t_{min})^\perp$). Donc les limites des
suites $(\delta^*_{min} \eta_n)$ et $(h_n)$ sont respectivement dans $\im
\delta^*_{min}$ et $\ker \beta_{max}$, et par suite $x = \lim x_n = \lim
h_n + \lim \delta^*_{min} \eta_n$ appartient \`a $\ker \beta_{max} \oplus \im
\delta^*_{min}$ qui est par cons\'equent ferm\'e.

\medskip

Pour conclure, comme $\ker \beta_{max} \oplus \im \delta^*_{min}$ est
ferm\'e et contient $C^\infty_0$, il contient aussi son adh\'erence
qui est l'espace $L^2$ tout entier, donc $\ker \beta_{max} \oplus \im
\delta^*_{min} = L^2$.
\end{proof}

\bigskip

Dans le cas o\`u la m\'etrique est hyperbolique, on peut encore raffiner un peu ce r\'esultat.

\begin{Prop}\label{minmax}
Si $M$ est une c\^one-vari\'et\'e hyperbolique, alors 
$\beta_{max}= \beta_{min}$, $\delta^*_{max}=\delta^*_{min}$. On note ces op\'erateurs simplement
$\beta$ et $\delta^*$, et le th\'eor\`eme pr\'ec\'edent devient 
$L^2(S^2M) = \ker \beta \oplus \im \delta^*$.
\end{Prop}

\bigskip

La d\'emonstration de cette proposition demande de conna\^itre le comportement des solutions de l'\'equation $\beta \circ \delta^* \eta = 0$, ce qui est l'objet des sections suivantes. Pour cette raison, la preuve est report\'ee en \ref{demminmax}, p. \pageref{demminmax}.

\bigskip

Ces r\'esultats sont encourageants : ils nous montrent que la condition de jauge de Bianchi est une ``bonne'' condition de jauge, puisqu'elle permet de normaliser toutes les d\'eformations infinit\'esimales $L^2$. Malheureusement, ils ne sont pas suffisants, au sens o\`u ils ne nous apprennent rien sur la r\'egularit\'e de la d\'eformation normalis\'ee. En particulier, si l'on part d'une d\'eformation infinit\'esimale $h \in L^{1,2}$ (ce qui est le cas d'une d\'eformation pr\'eservant les angles), alors on a aucune garantie que la d\'eformation normalis\'ee correspondante soit encore $L^{1,2}$. Ce ph\'enom\`ene de perte de r\'egularit\'e est due au caract\`ere singulier de la m\'etrique, qui interdit l'utilisation des r\'esultats classiques de r\'egularit\'e elliptique. 

Pour bien comprendre ce ph\'enom\`ene, on va proc\'eder \`a une \'etude d\'etaill\'ee de l'\'equation de normalisation et de l'op\'erateur correspondant $L = \na\n + (n-1) Id$ au voisinage du lieu singulier, pour une m\'etrique hyperbolique. On verra que le comportement des solutions est intimement reli\'e \`a la valeur des angles coniques. On pourra alors donner des r\'esultats plus pr\'ecis sur la r\'egularit\'e des d\'eformations normalis\'ees dans le cas o\`u les angles coniques sont suffisamment petits.

\subsection{Etude de l'\'equation de normalisation}\label{secLap}

On supposera d\'esormais dans toute la suite de cet article que la m\'etrique conique $g$ est {\em
 hyperbolique}. Le fait de conna\^itre explicitement la forme de la m\'etrique au voisinage du lieu singulier va nous permettre de travailler en coordonn\'ees cylindriques locales et d'effectuer des d\'ecompositions du type s\'eries de Fourier, ramenant ainsi une \'equation aux d\'eriv\'ees partielles \`a des \'equations diff\'erentielles ordinaires sur les coefficients.

Consid\'erons donc  l'op\'erateur $L : \eta  \mapsto \na \n \eta +
(n-1) \eta$.
La premi\`ere chose \`a remarquer sur l'op\'erateur $L$ est qu'il est {\em elliptique}. En
particulier, si $\phi$ est $C^\infty$ et que $L \eta = \phi$ au sens des
distributions, alors $\eta$ est $C^\infty$. Cependant, le caract\`ere
singulier d'une c\^one-vari\'et\'e nous emp\^eche d'utiliser directement les
in\'egalit\'es de type Schauder ou G\r{a}rding. Par exemple, on peut montrer qu'il
existe des 1-formes $\eta$ appartenant \`a $L^2$
telles que $L\eta = 0$ au sens des distributions avec $\n \eta$ qui n'est
pas dans $L^2$. Il va donc falloir faire attention aux domaines sur lesquels on se place.

Le th\'eor\`eme \ref{1+*} nous un tel domaine : avec
les conventions de la section \ref{rapop}, l'op\'erateur $\na \circ \n + (n-1) Id$
est auto-adjoint et inversible. \label{LapDeb} Son domaine est par d\'efinition
$$D = \left\{ \eta \in D(\n)\ |\ \n \eta \in D(\na)\right\} = \left\{ \eta
 \in  L^2\ |\ \n \eta, \ \na \n \eta \in L^2 \right\}$$
(dans la deuxi\`eme expression, il faut consid\'erer $\n$ et $\na \n$ au sens des
distributions). Il s'agit en fait de l'extension de Friedrichs de $L$. Ainsi, pour tout $\phi \in L^2(T^*M)$, il existe une unique $1$-forme $\eta$ telle que $L\eta = \phi$ et que $\eta$, $\n \eta$ et $\na \n \eta$ soient dans $L^2$ (ce qui peut aussi se d\'emontrer directement). Mais ces propri\'et\'es ne sont toujours pas satisfaisantes. On va dans cette section proc\'eder \`a l'\'etude de l'op\'erateur $L$, dans le but d'arriver \`a d\'emontrer, au moins dans certains cas, des propri\'et\'es suppl\'ementaires sur les solutions de l'\'equation de normalisation.

\subsubsection{Expression du laplacien de connexion en coordonn\'ees
  cylindriques}\label{vois1}

Soit $a$ un r\'eel positif suffisamment petit pour que le $a$-voisinage ferm\'e de
$\Sigma$ dans $M$ soit un voisinage tubulaire. Si $r$ est plus petit que $a$,
on note $U_r$ le $r$-voisinage de $\Sigma$ dans $M$ et $\Sigma_r$ le bord de
$U_r$.

Par d\'efinition, si $x$ est un point de $\Sigma$, il existe un voisinage $V$ de
$x$ dans $U_a$ et un voisinage $U$ de $x$ dans $\Sigma$, tels que $U=V\cap
\Sigma$ et $V \simeq U\times D^2$, et dans les coordonn\'ees cylindriques locales
adapt\'ees \`a la d\'ecomposition $V \simeq U\times D^2$, la m\'etrique est de la
forme
$$g=dr^2 + \sinh(r)^2 d\theta^2 + \cosh(r)^2 g_\Sigma,$$
o\`u $\theta$ est d\'efini non pas modulo $2\pi$ mais modulo l'angle conique
$u$. On utilisera les notations
suivantes : $e^r= dr$, $e^\theta = \sinh(r) d\theta$, $e_r = (e^r)^\sharp =
\frac{\d}{\d r}$, et $e_\theta = (e^\theta)^\sharp =
\frac{1}{\sinh(r)}\frac{\d}{\d \theta}$.

Soit $\eta$ une section de $T^*M$. Au voisinage de $\Sigma$ on peut faire une
d\'ecomposition orthogonale et on \'ecrit
$$\eta = f e^r + g e^\theta + \omega,$$
 avec $f$, $g$, deux fonctions de $M$ dans $\mathbb{R}$ (ou $\mathbb{C}$, on
 sera souvent amen\'e dans la suite \`a complexifier les fibr\'es sur lesquels on
travaille), et $\omega$ une 1-forme. On
remarque que bien que les coordonn\'ees ne soient que
locales, les formes $e^r$ et $e^\theta$ sont bien d\'efinies sur tout $U_a$,
ainsi que la d\'ecomposition orthogonale pr\'ec\'edente.

Au vu de la forme de notre voisinage tubulaire, sur tout ouvert $V$ de $U_a$ du
type ci-dessus et suffisamment petit, on peut d\'efinir localement des
champs de vecteurs $e_1,\ldots e_{n-2}$ de telle sorte que
$(e_r,e_\theta,e_1,\ldots e_{n-2})$ forme un rep\`ere mobile orthonorm\'e (local),
v\'erifiant
$$\n_{e_r} e_k = \n_{e_\theta} e_k = 0$$
pour tout $k$ dans $1\ldots n-2$. On d\'efinit de m\^eme des 1-formes locales
$e^1,\ldots e^{n-2}$ telles
que le rep\`ere $(e^r,e^\theta,e^1,\ldots e^{n-2})$ soit le rep\`ere mobile dual du
pr\'ec\'edent.

Avant de commencer les calculs, introduisons encore quelques notations. On
note $N$ le (sous-)fibr\'e vectoriel au-dessus de $U_a$, dont la fibre
au-dessus de $x \in U_a$ est le sous-espace vectoriel de $T^*_x M$
orthogonal \`a $e^\theta$ et $e^r$, et $N^*$ le (sous-)fibr\'e vectoriel au-dessus
de $U_a$, dont la fibre au-dessus de $x \in U_a$ est le sous-espace vectoriel
de $T_x M$ orthogonal \`a $e_\theta$ et $e_r$. La $1$-forme $\omega$ introduite
plus haut est naturellement une section de $N$. Les sections $(e_1,\ldots,
e_{n-2})$ forment localement une base de $N^*$, de m\^eme pour $(e^1,\ldots,
e^{n-2})$ et $N$. Si $s$ est une section de $N^*$,
et $t$ une section de $N$ ou de $N^*$, on note $\n_{\Sigma\,s}t$, ou de fa\c{c}on plus lisible 
$(\n_\Sigma)(s,t)$, la projection orthogonale sur $N$ ou sur $N^*$ de $\n_s t$.

Si $f$ est une fonction de $U_a$, on pose $$d_\Sigma f =
\sum_{k=1}^{n-2} (e_k.f)e^k,$$ et $$\Delta_\Sigma f =  \cosh^2 r
\sum_{k=1}^{n-2} (\n_\Sigma)(e_k,e_k).f - e_k.e_k.f$$ (c'est \`a un
facteur pr\`es l'oppos\'e de la trace de la hessienne de $f$
restreinte \`a $N^*$). Ces deux op\'erateurs sont ind\'ependants du
choix des $e_k$. En fait, avec les notations ci-dessus, dans $V$
on a une identification, \`a $r$ et $\theta$ fix\'e, de $U\times
\{r,\theta\}$ \`a $U\subset \Sigma$, et $N^*$ et $N$ restreints \`a
$U\times \{r,\theta\}$ s'identifient de la m\^eme fa\c{c}on \`a $TU$ et
$T^*U$. Les op\'erateurs ci-dessus correspondent via ces
identifications \`a la diff\'erentielle et au laplacien de $U\subset
\Sigma$ (en fait la m\'etrique $g_\Sigma$ sur $U$ et la m\'etrique $\cosh(r)^2 g_\sigma$ sur $U 
\times \{r,\theta\}$ diff\`erent d'un facteur $\cosh(r)^2$, qui se retrouve dans l'expression de 
$\Delta_\Sigma$).

Il en est de m\^eme pour $\n_\Sigma$, et pour les deux op\'erateurs suivants. Si
$\omega$ est une section de $N$, on pose
\begin{eqnarray*}\delta_\Sigma (\omega) & = & - \cosh^2 r
\sum_{k=1}^{n-2} \n_\Sigma(e_k,\omega)(e_k) \\
& = & \cosh^2 r \sum_{k=1}^{n-2}
  \omega(\n_\Sigma(e_k,e_k)) - e_k.\omega(e_k),
\end{eqnarray*}
et
$$(\na \n)_\Sigma \omega = \cosh^2 r \sum_{k=1}^{n-2}
\n_\Sigma(\n_\Sigma(e_k,e_k), \omega) - \n_\Sigma (e_k, \n_\Sigma( e_k,
\omega)),$$
qui correspondent \`a la codiff\'erentielle et au laplacien de connexion pour les
1-formes de $\Sigma$.

\bigskip

On est maintenant arm\'e pour le calcul explicite de $\na \n \eta$. En utilisant
notre rep\`ere mobile, on a
\begin{eqnarray*}
\na \n \eta & = & -\n_{e_r}\n_{e_r}\eta - \n_{e_\theta}\n_{e_\theta}\eta +
\n_{\n_{e_r}e_r} \eta + \n_{\n_{e_\theta}e_\theta}\eta + \sum_{k=1}^{n-2}
\n_{\n_{e_k}e_k} \eta -\n_{e_k}\n_{e_k}\eta.
\end{eqnarray*}
Comme la m\'etrique conique est hyperbolique, on a les expressions suivantes :
$$\n_{e_r}e_r = 0,\ \n_{e_\theta}e_\theta = -\frac{1}{\tanh(r)}e_r, {\rm \ et\
  } \n_{e_k}e_k = \n_\Sigma(e_k,e_k) - \tanh(r)\,e_r.$$
On v\'erifie aussi que   
$$\begin{array}{ccc}
\n_{e_r} e^r = 0 & \n_{e_r} e^\theta = 0 & \n_{e_r} e^j = 0 \\
\n_{e_\theta} e^r = \frac{1}{\tanh(r)} e^\theta & \n_{e_\theta} e^\theta = -
\frac{1}{\tanh(r)} e^r & \n_{e_\theta} e^j = 0 \\
\n_{e_i} e^r = \tanh(r) e^i & \n_{e_i} e^\theta = 0 & \n_{e_i} e^j = -
\tanh(r) \delta_{ij} e^r + \n_\Sigma (e_i,e^j).
\end{array}$$

On trouve alors que
$$\na \n \eta  =  -\n_{e_r}\n_{e_r}\eta - \n_{e_\theta}\n_{e_\theta}\eta -
 (\frac{1}{\tanh(r)} + (n-2)\tanh(r))\n_{e_r} \eta + \sum_{k=1}^{n-2}
\n_{\n_\Sigma(e_k,e_k)} \eta -\n_{e_k}\n_{e_k}\eta.$$

En rempla\c{c}ant $\eta$ par $f e^r + g e^\theta + \omega$, un calcul explicite nous
donne
l'expression suivante pour les composantes de
$\na \n \eta$; selon $e^r$ :

\begin{multline*} -\frac{\partial^2 f}{\partial r^2} -
\frac{1}{\sinh(r)^2}\frac{\partial^2 f}{\partial \theta ^2}
- \left(\frac{1}{\tanh(r)} + (n-2) \tanh(r)\right) \frac{\partial f}{\partial
  r} + \left(\frac{1}{\tanh(r)^2} + (n-2) \tanh(r)^2\right)f 
 + \frac{1}{\cosh(r)^2} \Delta_\Sigma f \\
 + \frac{2}{\sinh(r) \tanh (r)} \frac{\partial g}{\partial \theta} 
- \frac{2 \tanh(r)}{\cosh(r)^2} \delta_\Sigma \omega, \end{multline*}

selon $e^\theta$ :

$$ -\frac{\partial^2 g}{\partial r ^2} -
\frac{1}{\sinh(r)^2}\frac{\partial^2 g}{\partial \theta ^2}
-\left(\frac{1}{\tanh(r)} + (n-2) \tanh(r)\right)\frac{\partial g}{\partial r}
+ \frac{g}{\tanh(r)^2} + \frac{1}{\cosh(r)^2} \Delta_\Sigma g 
- \frac{2}{\sinh(r) \tanh(r)}\frac{\partial f}{\partial \theta},$$

et selon la composante incluse dans $N$ :

$$ -\n_{e_r} \n_{e_r} \omega - \n_{e_\theta}\n_{e_\theta} \omega
 - \left(\frac{1}{\tanh(r)} + (n-2) \tanh(r)\right) \n_{e_r} \omega +
 \tanh(r)^2 \omega 
 + \frac{1}{\cosh(r)^2} (\na \n)_\Sigma \omega - 2\tanh(r) \,d_\Sigma f$$

\bigskip

Pour pouvoir manipuler cette expression, on va effectuer dans la section
suivante une sorte de d\'ecomposition en s\'eries de Fourier g\'en\'eralis\'ees.

\subsubsection{D\'ecomposition en s\'erie de Fourier g\'en\'eralis\'ee}\label{vois2}

On sait qu'au voisinage du lieu singulier, la m\'etrique $g$ se met
localement sous la forme $$g = dr^2 + \sinh(r)^2 d\theta^2 +
\cosh(r)^2 g_\Sigma.$$ Si la coordonn\'ee $\theta$ \'etait d\'efinie
(toujours modulo l'angle conique $\alpha$) sur tout un voisinage
du lieu singulier, on pourrait faire des d\'ecompositions en s\'eries
de Fourier, du type $$f(r,\theta,z) = \sum f_n(r,z) \exp(2i\pi
n\theta/\alpha).$$ Mais en g\'en\'eral la coordonn\'ee d'angle $\theta$
n'est d\'efinie que localement, ce qui emp\^eche d'\'ecrire de telles
d\'ecompositions. On va donc proc\'eder \`a une autre sorte de
d\'ecomposition; on obtiendra finalement des \'ecritures du type
$$f(r,\theta,z) = \sum f_n(r) \psi_n(\theta,z)$$
 o\`u les $(\psi_n)$
forment une base hilbertienne bien choisie du bord d'un voisinage
tubulaire du lieu singulier.

{\bf Une base hilbertienne adapt\'ee}\label{secbashib}

On se place maintenant au voisinage d'une composante connexe du lieu singulier. Pour
simplifier les notations, la notation $\Sigma$ d\'esigne ici la composante connexe en question 
du lieu singulier, et $\alpha$ l'angle conique correspondant.

Comme pr\'e\-c\'e\-dem\-ment, on choisit un r\'eel positif $a$ suffisamment petit
pour que le $a$-voisinage ferm\'e de $\Sigma$ dans $M$ soit un voisinage
tubulaire. Si $r$ est inf\'erieur ou \'egal \`a $a$,
on note $U_r$ le $r$-voisinage de $\Sigma$ dans $M$ et $\Sigma_r$ le bord de
$U_r$. \label{a}

On va particuli\`erement s'int\'eresser \`a la sous-vari\'et\'e $\Sigma_a$. 
Pour pouvoir faire les d\'e\-com\-po\-si\-tions voulues, on veut trouver une
``bonne'' base hilbertienne sur $\Sigma_a$, pour les fonctions
comme pour les $1$-formes, ou plus pr\'ecis\'ement pour les sections
du sous-fibr\'e $N$ d\'efini pr\'ec\'edemment (pour m\'emoire, $N$ est le (sous-)fibr\'e vectoriel au-dessus 
de $U_a$, dont la fibre au-dessus de $x \in U_a$ est le sous-espace vectoriel de $T^*_x M$
orthogonal \`a $e^\theta$ et $e^r$). Pour faciliter les calculs, on sera amen\'e \`a consid\'erer 
plut\^ot les complexifi\'es de ces fibr\'es.

Tout point
$x$ de $\Sigma_a$ admet un voisinage $\mathcal{V}$ de la forme $U \times S^1$,
o\`u $U$ est un ouvert de
$\Sigma$. Dans ce voisinage, la m\'etrique de $\Sigma_a$, induite par celle
de $M$, s'exprime comme une m\'etrique produit ; plus pr\'ecisement on a, dans
les coordonn\'ees adapt\'ees,
$$g_a = \sinh(a)^2 d\theta^2 + \cosh(a)^2 g_\Sigma,$$
o\`u la variable $\theta$ est d\'efinie modulo l'angle conique $\alpha$.
La vari\'et\'e $\Sigma_a$ est donc localement un produit, avec la m\'etrique correspondante.
C'est un cas simple de submersion riemannienne \`a fibres totalement g\'eod\'esiques, la base \'etant 
$\Sigma$, munie de la m\'etrique $\cosh(a)^2 g_\Sigma$, et la fibre $S^1$, muni de la m\'etrique 
$\sinh(a)^2 d\theta^2$. La th\'eorie spectrale de telles vari\'et\'es a d\'ej\`a \'et\'e \'etudi\'e, voir par exemple 
\cite{BessBo}, \cite{BouBeBer}. En particulier, les laplaciens sur une telle vari\'et\'e se d\'ecompose en somme d'un ``laplacien vertical'' (ici $-\n_{e_\theta} \n_{e_\theta}$) et d'un ``laplacien horizontal'' (ici $\cosh(a)^{-2}\Delta_\Sigma$ ou $\cosh(a)^{-2}(\na \n)_\Sigma$), commutant entre eux. Le lien entre les fonctions et les $1$-formes est donn\'e par la relation de commutation suivante. La m\'etrique sur $\Sigma$ \'etant hyperbolique, on peut utiliser la formule de Weitzenb\"ock \eqref{0015}
$$(\na \n)_\Sigma = \Delta_\Sigma + (n-3) Id = d_\Sigma \delta_\Sigma + \delta_\Sigma d_\Sigma + (n-3) Id,$$
valable pour les sections de $N$. On a alors
$$(\na \n)_\Sigma \circ d_\Sigma = (\Delta_\Sigma +
(n-3)Id)\circ d_\Sigma = d_\Sigma \circ (\Delta_\Sigma + (n-3)Id),$$
et donc 
\begin{equation}(\na \n)_\Sigma (d_\Sigma f) = d_\Sigma(\Delta_\Sigma f ) +
(n-3)d_\Sigma f .\label{0030}\end{equation}

On pose $\gamma = \frac{2 \pi}{\alpha}$ ; cette quantit\'e interviendra tr\`es 
fr\'equemment dans la suite. La proposition suivante se d\'eduit directement des r\'esultats des articles \cite{BessBo}Êet \cite{BouBeBer} cit\'es et de la relation de commutation \eqref{0030}~:

\bigskip

\begin{Prop}\label{bashib} 
Il existe une base hilbertienne $(\psi_j)_{j\in \mathbb{N}}$ du complexifi\'e de
$L^2(\Sigma_a)$, telle que pour tout indice $j$, il existe un r\'eel $\lambda_j
\geq 0$ et un entier relatif $p_j$, pour lesquels
$$\begin{cases} \Delta_\Sigma \psi_j = \lambda_j \psi_j \\
e_\theta.\psi_j = \frac{ip_j\gamma}{\sinh(a)}\psi_j.
\end{cases}$$

\medskip

Soit $J$ l'ensemble des $j$ pour lesquels $\lambda_j > 0$. Il existe une base
hilbertienne $(\phi_j)_{j\in J } \cup (\varphi_j)_{j\in \mathbb{N}}$ du
complexifi\'e de $L^2(N)$, telle que :
\begin{itemize}
\item pour tout indice $j$ appartenant \`a $J$, $\phi_j = \frac{\cosh(a)}{
    (\lambda_j)^{1/2}} d_\Sigma \psi_j$, et donc
$$\begin{cases} (\na \n)_\Sigma \phi_j = (\lambda_j + n - 3) \phi_j \\
\n_{e_\theta}\phi_j = \frac{ip_j\gamma}{\sinh(a)}\phi_j\\
\delta_\Sigma \phi_j = \cosh(a) (\lambda_j)^{1/2} \psi_j;
\end{cases}$$
\item pour tout indice $j \in \mathbb{N}$, il existe un r\'eel $\mu_j$ et un
  entier relatif $p'_j$, pour lesquels
$$\begin{cases} (\na \n)_\Sigma \varphi_j = \mu_j \varphi_j \\
\n_{e_\theta}\varphi_j = \frac{ip'_j\gamma}{\sinh(a)}\varphi_j,
\end{cases}$$
et on a de plus $\delta_\Sigma \varphi_j = 0$.
\end{itemize}
\end{Prop}

\bigskip

On va utiliser ces r\'esultats pour proc\'eder \`a la
d\'ecomposition de $\eta = f e^r + g e^\theta + \omega$ sur tout $U_a$.

Pour passer de $\Sigma_a$ \`a $\Sigma_r$, on utilise le transport
parall\`ele et le flot le long des g\'eod\'esiques, int\'egrales du champ de vecteur
$e_r$. Cela
revient \`a \'etendre \`a tout $U_a$ les fonctions $\psi_j$ et les formes $\phi_j$,
$\varphi_j$, en demandant seulement que $e_r.\psi_j =0$, et que $\n_{e_r}\phi_j
= \n_{e_r}\varphi_j = 0$. On note encore $\psi_j$, $\phi_j$ et $\varphi_j$ ces
extensions.

En proc\'edant \`a de simples
changement d'\'echelle, on montre que ces fonctions et formes \'etendues se comportent de la fa\c{c}on 
suivante :
 $$\begin{cases}\displaystyle \frac{\d}{\d r}\psi_j = 0\\
 \displaystyle \frac{\d}{\d \theta}\psi_j = ip_j\gamma\,\psi_j\\
 \displaystyle \Delta_\Sigma \psi_j = \lambda_j \psi_j\\
 \displaystyle d_\Sigma \psi_j = \frac{(\lambda_j)^{1/2}}{\cosh(r)}
 \phi_j, \end{cases}$$
$$\begin{cases}
 \displaystyle \n_{\frac{\d}{\d r}} \phi_j = 0\\
 \displaystyle \n_{\frac{\d}{\d \theta}} \phi_j = ip_j\gamma\,\phi_j \\
 \displaystyle (\na \n)_\Sigma \phi_j =  (\lambda_j + n-3) \phi_j \\
 \displaystyle \delta_\Sigma \phi_j = \cosh(r) (\lambda_j)^{1/2} \psi_j,
\end{cases}$$
  et
$$\begin{cases} \displaystyle \n_{\frac{\d}{\d r}} \varphi_j = 0\\
 \displaystyle \n_{\frac{\d}{\d \theta}} \varphi_j = ip'_j\gamma\,\varphi_j\\
 \displaystyle (\na \n)_\Sigma \varphi_j = \mu_j \varphi_j\\
\displaystyle  \delta_\Sigma \varphi_j = 0\end{cases}.$$

\bigskip

Pour un $r$ fix\'e, on note $f_r$ la restriction de $f$ \`a
$\Sigma_r$. On peut de m\^eme \'etendre $f_r$ en une fonction
$\tilde{f}_r$ d\'efinie sur tout $U_a$ en utilisant le flot du champ
de vecteur $e_r$, c'est-\`a-dire en demandant seulement que
$e_r.\tilde{f}_r$ soit identiquement nul (et \'evidemment que
$\tilde{f}_r=f_r$ sur $\Sigma_r$). En particulier, on peut
regarder la restriction \`a $\Sigma_a$ de $\tilde{f}_r$, not\'ee
$\tilde{f}_r|_{\Sigma_a}$ . On peut maintenant utiliser les
r\'esultats de la proposition \ref{bashib} pour d\'ecomposer
$\tilde{f}_r|_{\Sigma_a}$ sous la forme d'une s\'erie :
$\tilde{f}_r|_{\Sigma_a} = \sum f_r^j \psi_j$. Finalement, en
r\'eutilisant le flot pour se ramener \`a $\Sigma_r$, on obtient la
d\'ecomposition suivante, valable sur $\Sigma_r$ : $f_r = \sum f_r^j
\psi_j$. En faisant cette manipulation pour tout $r$, et en posant
$f_j(r) = f_r^j$, on obtient $$f = \sum_{j\in \mathbb{N}} f_j(r)
\psi_j.$$

\medskip

On effectue \'evidemment une d\'ecomposition similaire pour la fonction $g$. Pour
la section $\omega$, le
m\^eme proc\'ed\'e fonctionne, en rempla\c{c}ant le flot par le transport parall\`ele, et
on obtient une d\'ecomposition $$\omega = \sum_{j\in J} \omega_j(r) \phi_j +
\sum_{j\in \mathbb{N}} \varpi_j(r) \varphi_j.$$

On peut v\'erifier facilement que si $\eta$ est $C^\infty$ alors les coefficients
$f_j$, $g_j$, $\omega_j$ et $\varpi_j$ le sont aussi (en effet,
$f_j(r) = \int_{\Sigma_a} \overline{\psi_j}\tilde{f}_r$ et on peut d\'eriver sous
l'int\'egrale; il en est de m\^eme pour les autres coefficients).

On a finalement obtenu l'expression suivante pour $\eta$ :
\begin{eqnarray*}
\eta & = & \sum_{j\in \mathbb{N}} f_j(r) \psi_j\, e^r + \sum_{j\in \mathbb{N}}
g_j(r) \psi_j\, e^\theta \\
& + & \sum_{j\in J} \omega_j(r) \phi_j + \sum_{j\in \mathbb{N}}
\varpi_j(r) \varphi_j.
\end{eqnarray*}
Il est plus judicieux de regrouper les termes de cette
d\'ecomposition de la fa\c{c}on suivante, faisant appara\^itre des ``blocs
\'el\'ementaires'' de m\^eme fr\'equence :
\begin{eqnarray}
\eta & = & \sum_{j\in J} \left( f_j(r) \psi_j\, e^r + g_j(r) \psi_j\, e^\theta
+ \omega_j(r) \phi_j \right) \nonumber \\
& + & \sum_{j \in \mathbb{N}\setminus J} \left( f_j(r) \psi_j\, e^r + g_j(r)
  \psi_j\, e^\theta \right) \nonumber \\
& + & \sum_{j\in \mathbb{N}} \varpi_j(r) \varphi_j. \label{0006}
\end{eqnarray}

La norme $L^2$ de $\eta$ sur $U_a$ s'exprime bien dans cette d\'ecomposition : on a 
$$||\eta_{|U_a}||^2 = \sum_j \int_0^a \left( |f_j|^2 + |g_j|^2 + |\omega_j|^2 + |\varpi_j|^2 
\right) \frac{\sinh(r)}{\sinh(a)}\left(\frac{\cosh(r)}{\cosh(a)}\right)^{n-2} dr.$$

{\bf Expression du laplacien dans cette d\'ecomposition}

En partant de cette expression pour $\eta$, on va effectuer la m\^eme
d\'ecomposition pour $\na \n \eta + (n-1) \eta$. On note toujours $\gamma$ pour
$\frac{2\pi}{\alpha}$. 

On obtient alors, pour la composante de $\na \n \eta + (n-1)\eta$ en $\psi_j
e^r$, si $j \in J$ : 
\begin{multline}\label{0007} -f_j'' - \left(\frac{1}{\tanh(r)} + (n-2)
\tanh(r)\right) f_j' + \left(\frac{1}{\tanh(r)^2} + (n-2)
\tanh(r)^2 + \frac{p_j^2\gamma^2}{\sinh(r)^2} +
\frac{\lambda_j}{\cosh(r)^2} + n-1 \right)f_j \\ \mbox{} + \frac{2ip_j\gamma}{\sinh(r)
\tanh(r)}g_j - \frac{2\tanh(r)(\lambda_j)^{1/2} }{\cosh(r)} \omega_j, \end{multline}

si $j \notin J$ :
\begin{multline}\label{0008} -f_j'' - \left(\frac{1}{\tanh(r)} + (n-2) \tanh(r)\right) f_j' +
\left(\frac{1}{\tanh(r)^2} + (n-2) \tanh(r)^2 +
\frac{p_j^2\gamma^2}{\sinh(r)^2}+ n-1 \right)f_j \\ 
+ \frac{2ip_j\gamma}{\sinh(r) \tanh(r)}g_j, \end{multline}

pour la composante en $\psi_j e^\theta$ :
\begin{multline}\label{0009} -g_j'' - \left(\frac{1}{\tanh(r)} + (n-2) \tanh(r)\right) g_j' +
\left(\frac{1}{\tanh(r)^2} + \frac{p_j^2\gamma^2}{\sinh(r)^2} + \frac{\lambda_j}{\cosh(r)^2} + 
n-1 \right)g_j \\ - \frac{2ip_j\gamma}{\sinh(r)\tanh(r)}f_j, \end{multline}

pour la composante en $\phi_j$ :
\begin{multline}\label{0010} -\omega_j'' - \left(\frac{1}{\tanh(r)} + (n-2) \tanh(r)\right)
\omega_j' + \left(\tanh(r)^2 + \frac{p_j^2\gamma^2}{\sinh^2(r)} +
\frac{\lambda_j + n-3}{\cosh(r)^2} +n-1 \right)\omega_j \\
- \frac{2\tanh(r) (\lambda_j)^{1/2}}{\cosh(r)} f_j, \end{multline}

et pour la composante en $\varphi_j$ :
\begin{eqnarray} - \varpi_j'' - \left(\frac{1}{\tanh(r)} +
(n-2) \tanh(r)\right) \varpi_j' + \left(\tanh(r)^2 +
\frac{{p'_j}^2 \gamma^2}{\sinh(r)^2} +
\frac{\mu_j}{\cosh(r)^2} +n-1 \right) \varpi_j. \label{0011} \end{eqnarray}

\bigskip

\subsubsection{Comportement des solutions de l'\'equation homog\`ene au voisinage de la
  singularit\'e}\label{vois3}

On va maintenant chercher \`a r\'esoudre l'\'equation $L\eta = 0$ au voisinage de
$\Sigma$. Si la $1$-forme $\eta$ v\'erifie $L\eta=0$ au voisinage du lieu singulier, alors par 
r\'egularit\'e elliptique $\eta$ est localement $C^\infty$ ; cela justifie la d\'ecomposition en s\'erie 
\eqref{0006}, qui converge uniform\'ement sur tout compact du voisinage du lieu singulier. 

La d\'e\-com\-po\-si\-tion de $L\eta$ ci-dessus permet alors de passer d'une
\'equation aux d\'eriv\'ees partielles \`a une infinit\'e d'\'equations diff\'erentielles
ordinaires. R\'esoudre l'\'equation $L \eta = 0$ au voisinage du lieu singulier revient donc \`a r\'esoudre 
une \'equation diff\'erentielle lin\'eaire pour chaque coefficient de la d\'ecomposition.
On peut ainsi \'etudier le comportement de chacun des termes du d\'eveloppement de $\eta$, et il est 
ensuite relativement ais\'e d'en d\'eduire des propri\'et\'es du type $L^2$ pour $\eta$ et ses d\'eriv\'ees au 
voisinage du lieu singulier. 

\medskip

Pour chaque indice $j$, l'\'equation (ou plut\^ot le syst\`eme) que l'on obtient pr\'esente une 
singularit\'e ``r\'eguli\`ere'' en $r=0$. On sait (voir \cite{Wasow}, cf aussi \cite{Mazzeo}) que les 
solutions d'une telle 
\'equation sont des combinaisons lin\'eaires de fonctions de la forme
$r^k f(r)$ avec $f$ une fonction analytique, o\`u les exposants $k$ s'obtiennent
comme racines de l'\'equation indicielle (en cas de racines multiples ou
s\'epar\'ees par des entiers, il faut \'eventuellement
rajouter des termes en $\ln r$ dans l'expression des solutions).

On pose donc, pour un entier $j$ donn\'e,
$$\left\{ \begin{array}{rcl}
f_j(r) & = & r^k(f_0 + f_1 r + f_2 r^2 + \cdots), \\
g_j(r) & = & r^k(g_0 + g_1 r + g_2 r^2 + \cdots), \\
\omega_j(r) & = & r^k(\omega_0 + \omega_1 r + \omega_2 r^2 + \cdots), \\
\varpi_j(r) & = & r^k(\varpi_0 +
\varpi_1 r + \varpi_2 r^2 + \cdots).
\end{array} \right.$$

A partir des expressions \eqref{0007} \`a \eqref{0011}, on obtient les syst\`emes d'\'equations 
indicielles suivants :  si $j \in J$,
$$\left\{ \begin{array}{rcccl}
(-k^2 + 1 + p_j^2\gamma^2)f_0 & + & 2ip_j\gamma g_0 & = & 0\\
-2ip_j\gamma f_0 & + & (-k^2 +1 + p_j^2\gamma^2) g_0 & = & 0\\
& &(-k^2 + p_j^2\gamma^2)\omega_0 & = & 0, \end{array} \right.$$
si $j \notin J$,
$$\left\{ \begin{array}{rcccl}
(-k^2 + 1 + p_j^2\gamma^2)f_0 & + & 2ip_j\gamma g_0 & = & 0\\
-2ip_j\gamma f_0 & + & (-k^2 +1 + p_j^2\gamma^2) g_0 & = & 0,
\end{array} \right.$$
et enfin
$$(-k^2 + {p_j'}^2\gamma^2)\varpi_0 = 0.$$

Commen\c{c}ons par \'etudier le premier syst\`eme, le plus compliqu\'e. Son d\'eterminant vaut (au signe 
pr\`es) $(k^2-p_j^2\gamma^2)(k^2-(p_j\gamma +1)^2)(k^2-(p_j\gamma-1)^2)$. Les valeurs de
l'exposant $k$ pour lesquelles le syst\`eme admet des solutions
non triviales (racines indicielles) sont donc $\pm p_j\gamma \pm 1$ et $\pm
p_j\gamma$. Plus pr\'ecis\'ement,
pour $k=\pm(p_j\gamma +1)$, les coefficients dominants
$(f_0,g_0,\omega_0)$ sont engendr\'es par $(1,-i,0)$, pour
$k= \pm (p_j\gamma -1)$, par $(1,i,0)$, et pour $k= \pm p_j\gamma$, par
$(0,0,1)$. On remarque que l'on a toujours des racines
s\'epar\'ees par des entiers, ce qui peut rajouter des termes logarithmiques, mais on
n'aura pas \`a en tenir compte car seul l'exposant dominant va nous
int\'eresser.

Le cas des racines doubles est un peu plus compliqu\'e. Elles apparaissent si $p_j
= 0$, $p_j\gamma = \pm 1$, ou $p_j\gamma = \pm \frac{1}{2}$. En fait si $p_j\gamma =
\frac{1}{2}$, les solutions correspondant \`a $k= p_j\gamma$ et \`a $k = 1 - p_j\gamma$
sont lin\'eairement ind\'ependantes, on n'a donc pas besoin de termes
logarithmiques; m\^eme chose pour $ p_j\gamma = - \frac{1}{2}$.

Pour $p_j\gamma = 1$, les solutions pour $k = p_j\gamma -1$ et $k = - (p_j\gamma -1)$
sont les m\^emes. On a donc besoin d'un terme logarithmique. M\^eme chose si
$p_j\gamma = -1$.

Enfin, pour $p_j = 0$, il y a trois d\'eg\'en\'erescence. Cependant pour $k = 1$ ou $k
= -1$, on n'a pas de perte de dimension et donc pas besoin de termes
logarithmiques. Par contre, pour $k = 0$, le terme en logarithme est
n\'ecessaire.

On remarque que les deux premiers cas de racines doubles ne se rencontrent que
pour des valeurs
particuli\`eres de l'angle conique. Par contre le dernier cas se rencontre quel
que soit l'angle. C'est l'existence de ces solutions logarithmiques, qui sont
dans $L^2$ mais dont la d\'eriv\'ee covariante ne l'est pas, qui fait que
l'op\'erateur $L$ n'est jamais essentiellement auto-adjoint dans notre cadre.

\bigskip

Les deux syst\`emes restant sont plus simples \`a \'etudier et ne pr\'esentent rien de
nouveau par rapport \`a ce qui pr\'ec\`ede. La proposition suivante regroupe tous ces
r\'esultats :

\medskip

\begin{Prop}\label{dvpt} Soit $\eta$ une solution de l'\'equation $L\eta=0$ sur un
  voisinage d'une composante connexe de $\Sigma$, d'angle conique $\alpha$.
Alors chacun des termes apparaissant dans la d\'e\-com\-po\-si\-tion
\begin{eqnarray*}
\eta & = & \sum_{j\in J} \left( f_j(r) \psi_j\, e^r + g_j(r) \psi_j\, e^\theta
+ \omega_j(r) \phi_j \right) \\
& + & \sum_{j \in \mathbb{N}\setminus J} \left( f_j(r) \psi_j\, e^r + g_j(r)
  \psi_j\, e^\theta \right) \\
& + & \sum_{j\in \mathbb{N}} \varpi_j(r) \varphi_j.
\end{eqnarray*}
est solution de l'\'equation $L\eta=0$ au voisinage de la composante connexe du lieu singulier.

\bigskip

Soit $j$ un indice appartenant \`a $J$. L'ensemble des solutions du
type $$f_j(r) \psi_j\, e^r + g_j(r) \psi_j\, e^\theta +
\omega_j(r) \phi_j$$ forme un espace vectoriel (de dimension 6).
Si $p_j\gamma \notin \{-1,0,1\}$, alors on dispose d'une base
constitu\'ee de solutions \'el\'ementaires pour lesquelles
$v(r)=(f_j(r),g_j(r),\omega_j(r))$ est de la forme $r^k(v_0+v_1
r+\cdots)$, avec $k\in \{\pm p_j \gamma \pm 1, \pm p_j \gamma \}$.
Pour $k=\pm(p_j\gamma +1)$, on peut prendre $v_0 = (1, -i,0)$, pour
$k= \pm(p_j\gamma-1)$, $(1,i,0)$, et pour $k=\pm p_j\gamma$,
$(0,0,1)$. Si $p_j\gamma=-1$, resp. $1$, resp. $0$, les deux
solutions \'el\'ementaires ci-dessus correspondant \`a $k=0$ sont
identiques, il faut donc rajouter une solution de la forme
$\ln(r)(v_0+v_1r+\cdots)$ avec $v_0=(1,-i,0)$, resp. $(1,i,0)$,
resp. $(0,0,1)$.

\bigskip

Maintenant si l'indice $j$ n'appartient pas \`a $J$, l'ensemble des
solutions du type $$f_j(r) \psi_j\, e^r + g_j(r) \psi_j\,
e^\theta$$ forme un espace vectoriel (de dimension 4). Si
$p_j\gamma \notin \{-1,1\}$, alors on dispose d'une base constitu\'ee
de solutions \'el\'ementaires pour lesquelles $v'(r)=(f_j(r),g_j(r))$
est de la forme $r^k(v'_0+v'_1 r+\cdots)$, avec $k = \pm p_j \gamma
\pm 1$. Pour $k=\pm(p_j\gamma +1)$, on peut prendre $v'_0 = (1,
-i)$, et pour $k= \pm(p_j\gamma-1)$, $v'_0=(1,i)$. Si
$p_j\gamma=-1$, resp. $1$, les deux solutions \'el\'ementaires
ci-dessus correspondant \`a $k=0$ sont identiques, il faut donc
rajouter une solution de la forme $\ln(r)(v'_0+v'_1r+\cdots)$ avec
$v'_0=(1,-i)$, resp. $(1,i)$.

\bigskip

Enfin, pour tout indice $j$, l'ensemble des solutions du type
$\varpi_j(r) \varphi_j$ forme un espace vectoriel (de
dimension 2). Si $p'_j \neq 0$, alors on dispose d'une base
constitu\'ee de deux solutions \'el\'ementaires pour lesquelles
$\varpi_j(r) = r^k(1 + \varpi_1 r +\cdots)$,
avec $k = \pm p_j \gamma$. Si $p'=0$ les deux solutions
\'el\'ementaires ci-dessus sont identiques,  il faut donc rajouter une
solution pour laquelle $\varpi_j(r) = \ln(r)(1 +
\varpi_1 r +\cdots)$.
\end{Prop}

\bigskip

Dans la suite, on supposera toujours que $p_j$ et $p'_j$ sont {\bf positifs} ; en effet une 
simple conjugaison permet de passer de $p_j$ \`a $-p_j$.

\subsection{Cas des angles coniques inf\'erieurs \`a $2\pi$}\label{2pi}

Dans cette sous-section, {\bf tous les angles coniques seront suppos\'es strictement
inf\'erieurs \`a $2\pi$}. En particulier, si $p$ est un entier, alors
soit $p\gamma = 0$, soit $|p\gamma|>1$.

On va \'etudier maintenant quels sont les exposants dominants
possibles pour une solution de l'\'equation $L\eta=0$ au voisinage du
lieu singulier, en fonction des diff\'erentes conditions impos\'ees \`a
$\eta$.

Le premier r\'esultat est le lemme suivant :

\begin{Lem}\label{nnd} Soit $M$ une c\^one-vari\'et\'e hyperbolique dont tous les
  angles coniques sont strictement inf\'erieurs \`a $2\pi$. Soit $\eta$ une 1-forme telle que
  $L\eta$ soit \'egal \`a $0$ au voisinage du lieu singulier et que $\eta$ et $\n \eta$
  soient dans $L^2$. Alors $\n d \eta$ est dans $L^2$.
\end{Lem}

\begin{proof} Commen\c{c}ons par montrer que $\na \n d\eta$ est dans $L^2$. Pour cela on utilise la formule de Weitzenb\"ock suivante, valable pour une m\'etrique hyperbolique, qui est un
analogue de la formule \eqref{0015} pour les $1$-formes que l'on a d\'ej\`a
utilis\'ee \`a plusieurs reprises (voir \cite{Besse} \S 1.I) :
$$\forall \omega \in \Omega^2M,\
\na \n \omega = \Delta \omega + 2 (n-2) \omega.$$
En particulier, 
$$\na \n d\eta = \Delta d\eta + 2 (n-2) d\eta = d \Delta \eta + 2(n-2) \eta = d(L \eta) -2 d\eta.$$
Au voisinage du lieu singulier, on a alors $\na \n d\eta  = - 2 d\eta$, et donc $\na \n d\eta$ est dans 
$L^2$. 

D'apr\`es la proposition \ref{derl2}, comme $d\eta$ et $\na \n d\eta$ sont $L^2$, il suffit juste
de montrer que $\n_{e_r} d\eta$ est dans $L^2$ pour prouver que $\n d\eta$ est dans $L^2$. 

\bigskip

Pour cela, on regarde comment les conditions $\eta \in L^2$, $\n \eta \in L^2$ se traduisent sur la 
d\'eveloppement de $\eta$. On choisit un r\'eel $a$ suffisamment petit pour que $L\eta$ soit nul sur 
$U_a$ ; c'est ce $a$ que l'on utilise pour la d\'ecomposition de $\eta$ (voir proposition 
\ref{bashib}). On \'ecrit alors
\begin{eqnarray*}
\eta & = & \sum_{j\in J} \left( f_j(r) \psi_j\, e^r + g_j(r) \psi_j\, e^\theta
+ \omega_j(r) \phi_j \right) \\
& + & \sum_{j \in \mathbb{N}\setminus J} \left( f_j(r) \psi_j\, e^r + g_j(r)
  \psi_j\, e^\theta \right) \\
& + & \sum_{j\in \mathbb{N}} \varpi_j(r) \varphi_j.
\end{eqnarray*}

La norme $L^2$ de $\eta$ sur $U_a$ est donn\'ee par
$$||\eta_{|U_a}||^2 = \sum_j \int_0^a \left( |f_j|^2 + |g_j|^2 + |\omega_j|^2 + |\varpi_j|^2 
\right) \frac{\sinh(r)}{\sinh(a)}\left(\frac{\cosh(r)}{\cosh(a)}\right)^{n-2} dr.$$

Comme $L\eta = 0$ au voisinage du lieu singulier, les fonctions $f_j$, $g_j$ etc. sont 
\'equivalentes \`a $c\,r^{k_j}$ (ou $c\,r^{k_j}\ln(r)$) quand $r$ tend vers $0$.  Pour que la 
quantit\'e ci-dessus soit finie, tous les exposants dominants apparaissant dans la d\'ecomposition 
de $\eta$ (donn\'e par la proposition \ref{dvpt}) doivent \^etre strictement plus grand que $-1$. Or on 
a vu que $k$ est de la forme $\pm p_j\gamma \pm 1$ ou $\pm p_j\gamma$, o\`u $p_j$ est un entier 
que l'on peut supposer positif, et $\gamma$ vaut $2\pi$ divis\'e par l'angle conique de la 
composante connexe du lieu singulier. Par cons\'equent le fait que $\eta$ soit dans $L^2$ \'elimine 
d'embl\'ee les solutions 
avec $k = -p_j\gamma -1$, avec $k = -p_j\gamma$ pour $p_j\neq 0$, et avec $k = -p_j\gamma +1$ 
pour $p_j>1$ (et aussi $p_j=1$ si $\gamma\geq 2$, c'est-\`a-dire si l'angle conique est inf\'erieur 
ou \'egal \`a $\pi$). 

\bigskip

Le fait que $\n \eta$ soit dans $L^2$ impose aussi des conditions sur les exposants possibles.
A partir de la d\'ecomposition ci-dessus de $\eta$, on obtient l'\'ecriture suivante pour $\n \eta$ :
\begin{eqnarray}
\n \eta & = & \sum_{j\in \N} \bigg( f_j' \psi_j\, e^r \otimes e^r + g_j' \psi_j\, e^r \otimes e^\theta + 
\left(\frac{ip_j\gamma}{\sinh(r)} f_j - \frac{1}{\tanh(r)} g_j \right) \psi_j\, e^\theta \otimes e^r \nonumber \\
& & + \left(\frac{1}{\tanh(r)} f_j + \frac{ip_j\gamma}{\sinh(r)} g_j \right) \psi_j\, e^\theta \otimes e^\theta
+ f_j \tanh(r) \psi_j\, \cosh(r)^2 g_\Sigma \bigg) \nonumber \\
& + & \sum_{j \in J} \bigg( \omega_j'\, e^r \otimes \phi_j + \frac{ip_j\gamma}{\sinh(r)} \omega_j\, e^\theta \otimes \phi_j + \bigg(\frac{(\lambda_j)^{1/2}}{\cosh(r)} f_j - \tanh(r) \omega_j \bigg) \phi_j \otimes e^r + \frac{(\lambda_j)^{1/2}}{\cosh(r)} g_j\, \phi_j \otimes e^\theta \nonumber \\
& & + \omega_j\, \delta^*_\Sigma \phi_j \bigg) \nonumber \\
& + & \sum_{j \in \N} \left( \varpi_j'\, e^r \otimes \varphi_j + \frac{ip'_j\gamma}{\sinh(r)} \varpi_j\, e^\theta \otimes \varphi_j - \tanh(r) \varpi_j \varphi_j \otimes e^r + \varpi_j\, \n_\Sigma \phi_j
\right).  \label{0031}
\end{eqnarray}
Chacun des termes de cette expression doit \^etre $L^2$ ; de la m\^eme fa\c{c}on que ci-dessus, cela impose que tous les exposants $k$ apparaissant dans la 
d\'ecomposition de $\eta$ (donn\'e par la proposition \ref{dvpt}) doivent \^etre sup\'erieurs ou \'egaux \`a 
$0$, et qu'il n'y ait pas de termes logarithmiques. Comme les angles coniques sont suppos\'es inf\'erieurs \`a $2\pi$, cela revient \`a dire que les fonctions $f_j$, $g_j$, etc. sont des combinaisons lin\'eaires des seules solutions \'el\'ementaires d\'efinies \`a la section pr\'ec\'edente pour lesquelles 
 l'exposant $k$ vaut $0$, $1$, $p_j\gamma -1$, $p_j\gamma$, $p'_j\gamma$ ou
$p_j\gamma+1$.

\medskip

Regardons maintenant $\n_{e_r} d\eta \in L^2$. De la d\'ecomposition en s\'erie \eqref{0031} on d\'eduit facilement l'\'ecriture suivante de $d\eta$ :
\begin{eqnarray*}
d \eta & = & \sum_{j\in \N} \left( g_j' + \frac{1}{\tanh(r)} g_j - \frac{ip_j\gamma}{\sinh(r)} f_j \right)\psi_j\, e^r \wedge e^\theta \\
& & \mbox{} + \sum_{j \in J} \bigg( \bigg(\omega_j' + \tanh(r) \omega_j - \frac{(\lambda_j)^{1/2}}{\cosh(r)} f_j\bigg)\, e^r \wedge \phi_j + \bigg( \frac{ip_j\gamma}{\sinh(r)} \omega_j - \frac{(\lambda_j)^{1/2}}{\cosh(r)} g_j\bigg) \, e^\theta \wedge \phi_j \bigg)\\
& & \mbox{} + \sum_{j \in \N} \bigg( \bigg( \varpi_j'+ \tanh(r) \varpi_j \bigg) \, e^r \wedge \varphi_j + \frac{ip'_j\gamma}{\sinh(r)} \varpi_j\, e^\theta \wedge \varphi_j + \varpi_j\, \d_\Sigma \phi_j \bigg).
\end{eqnarray*}

En d\'erivant par rapport \`a $e_r$, on aboutit \`a l'expression suivante :
\begin{eqnarray*}
\n_{e_r} d\eta  & = & \sum_{j \in \mathbb{N}} \left( g_j'' + \frac{1}{\tanh(r)} g_j' + (1 - \frac{1}{ 
\tanh(r)^2}) g_j - \frac{i p_j\gamma}{\sinh(r)} f_j' + \frac{i p_j\gamma }{ \sinh(r) \tanh(r)} f_j 
\right) \psi_j e^r \wedge e^\theta \\
& + & \sum_{j \in \mathbb{N}\setminus J} \left( \omega_j'' + \tanh(r) \omega_j' + (1-\tanh(r)^2) 
\omega_j - \frac{(\lambda_j)^{1/2} }{ \cosh(r)} f_j' + \frac{(\lambda_j)^{1/2} \tanh(r) }{ \cosh(r)} 
f_j \right) e^r \wedge \phi_j \\
& + & \sum_{j \in \mathbb{N}\setminus J} \left( \frac{i p_j\gamma }{ 
\sinh(r)} \omega_j' - \frac{i p_j\gamma }{ 
\sinh(r) \tanh(r)} \omega_j - \frac{(\lambda_j)^{1/2} }{ \cosh(r)} g_j' + \frac{(\lambda_j)^{1/2} 
\tanh(r) }{ \cosh(r)} g_j \right) e^\theta \wedge \phi_j \\
& + & \sum_{j \in \mathbb{N}} \bigg( \left(\varpi_j'' + \tanh(r) \varpi_j' + (1 - \tanh(r)^2) 
\varpi_j\right) e^r \wedge \varphi_j \\
& & \left. \mbox{} + \left(\frac{i p_j'\gamma }{ \sinh(r)} \varpi_j' - \frac{i 
p_j'\gamma }{ \sinh(r) \tanh(r)} \varpi_j\right) e^\theta \wedge \varphi_j + (\varpi_j' - 
\tanh(r)\varpi_j) d_\Sigma \varphi_j \right)
\end{eqnarray*}

Dans cette somme, seules des combinaisons de solutions \'el\'ementaires dont les exposants dominants 
sont positifs apparaissent. Un simple calcul montre alors que tous les termes de cette s\'erie sont d'exposant sup\'erieur \`a $-1$, et donc appartiennent $L^2$. Il reste cependant \`a v\'erifier que la s\'erie converge en norme $L^2$. 

Une m\'ethode pour cela est d'utiliser les r\'esultats tr\`es g\'en\'eraux de la th\'eorie des ``op\'erateurs d'ar\^etes'' ou op\'erateurs sur des vari\'et\'es \`a bord d\'eg\'en\'er\'e, dont $L$ est un exemple typique. L'article \cite{Mazzeo} traite en particulier des op\'erateurs elliptiques, qui est le cas qui nous int\'eresse ici.

Une autre m\'ethode consiste \`a majorer directement les termes apparaissant dans la s\'erie, en utilisant les majorations des solutions \'el\'ementaires d'exposant positif donn\'ees dans \cite{These}.
Comme l'op\'erateur $L$ est elliptique, la $1$-forme $\eta$ est $C^\infty$ au voisinage du lieu 
singulier. En particulier, comme $\Sigma_a$ est compacte, toutes les d\'eriv\'ees de $\eta$ sont de norme 
$L^2$ finie sur $\Sigma_a$. On en d\'eduit que pour tout polyn\^ome (ou fonction major\'ee par un 
polyn\^ome) \`a deux variables $P$, la s\'erie
$$\sum_{j,k} |P(p_j\gamma,\lambda_j) f_j(a)|^2$$
converge, et qu'il en est de m\^eme en rempla\c{c}ant $f_j(a)$ par $g_j(a)$, $\omega_j(a)$ ou 
$\varpi_j(a)$. Ce r\'esultat, combin\'e aux estimations de \cite{These}, lemme 2.1.3, montre directement que la s\'erie donnant l'expression de $\n_{e_r} d\eta$
converge en norme $L^2$ sur $U_a$. 

\medskip

R\'ecapitulons : si $\eta$ est une solution de l'\'equation $L\eta=0$ au
voisinage de $\Sigma$, telle que $\eta$ et $\n \eta$ soient dans $L^2$, on
a vu que seuls les solutions \'el\'ementaires ayant un exposant dominant 
sup\'erieur ou \'egal \`a $0$ apparaissent dans la d\'ecomposition de $\eta$. 
Par cons\'equent les termes apparaissant dans la d\'ecomposition de $\n d\eta$ ont tous
des exposants dominants strictement sup\'erieurs \`a $-1$.
Ce fait, et le caract\`ere $C^\infty$ de $\eta$ pr\`es du lieu singulier, suffit \`a prouver que $\n d \eta$ est  dans $L^2$.
\end{proof}

\bigskip

A partir de ce lemme, on va pouvoir passer de l'\'etude des solutions de
l'\'equation $L\eta=0$ au voisinage du lieu singulier \`a celle des solutions de
l'\'equation $L\eta=f$ sur $M$ enti\`ere. Le th\'eor\`eme suivant montre que quand les angles coniques sont inf\'erieurs \`a $2\pi$, on peut contr\^oler le comportement de certaines combinaisons des d\'eriv\'ees premi\`eres et secondes des solutions de l'\'equation de normalisation ; ce contr\^ole sera en fait suffisant pour d\'emontrer la rigidit\'e infinit\'esimale.

\begin{Thm}\label{nablad} Soit $M$ une c\^one-vari\'et\'e hyperbolique dont tous les
angles coniques sont strictement inf\'erieurs \`a $2\pi$. Soit $\phi$
une section de $L^2(T^*M)$. Alors il existe une unique section
$\eta$ de $L^2(T^*M)$, solution de l'\'equation $L\eta =\phi$,
telle que $\eta$, $\n \eta$, $d\delta \eta$, et $\n d\eta$
(au sens des distributions) soient dans $L^2$.
\end{Thm}

\begin{proof} On sait depuis le d\'ebut de la section \ref{LapDeb} (p. \pageref{LapDeb} que l'on peut
r\'esoudre de fa\c{c}on unique l'\'equation $L\eta = \phi$ avec $\eta$, $\n \eta$ et 
$\na \n \eta$ dans $L^2$.
Le seul point qui reste \`a montrer est que $\n d\eta$ est aussi
$L^2$ ; en effet, si $\n d \eta$ est dans $L^2$, alors automatiquement $\delta d \eta$ est dans $L^2$, et donc c'est aussi le cas pour $d \delta \eta = \na \n \eta - \delta d \eta - (n-1) \eta$, voir la formule de Weitzenb\"ock \eqref{0015}.

\bigskip

Les formes $C^\infty$ \`a support compact \'etant dense dans $L^2$,
on peut trouver une suite $(\phi_k)$ de 1-formes $C^\infty$ \`a support
  compact telle que $\phi_k \to \phi$ dans $L^2$ quand $k \to
  \infty$. Soit $(\eta_k)$ la suite d'\'el\'ements de $D$ telle que pour tout
  entier $k$, $L \eta_k =\phi_k$. On applique alors le th\'eor\`eme \ref{1+*}
  (avec, \`a un facteur pr\`es,
  $A=\n$ et $A^* = \na$) : les transformations $(\na \n + (n-1)Id)^{-1}$ et
  $\n(\na \n + (n-1)Id)^{-1}$ sont continues, donc
$$\lim_{k\to \infty} \eta_k = \lim_{k\to \infty} (\na \n + (n-1)Id)^{-1}
  (\phi_k) = (\na \n + (n-1)Id)^{-1} (\phi) = \eta$$
et
$$\lim_{k\to \infty} \n \eta_k = \lim_{k\to \infty} \n ((\na \n +
  (n-1)Id)^{-1} (\phi_k)) = \n ((\na \n + (n-1)Id)^{-1} (\phi)) =
  \n \eta,$$
les limites \'etant au sens $L^2$. Comme $d\eta_k$ est la partie
antisym\'etrique de $\n \eta_k$, la suite $(d\eta_k)$ est aussi
convergente, avec $\lim_{k\to \infty} d \eta_k = d\eta$.
Maintenant, comme $\phi_k$ est \`a
  support compact, $L\eta_k$ est identiquement nul au voisinage du lieu
  singulier, et $\eta_k$ rentre donc dans le cadre de la proposition
  \ref{dvpt}. Comme $\eta_k$ appartient \`a $D(=D')$, $\eta_k$ ainsi que
  $d\eta_k$ sont dans $L^2$, et on a vu au lemme \ref{nnd} qu'alors $\n d
  \eta_k \in L^2$. On va maintenant montrer que $(\n d
  \eta_k)$, suite de sections du fibr\'e $T^*M\otimes \Lambda^2M$, est {\em
  born\'ee} dans $L^2$.

Pour cela, on consid\`ere $\xi$, section $C^{\infty}$ \`a support compact
de $T^*M\otimes \Lambda^2M$ (``section test''), et on s'int\'eresse au produit
scalaire $\<\n d \eta_k, \xi \>$. Le but est d'arriver \`a monter que $$|\<\n d
\eta_k, \xi \>| \leq M ||\xi||,$$ o\`u $M$ ne d\'epend pas de $k$ ni de $\xi$.

\bigskip

La restriction de la d\'eriv\'ee covariante \`a $\Omega^2M$ nous donne
un op\'erateur (non born\'e) $\n : L^2(\Lambda^2M) \to
L^2(T^*M\otimes\Lambda^2M)$; son adjoint est la restriction de
$\na$ \`a $L^2(T^*M\otimes\Lambda^2M)$, et les r\'esultats de la section
\ref{secIPP} s'appliquent. Maintenant, en utilisant la d\'efinition
de l'adjoint d'un op\'erateur, on a l'\'egalit\'e $\ker \na = (\im
\n)^\perp$. Le th\'eor\`eme \ref{imferme} nous garantie que l'image de $\n$ est ferm\'ee,
et on a donc la d\'ecomposition orthogonale suivante :
 $$L^2(T^*M\otimes \Lambda^2M) = \ker \na \oplus \im \n.$$
On peut donc \'ecrire $\xi = k + \n \zeta$ dans cette
d\'ecomposition, et d'apr\`es le corollaire \ref{corimferme} on peut m\^eme choisir $\zeta$ de telle 
sorte que 
$$||\zeta|| \leq c\,||\n \zeta|| \leq c\,||\xi||$$
pour une constante $c$ donn\'ee ne d\'ependant pas de $\xi$.

Retournons au produit scalaire :
\begin{eqnarray*}
\<\n d \eta_k, \xi \> & = & \<\n d \eta_k, \n \zeta + k\> \\
& = & \<\n d \eta_k, \n \zeta\>.
\end{eqnarray*}
Pour pouvoir faire une int\'egration par parties, il faut v\'erifier
que tous les termes impliqu\'es sont $L^2$. On sait d\'ej\`a que
$\zeta$, $\n \zeta$, $\n d \eta_k$ le sont, reste \`a montrer que c'est aussi le
cas de $\na \n d \eta_k$. Pour cela on utilise la formule de Weitzenb\"ock
suivante, valable pour une m\'etrique hyperbolique, qui est un
analogue de la formule pour les $1$-formes que l'on a d\'ej\`a
utilis\'ee \`a plusieurs reprises (voir \cite{Besse} \S 1.I) :
\begin{eqnarray} \forall \omega \in \Omega^2M,\
\na \n \omega = \Delta \omega + 2 (n-2) \omega. \end{eqnarray}

En l'appliquant \`a $\eta_k$, on trouve
$$\na \n d \eta_k = \Delta d
\eta_k + 2(n-2) d \eta_k = d \Delta \eta_k + 2(n-2) d
\eta_k,$$
car $d$ et $\Delta = d\delta +\delta d$ commutent.
D'autre part
$$\Delta \eta_k = L\eta_k - 2(n-1)\eta_k =
\phi_k - 2 (n-1) \eta_k.$$
 Finalement,
$$\na \n d \eta_k = d \phi_k - 2 d\eta_k.$$

Comme $\phi_k$ est \`a support compact et que $\eta_k \in D$, les
formes $d\phi_k$ et $d\eta_k$ sont $L^2$, donc $\na \n d
\eta_k$ est $L^2$, donc on peut donc int\'egrer par parties
(th\'eor\`eme \ref{ipp}) :
\begin{eqnarray*}
\<\n d \eta_k, \xi \> & = & \<\n d \eta_k, \n \zeta\>\\ & = &
\<\na \n d \eta_k, \zeta\>\\ & = & \<d \phi_k - 2 d\eta_k,
\zeta\>.
\end{eqnarray*}

Comme $\n \zeta$ est $L^2$, $\delta \zeta = - \tr_g \n \zeta$ est
aussi $L^2$, on a m\^eme $||\delta \zeta|| \leq \sqrt{n} ||\n
\zeta||$. D'autre part $\phi_k$, $d\phi_k$ et $\zeta$ sont $L^2$,
on peut encore int\'egrer par parties :
\begin{eqnarray*}
\<\n d \eta_k, \xi \> & = & \<d \phi_k - 2 d\eta_k, \zeta\>\\
& = & \<\phi_k,\delta \zeta\>- 2\<d\eta_k,\zeta\>.
\end{eqnarray*}

Pour finir on majore avec Cauchy-Schwarz:
\begin{eqnarray*}
|\<\n d \eta_k, \xi \>| & \leq & ||\phi_k||\,||\delta \zeta|| +
2 ||d\eta_k||\,||\zeta|| \\ & \leq & (\sqrt{n} ||\phi_k|| + 2 c
||d\eta_k||)\,||\n \zeta||\\ & \leq & M ||\xi||
\end{eqnarray*}
car les suites $(\phi_k)$ et $(d\eta_k)$ sont convergentes, donc born\'ees,
dans $L^2$. Cette majoration, valable pour toute section test $\xi$, implique
directement que la suite $(\n d\eta_k)$ est born\'ee dans $L^2$.

\bigskip

Par cons\'equent, on peut extraire une sous-suite, encore not\'ee $(\n
d\eta_k)$, qui converge faiblement vers une limite $l\in L^2$ : c'est-\`a-dire
que quel que soit $\xi \in L^2(T^*M\otimes \Lambda^2M)$,
$$\lim_{k \to \infty}\<\n d\eta_k, \xi\> = \<l,\xi\>.$$ Mais alors, si $\xi$
est $C^\infty$ \`a
support compact,
$$\<\n d\eta_k, \xi\> = \<d\eta_k, \na \xi\>,$$ et
$$\lim_{k \to \infty} \<d\eta_k, \na \xi\> = \<d\eta,\na \xi\>$$ car
$(d\eta_k)$ converge dans $L^2$ vers $d\eta$. Par cons\'equent, on a
$$\<d\eta,\na \xi\> = \<l,\xi\>$$ pour tout $\xi \in C^\infty_0$, ce qui
signifie exactement que $$l=\n_{max} d\eta = \n d\eta,$$
 et par suite $\n d\eta$ appartient \`a $L^2$.
\end{proof}

\medskip

Notons que si en plus $\phi$ est $C^\infty$, alors par r\'egularit\'e
elliptique la solution $\eta$ ci-dessus est aussi de classe
$C^\infty$.

\bigskip

\subsection{Cas des angles coniques inf\'erieurs \`a $\pi$ }\label{pi}
 
Si on suppose que tous les angles coniques sont inf\'erieurs \`a $\pi$, il est alors possible d'avoir un meilleur contr\^ole sur le comportement des solutions de l'\'equation de normalisation $L\eta = \phi$, o\`u $\phi \in L^2$. Plus pr\'ecis\'ement, on a le th\'eor\`eme suivant :

\medskip

\begin{Thm}\label{nnu}
Soit $M$ une c\^one-vari\'et\'e hyperbolique dont tous les
angles coniques sont strictement inf\'erieurs \`a $\pi$. Alors l'op\'erateur 
$$L = \na \n + (n-1)Id\ :\ L^{2,2}(T^*M) \to L^2(T^*M)$$
est un isomorphisme.
\end{Thm}

\medskip

On retrouve donc quand les angles sont suffisamment petits une propri\'et\'e toujours valable sur 
une vari\'et\'e compacte. Une application imm\'ediate est que si l'on part d'une d\'eformation infinit\'esimale $L^{1,2}$, par exemple pr\'eservant les angles, alors la d\'eformation normalis\'ee correspondante sera aussi dans $L^{1,2}$.

La d\'emonstration est en fait compl\`etement analogue \`a celle du th\'eor\`eme \ref{nablad} ; on ne redonnera donc pas tous les d\'etails. Elle repose elle aussi sur un lemme :

\begin{Lem}\label{nnul2} Soit $M$ une c\^one-vari\'et\'e hyperbolique dont tous les
angles coniques sont strictement inf\'erieurs \`a $\pi$. Soit $\eta$
une $1$-forme telle que $L\eta$ soit \'egal \`a $0$ au voisinage du lieu singulier et que $\eta$ et $\n \eta$ 
soient dans $L^2$. Alors $\n \n \eta$ appartient \`a $L^2$.
\end{Lem}

\bigskip

\begin{proof}[D\'emonstration du lemme] La preuve est en tout point similaire \`a celle du lemme \ref{nnd}. Comme on sait d\'ej\`a que $\n d \eta$ est dans $L^2$, il suffit de le montrer pour $\n \delta^*\eta$. D'apr\`es la proposition \ref{derl2}, on a juste \`a montrer que $\na \n (\delta^* \eta)$ et $\n_{e_r} \delta^*\eta$ sont dans $L^2$.

Pour $\na \n (\delta^*\eta)$, cela r\'esulte d'une relation de commutation entre $\na \n$ et
$\delta^*$. Soit $\sigma$ une $1$-forme (lisse) sur $M$. On sait que la d\'eformation
$\delta^* \sigma$ est triviale, c'est-\`a-dire que l'on a toujours $E'(\delta^*
\sigma)=0$, soit    
$$\na \n (\delta^* \sigma) - 2 \RO (\delta^*\sigma) - 2 \delta^*
\beta \delta^* \sigma = 0,$$
voir \eqref{0019}. La m\'etrique \'etant hyperbolique, on peut utiliser les relations \eqref{0020} 
et \eqref{0022} pour simplifier cette expression. On obtient alors la relation suivante :
\begin{eqnarray}\na \n (\delta^* \sigma) =  2 \delta^*\sigma + 2 (\delta \sigma) g +
\delta^*(\na \n \sigma + (n-1) \sigma). \label{0001}\end{eqnarray}
En appliquant cette formule \`a $\eta$, on obtient au voisinage du lieu singulier
$$\na \n (\delta^* \eta) =  2 \delta^*\eta + 2 (\delta \eta) g$$
car $L\eta = \na \n \eta + (n-1) \eta = 0$ pr\`es du lieu singulier. Comme $\eta$ est $L^{1,2}$, on v\'erifie facilement que le terme de droite, et donc $\na \n \delta^* \eta$, appartiennent \`a $L^2(S^2M)$.

Pour montrer que $\n_{e_r} \delta^*\eta$ est dans $L^2$, l'id\'ee est la suivante.
Le fait que $\eta$ appartient \`a 
$L^{1,2}$ impose que seules les solutions \'el\'ementaires d'exposant dominant positif apparaissent 
dans la d\'ecomposition de $\eta$ donn\'ee par la proposition \ref{dvpt}. Et comme tous les angles coniques sont suppos\'es plus petits que 
$\pi$, les exposants dominants positifs sont soit nuls, soit plus grands que $1$. Cela 
donne suffisamment de contr\^ole sur les solutions \'el\'ementaires pour montrer que tous les termes 
apparaissant dans la d\'ecomposition de $\n \n \eta$ sont dans $L^2$.
On utilise ensuite 
le fait que $P$ est elliptique, donc que $u$ est $C^\infty$ pr\`es du lieu singulier, pour assurer comme dans la preuve du lemme \ref{nnd}
la convergence en norme $L^2$ de la d\'ecomposition en s\'erie sur un voisinage du lieu singulier de 
$\n_{e_r} \delta^* \eta$, \`a l'aide des r\'esultats de \cite{Mazzeo} ou de \cite{These}. 
\end{proof}

\bigskip

Pour passer de $L\eta = 0$ au voisinage du lieu singulier \`a $L\eta = \phi$ sur $M$, l'id\'ee est encore de prendre une suite de $1$-formes $\phi_k$, $C^\infty$ \`a support compact, convergeant vers $\phi$ en norme $L^2$. On prend alors une suite $\eta_k$ telle que $L\eta_k = \phi_k$ et $\eta_k \in L^{1,2}$. On sait d\'ej\`a que $\eta_k$ et $\n \eta_k$ convergent en norme $L^2$ vers $\eta$ et $\n \eta$ respectivement, et que $\n d \eta_k$ converge faiblement vers $\n d \eta \in L^2$. Il ne reste plus qu'\`a consid\'erer $\n \delta^* \eta_k$. Pour montrer que cette suite est born\'ee, on regarde encore le produit scalaire contre une section test (i.e. $C^{\infty}$ \`a support compact) $\xi$ de $T^*M\otimes S^2M$.
Gr\^ace \`a plusieurs int\'egrations par partie et \`a la relation de commutation \eqref{0001}, on aboutit \`a l'existence d'une constante $M$ telle que $|\<\n \delta^* \eta_k ,\xi\>| \leq M ||\xi||$ pour toute section test $\xi$, ce qui justifie que la suite $\n \delta^* \eta_k$ est born\'ee. On finit en montrant qu'elle converge faiblement vers une limite, qui est alors \'egale \`a $\n \delta^*\eta$, ce qui termine de d\'emontrer que $\n \delta^*\eta \in L^2$.

\bigskip

{\em Remarque :} Bien qu'on ne le montre pas ici, il est
int\'eressant de noter que les th\'eor\`emes \ref{nablad} et \ref{nnu} cessent d'\^etre vrai d\`es qu'un angle conique est plus grand
que respectivement $2\pi$ ou $\pi$. Il
devient alors beaucoup plus difficile de trouver un ``bon'' domaine
pour r\'esoudre l'\'equation de normalisation, cf \cite{HK2} et
\cite{BB} pour des r\'esultats dans cette direction.

\bigskip

\section{D\'eformations Einstein infinit\'esimales}\label{demo}

\subsection{Rigidit\'e infinit\'esimale des c\^one-vari\'et\'es}\label{rigidite}

Nous avons maintenant en main tous les outils pour montrer le th\'eor\`eme suivant~:

\begin{Thm}\label{thmppal} Soit $M$ une c\^one-vari\'et\'e hyperbolique dont tous
  les angles coniques sont strictement inf\'erieurs \`a $2\pi$.
Soit $h_0$ une d\'eformation Einstein infinit\'esimale (i.e. v\'erifiant
  l'\'equation $E'_g(h_0)=0$) telle que $h_0$ et $\n h_0$ soient dans
  $L^2$. Alors la d\'eformation $h_0$ est triviale, i.e. il existe une forme
  $\eta \in \Omega^1 M$ telle que $h_0 = \delta^* \eta$.
\end{Thm}

Dans toute cette sous-section on supposera donc que les angles coniques sont
toujours inf\'erieurs \`a $2\pi$.

\bigskip

\begin{proof}
La premi\`ere \'etape de la d\'emonstration consiste \`a normaliser $h_0$,
c'est-\`a-dire \`a chercher $\eta$ tel que $h = h_0 - \delta^*
\eta$ v\'erifie la condition de jauge $\beta(h)=0$, ce qui revient
\`a r\'esoudre l'\'equation $\beta \circ \delta^* u = \beta h_0$.
Comme $\n h_0$ est dans $L^2$, $\beta h_0$ l'est aussi, et d'apr\`es
le th\'eor\`eme \ref{nablad} cette \'equation admet une unique solution
$\eta$ telle que $\eta$, $\n \eta$, $d \delta \eta$ et $\n
d\eta$ soient dans $L^2$. On pose $h = h_0 - \delta^*\eta$.
Notons que l'on a perdu des informations en normalisant : en
effet, rien ne garantit que la d\'eformation normalis\'ee $h$ v\'erifie
encore $\n h \in L^2$, puisqu'on ne conna\^it rien pour l'instant sur $\n
\delta^* \eta$.

La d\'eformation $h$ v\'erifie alors :
$$\begin{cases} \na \n h - 2 \mathring{R}h = 0 \\ \delta h + d \tr h = 0
\end{cases}$$

En prenant la trace par rapport \`a $g$ de la premi\`ere \'equation, on
obtient $$\Delta(\tr h) + 2(n-1) \tr h =0,$$ ce qui incite \`a
int\'egrer par parties, mais pour le faire il faut d'abord v\'erifier
que les termes impliqu\'es sont $L^2$, avant de pouvoir appliquer le
th\'eor\`eme \ref{Cheeg}. Comme $h_0$ et $\delta^* \eta$ sont $L^2$,
$h$ est bien $L^2$, donc $\tr h$ aussi, et donc $\Delta \tr h$
aussi. Maintenant, $$d \tr h = d \tr h_0 + d \tr \delta^* \eta =
d \tr h_0 - d\delta \eta,$$
 donc $d \tr h$ est $L^2$ ($d \tr h_0$ est $L^2$ car $\n h_0$ l'est). Par
suite, on trouve en int\'egrant contre $\tr h$ :
\begin{eqnarray*} 0 & = & \< \tr h, \Delta (\tr h) + 2 (n-1)
\tr h \> \\
 & = & ||d\tr h||^2 + 2(n-1)||\tr h||^2
\end{eqnarray*}
et donc $\tr h = 0$, ce qui, avec $\beta(h)=0$, implique aussi
$\delta h = 0$. Finalement, on a
 $$\begin{cases} \na \n h - 2
\mathring{R}h = 0 \\ \delta h =0 \\ \tr h = 0
\end{cases}$$

\bigskip

La deuxi\`eme \'etape de la d\'emonstration consiste \`a utiliser une
autre formule de Weitzenb\"ock (cf \cite{Besse}, \S 12.69). Un
$2$-tenseur peut toujours se voir comme une $1$-forme \`a valeur
dans le fibr\'e cotangent $T^*M$. Ce fibr\'e \'etant muni de la
connexion de Levi-Civit\`a $\n$, on note $d^\n$ la diff\'erentielle
ext\'erieure associ\'ee sur les formes \`a valeurs dans $T^*M$.
L'op\'erateur adjoint est la codiff\'erentielle not\'ee $\delta^\n$.
Notons que si $u$ est une $0$-forme \`a valeurs dans $T^*M$
(c'est-\`a-dire une $1$-forme usuelle), alors $d^\n u = \n
u$; de m\^eme pour une $1$-forme \`a valeurs dans $T^*M$,
$\delta^\n h = \na h$. On a alors la formule suivante, valable
pour tout $2$-tenseur sym\'etrique :
\begin{eqnarray} \na \n h = (\delta^\n d^\n + d^\n \delta^\n) h + \mathring{R}h - h\circ
 ric. \label{0004} \end{eqnarray}
Pour une m\'etrique hyperbolique, cela se simplifie en
$$\na \n h = (\delta^\n d^\n + d^\n \delta^\n) h + n h - (\tr h)g.$$
En combinant avec ce qui pr\'ec\`ede, on obtient
$$\begin{cases}  \delta^\n d^\n h + (n-2) h = 0 \\ \delta h =0 \\ \tr h = 0
\end{cases}$$

Pour conclure, ``il suffit'' d'une int\'egration par parties contre $h$. Comme $h$
est dans
$L^2$, $\delta^\n d^\n h$ est aussi dans $L^2$; si $\n h$, ou m\^eme seulement
$\n_{e_r} h$, \'etait $L^2$ on pourrait conclure en utilisant une m\'ethode
analogue \`a celle employ\'ee dans la d\'emonstration du th\'eor\`eme
\ref{ipp}. Malheureusement on ne sait rien sur le caract\`ere $L^2$ ou non de
$\n \delta^* \eta$. On va donc devoir contourner cette difficult\'e pour
montrer qu'on a bien $\<\delta^\n d^\n h, h\> = ||d^\n h||^2$.

Avant toutes choses, il faut montrer que $d^\n h$ est bien $L^2$. Comme $\n
h_0$ est $L^2$, $d^\n h_0$ est $L^2$; il ne reste qu'\`a regarder $d^\n
\delta^*\eta$. Or
$$\delta^*\eta = \n \eta - \frac{1}{2}d\eta = d^\n
\eta - \frac{1}{2}d\eta,$$
 donc
$$d^\n \delta^*\eta = (d^\n)^2\eta -
\frac{1}{2} d^\n d\eta.$$
L'op\'erateur $(d^\n)^2$ est bien connu, ce n'est rien d'autre
que l'oppos\'e de la courbure, i.e.
$$(d^\n)^2\eta (x,y) = - R(x,y)\eta = \n_x \n_y \eta - \n_y \n_x \eta
 - \n_{[x,y]}\eta.$$
 C'est un op\'erateur born\'e, c'est-\`a-dire continu, pour les normes $L^2$; par
 cons\'equent
$(d^\n)^2\eta$ est $L^2$. Il ne nous reste donc que la terme $d^\n d\eta$;
or le th\'eor\`eme \ref{nablad} nous garantit que $\n d\eta$, et donc $d^\n
d\eta$, sont bien $L^2$.

\bigskip

Le tenseur $d^\n h$ est donc bien dans $L^2$. Malheureusement,
on n'a pas d'analogue du r\'esultat de Cheeger (th\'eor\`eme \ref{Cheeg}) pour les
formes \`a valeurs dans un fibr\'e, du fait que $d^\n\circ d^\n$ ne s'annule pas
n\'ecessairement, \`a la diff\'erence de $d\circ d$. Cependant, en \'ecrivant
$$h = h_0 - \delta^*\eta = h_0 + \frac{1}{2}d\eta - d^\n \eta,$$
 on a
\begin{eqnarray}\<h,\delta^\n d^\n h\> = \<h_0 + \frac{1}{2}d\eta, \delta^\n d^\n h\> -
\<d^\n \eta, \delta^\n d^\n h\>.\label{0013}\end{eqnarray}

Le th\'eor\`eme \ref{nablad} nous assure que $\n(h_0 + \frac{1}{2}d\eta)$ est
dans $L^2$. Ceci
nous permet de montrer, exactement de la m\^eme fa\c{c}on que dans la d\'emonstration
du th\'eor\`eme \ref{ipp}, qu'on a bien
\begin{eqnarray}\<h_0 + \frac{1}{2}d\eta, \delta^\n d^\n h\> = \<d^\n(h_0 +
\frac{1}{2}d\eta), d^\n h\>.\label{0014}\end{eqnarray}

Pour le terme qui reste, comme $\eta$ et $\n \eta$ sont $L^2$, on peut
trouver d'apr\`es le corollaire \ref{approx} une suite $(\eta_k)$, $C^\infty$
\`a support compact, telle que $\lim_{k\to\infty} \eta_k = \eta$ et
$\lim_{k\to\infty} \n \eta_k = \n \eta$. On a alors
$$\lim_{k\to \infty} \<d^\n \eta_k, \delta^\n d^\n h\> =
\<d^\n \eta, \delta^\n d^\n h\>.$$
 On peut faire l'int\'egration par parties avec
$\eta_k$ :
$$\<d^\n \eta_k, \delta^\n d^\n h\> = \<(d^\n)^2 \eta_k, d^\n
h\>.$$
 Mais comme $(d^\n)^2$ est continue, on a
$$\lim_{k\to \infty} (d^\n)^2 \eta_k = (d^\n)^2 \eta,$$
 et donc
$$\lim_{k\to \infty} \<(d^\n)^2 \eta_k, d^\n h\> = {\<(d^\n)^2 \eta, d^\n
 h\>}.$$
 On en d\'eduit que
$$\<d^\n \eta, \delta^\n d^\n h\> = \<(d^\n)^2 \eta, d^\n h\>,$$
et avec \eqref{0013} et \eqref{0014} on a \'etabli l'\'egalit\'e
$$\<h,\delta^\n d^\n h\> = ||d^\n h||^2.$$

Par cons\'equent, comme $\delta^\n d^\n h + (n-2) h =0$, on a
\begin{eqnarray*} 0 & = & \<h,\delta^\n d^\n h + (n-2) h\> \\
& = & ||d^\n h||^2 + (n-2) ||h||^2
\end{eqnarray*}
et donc le tenseur $h$ est identiquement nul. Par suite $h_0 = \delta^*
\eta$, la d\'eformation est triviale.
\end{proof}

\bigskip

\begin{Cor}[Rigidit\'e infinit\'esimale] Soit $M$ une c\^one-vari\'et\'e hyperbolique dont tous
  les angles coniques sont strictement inf\'erieurs \`a $2\pi$. Alors $M$ est
  infinit\'esimalement rigide parmi les c\^ones-vari\'et\'es Einstein \`a angles
  coniques fix\'es.
\end{Cor}

\begin{proof} En effet, on a vu que toute d\'eformation infinit\'esimale de la
  structure de c\^one-vari\'et\'e pr\'eservant les angles pouvait se mettre sous la
  forme d'un $2$-tenseur sym\'etrique $h_0$ appartenant \`a $L^2$, dont la d\'eriv\'ee
  covariante $\n h_0$ est aussi dans $L^2$. On peut alors appliquer le th\'eor\`eme
  ci-dessus pour montrer que toutes les d\'eformations Einstein de ce type sont
  triviales.
\end{proof}

\subsection{Constructions de d\'eformations Einstein modifiant les angles}\label{constr}

On vient de voir que sur une c\^one-vari\'et\'e hyperbolique dont tous les angles coniques sont inf\'erieurs \`a $2\pi$, les d\'eformations Einstein infinit\'esimales pr\'eservant les angles \'etaient triviales. A l'inverse, on va maintenant s'int\'eresser aux d\'eformations Einstein infinit\'esimales r\'ealisant (au premier ordre) une variation donn\'ee des angles coniques. Pour cela, on aura besoin du th\'eor\`eme \ref{nnu} ; en cons\'equence on se restreindra au cas o\`u les angles coniques sont tous inf\'erieurs \`a $\pi$. On aura aussi besoin du th\'eor\`eme \ref{normgen}, ou plut\^ot de sa variante donn\'ee par la proposition 
\ref{minmax}, et dont la d\'emonstration suit :

\begin{nThm}[\ref{normgen}, \ref{minmax}]
Si $M$ est une c\^one-vari\'et\'e hyperbolique, alors les extensions minimales et maximales des op\'erateurs 
$\beta$ et $\delta^*$ sont \'egales :  
$\beta_{max}= \beta_{min}$, $\delta^*_{max}=\delta^*_{min}$. On note ces op\'erateurs simplement
$\beta$ et $\delta^*$. On a alors la d\'ecomposition en somme directe  
$L^2(S^2M) = \ker \beta \oplus \im \delta^*$.
\end{nThm}

\begin{proof} \label{demminmax} Le seul point laiss\'e de c\^ot\'e \`a la section \ref{premres} est l'\'egalit\'e des extensions minimales et maximales dans le cas o\`u la m\'etrique est hyperbolique.Prenons $\eta \in D(\delta^*_{max})$. D'apr\`es le th\'eor\`eme \ref{normgen}, il existe $k 
\in \ker \beta_{max}$ et $\eta' \in D(\delta^*_{min})$ tels que $\delta^*_{max} \eta = k + 
\delta^*_{min} \eta' = k$. La section $u=\eta - \eta'$ v\'erifie $\beta_{max}(\delta^*_{max} 
(u))=0$, et donc (au sens des distributions)
$$\na \n u + (n-1) u = 0.$$
La section $u$ admet alors un d\'eveloppement du type donn\'e \`a la proposition \ref{dvpt}. Le fait 
que $u$ et $\delta^*_{max} u$ soient dans $L^2$ impose que seuls les termes ayant un exposant 
dominant sup\'erieur ou \'egal \`a $0$ apparaissent dans cette d\'ecomposition. On peut alors appliquer 
le raisonnement utilis\'e dans les preuves des lemmes \ref{nnd} et \ref{nnul2} pour montrer que  
$\n u$ est dans $L^2$, et donc que $u$ appartient \`a $L^{1,2}$. Or le noyau de $\na \n + (n-1) 
Id$ dans $L^{1,2}$ est r\'eduit \`a $\{0\}$ ; par cons\'equent $u=0$, et $\eta = \eta'$ appartient \`a 
$D(\delta^*_{min})$. Cela montre que $D(\delta^*_{max})=D(\delta^*_{min})$, et donc que 
$\delta^*_{max} = \delta^*_{min}$

D'autre part on a vu (proposition \ref{prop} et la remarque p.\pageref{rem}) que 
$D(\beta^t_{min}) = D(\delta^*_{min})$ et que $D(\beta^t_{max}) = D(\delta^*_{max})$.
Alors $D(\beta^t_{min}) =D(\beta^t_{max})$, soit $\beta^t_{min} = \beta^t_{max}$. 
En passant \`a l'adjoint, on trouve directement que $\beta_{max} = \beta_{min}.$
\end{proof}

\bigskip

Le principe de construction d'une d\'eformation infinit\'esimale Einstein est assez simple. On part d'une d\'eformation infinit\'esimale (non Einstein) $h_0$ r\'ealisant la variation voulue des angles coniques, et on cherche \`a la ``corriger'' par un $2$-tenseur sym\'etrique $h$ dans $L^{1,2}$ (donc ne modifiant pas les angles coniques), de telle sorte que la nouvelle d\'eformation $h_0-h$ soit Einstein, c'est-\`a-dire que $E'(h_0-h)=0$. Cela revient donc \`a r\'esoudre l'\'equation 
$$E'(h) = E'(h_0),$$
avec $h$ dans $L^{1,2}$. Le probl\`eme est qu'il sagit d'une ``mauvaise'' \'equation : l'op\'erateur Einstein lin\'earis\'e $E' = \na\n - 2 \RO - 2 \delta^* \circ \beta$ est fortement d\'eg\'en\'er\'e, au sens o\`u son noyau, qui contient l'espace des d\'eformations triviales $\im \delta^*$, est de dimension infini. L'id\'ee est alors de chercher \`a r\'esoudre \`a la place le syst\`eme suivant :
\begin{equation} \begin{cases} \na \n h -2 \RO h = E'(h_0) \\ \beta(h) = 0 \end{cases} \label{0032} \end{equation}
L'op\'erateur $P=\na \n - 2 \RO$ est nettement plus sympathique : il est elliptique, sym\'etrique, coercif, non d\'eg\'en\'er\'e. En particulier, on peut montrer sans trop de difficult\'es que pour toute section $\phi$ de $L^2(S^2M)$, il existe un unique $h$ dans $L^{1,2}(S^2M)$ tel que $Ph=\phi$. Pour pouvoir r\'esoudre le syst\`eme \eqref{0032}, on va consid\'erer la restriction de $P$ au sous-espace vectoriel ferm\'e $\ker \beta$, dont les propri\'et\'es qui nous int\'eressent sont donn\'ees par la proposition suivante :

\bigskip

\begin{Prop} Soit $D(P) = \{ h \in \ker \beta\ |\ \n h \in L^2 {\rm \ et\ }
  \na \n h \in L^2 \}$. L'op\'erateur 
\begin{eqnarray*}
 P : D(P) \subset \ker \beta & \to & \ker \beta\\
h & \mapsto & \na \n h - 2 h + 2 (\tr h) g
\end{eqnarray*}
 est auto-adjoint et positif, et donc inversible.
\end{Prop} 

\begin{proof} V\'erifions avant toute autre chose que l'image de $P$ est bien
incluse dans $\ker \beta$. On sait que pour toute section lisse $h$ de $S^2M$, 
\begin{equation} \beta (E'(h)) = \beta (\na \n h - 2 \RO h - 2 \delta^* \beta h) = 0. \label{0033}\end{equation}
Donc $\beta (\na \n h - 2 \RO h) = 2 \beta \delta^* \beta h$, et si $\beta h =0$, on a 
imm\'ediatement 
$\beta (\na \n h - 2 \RO h) = 0.$
On va montrer que c'est aussi le cas (au moins au sens des distributions) si $h$ appartient \`a 
$D(P)$.

Soit $\xi$ une section test de $S^2M$, $C^\infty$ \`a support compact, et $h$ un \'el\'ement de 
$D(P)$. Alors 
$$\<Ph, \beta^t \xi\> = \<h, (\na \n - 2 \RO)(\beta^t \xi)\> = \<h, 2 \beta^t \delta \beta^t \xi 
\> = 0$$
car $h$ appartient \`a $\ker \beta$. Cela montre que $Ph \in \ker \beta$ pour tout $h \in D(P).$

\medskip

On aura ensuite besoin du lemme suivant :

\begin{Lem} Le domaine de l'op\'erateur $\n$, restreint \`a $\ker \beta$, est
  dense dans $\ker \beta$ pour la norme $L^2$ : $\overline{D(\n|_{\ker
  \beta})} = \ker \beta$. 
\end{Lem}

\begin{proof}[D\'emonstration du lemme] Ce r\'esultat d\'ecoule des deux th\'eor\`emes
  \ref{normgen} et \ref{nnu}. On utilise la 
  d\'e\-com\-po\-si\-tion $L^2 = \ker \beta \oplus \im \delta^*$ et le fait que
la projection $p_1$ sur le premier facteur est continue. Comme $C^\infty_0$ est
  un sous-espace dense de $L^2$, son image $p_1(C^\infty_0)$ est dense dans
  $\ker \beta$. 

Maintenant soit $k \in p_1(C^\infty_0)$. Il existe $\phi \in C^\infty_0$ et 
$\eta \in D(\delta^*)$ tels que $k = p_1(\phi)$ et $\phi = k + \delta^*
\eta$. On alors $\beta(\phi) = \beta (\delta^* \eta)$. Or $\beta(\phi)$
est \`a support compact, donc dans $L^2$. D'apr\`es le th\'eor\`eme \ref{nnu}, $\n \n \eta$ est
$L^2$, et par cons\'equent $\delta^* \eta \in D(\n)$. Comme $\phi$ est
$C^\infty$ \`a support compact, elle est aussi dans $D(\n)$, et donc $k=\phi -
\delta^* \eta \in D(\n)$. Par suite $p_1(C^\infty_0) \subset D(\n)$, et
$D(\n) \cap \ker \beta = D(\n|_{\ker \beta})$ est dense dans $\ker \beta$.
\end{proof}

\bigskip

Le fait que le domaine de $\n|_{\ker \beta}$ soit dense permet de consid\'erer
son adjoint $(\n|_{\ker \beta})^*$. On note $i$ l'inclusion $\ker \beta
\hookrightarrow L^2$; on a \'evidemment $\n|_{\ker \beta} = \n
\circ i$. Son adjoint est donc $p \circ \na$, o\`u $p$ est la projection
orthogonal de $L^2$ sur $\ker \beta$. De plus on d\'eduit facilement que
l'op\'erateur $\n|_{\ker \beta}$ est ferm\'e du fait que $\n$ et $\ker \beta$ sont
ferm\'es. Cela implique que $(\n|_{\ker \beta})^* \circ \n|_{\ker \beta} =
p\circ \na \circ \n|_{\ker \beta}$ est auto-adjoint.

Maintenant, l'op\'erateur $p\circ \mathring{R}|_{\ker \beta} : h \mapsto p(h -
(\tr h) g)$ est un endomorphisme born\'e (i.e. continue) et auto-adjoint de
$\ker \beta$, donc 
$p\circ \na \circ \n|_{\ker \beta} - 2 p\circ \mathring{R}|_{\ker \beta}$ est
encore auto-adjoint. Or 
\begin{eqnarray*}
p\circ \na \circ \n|_{\ker \beta} - 2 p\circ \mathring{R}|_{\ker \beta} 
& = & p\circ (\na \circ \n|_{\ker \beta} - 2
\mathring{R}|_{\ker \beta}) \\
 & = & p\circ P\\
 & = & P
\end{eqnarray*}
 puisque l'image de $P$ est incluse dans $\ker \beta$. On a ainsi montr\'e que
 $P$ \'etait auto-adjoint.

\bigskip

La positivit\'e d\'ecoule de la formule de Weitzenb\"ock utilis\'ee pour d\'emontrer la rigidit\'e infinit\'esimale 
\ref{0004} : pour tout $2$-tenseur sym\'etrique $h$,
$$\na \n h = d^\n \delta^\n h + \delta^\n d^\n h + n h - (\tr h) g.$$
Si $h$ est $C^\infty_0$, on a 
\begin{eqnarray*}
||\n h||^2 & = & \< \na \n h, h \>\\ 
 & = & \< d^\n \delta^\n h + \delta^\n
d^\n h + n h - (\tr h) g, h\> \\
 & = & ||\delta^\n h||^2 + ||d^\n h||^2 + n||h||^2 - ||\tr h||^2 .
\end{eqnarray*}
Cette \'egalit\'e est encore vraie si $h \in D(\n)$ ; il suffit de consid\'erer une
suite $h_n$ de $2$-tenseurs sym\'etriques $C^\infty_0$ telle que pour la norme
$L^2$, $h_n$ converge vers $h$ et $\n h_n$ vers $\n h$.

\bigskip

Soit $h \in D(P)$. Comme $h$, $\n h$, et $\na \n h$ sont dans $L^2$, on peut
int\'egrer par partie :
\begin{eqnarray*}
\<P h,h\> & = & \< \na \n h - 2 h + 2 (\tr h) g , h\> \\
 & = & ||\n h||^2 - 2 ||h||^2 + 2 ||\tr h||^2 \\
 & = & ||\delta^\n h||^2 + ||d^\n h||^2 + (n-2) ||h||^2 + ||\tr h||^2, 
\end{eqnarray*}
et donc $\<P h,h\> \geq (n-2) ||h||^2$.
\end{proof}

\bigskip

\begin{Cor}
Soit $h_0$ une d\'eformation infinit\'esimale telle que $E'(h_0)$ soit dans $L^2$. Alors il existe 
un unique $2$-tenseur sym\'etrique $h \in D(P)$ tel que la d\'eformation $h_0 + h$ soit Einstein, 
i.e. $E'(h_0 + h) = 0.$
\end{Cor}

\begin{proof} Il suffit de remarquer que $E'(h_0) \in \ker \beta$ et de r\'esoudre $P h = - 
E'(h_0)$ dans $D(P)$. Et comme $h \in \ker \beta$, $E'(h) = \na \n h - 2 \RO h (- 2 \delta^* 
\beta h) = P h$.
\end{proof}

Ce corollaire n'a d'int\'er\^et que si la d\'eformation $h_0$ n'est pas $L^{1,2}$. Sinon, d'apr\`es le 
th\'eor\`eme de rigidit\'e infinit\'esimale, $h + h_0$ est une d\'eformation triviale, c'est-\`a-dire que 
$h$ est (au signe pr\`es) la composante dans $\ker \beta$ de $h_0$. Mais si $h_0$ est bien choisi, 
on peut montrer ainsi le th\'eor\`eme suivant :

\bigskip

\begin{Thm}\label{consdef}
Soit $M$ une c\^one-vari\'et\'e hyperbolique dont tous les angles coniques $\alpha_1,\ldots \alpha_p$ 
sont strictement inf\'erieurs \`a $\pi$. Soit $\dot{\alpha} = (\dot{\alpha_1}, \ldots 
\dot{\alpha_p})$ une variation donn\'ee du $p$-uplet des angles coniques. Alors il existe une 
d\'eformation Einstein infinit\'esimale normalis\'ee $h$ (i.e. telle que $E'(h) = 0$ et $\beta h =0$) 
induisant la variation des angles coniques donn\'ee.
\end{Thm} 

\bigskip

\begin{proof}
Sans perdre de g\'en\'eralit\'es, on peut se limiter au cas o\`u un seul des $\dot{\alpha_i}$ est non 
nul, ce qui revient \`a ne consid\'erer qu'une seule composante connexe $\Sigma_k$ du lieu 
singulier. 

Soit $h_0$ une d\'eformation infinit\'esimale telle que $h_0$ soit \'egal \`a $\sinh(r)^2 d\theta^2$ 
pr\`es de $\Sigma_k$, et que $h_0$ soit nul en dehors 
d'un voisinage de cette composante connexe. C'est le mod\`ele de d\'eformation conique modifiant un 
(seul) 
angle conique. C'est aussi une d\'eformation hyperbolique (donc Einstein) pr\`es du lieu singulier : 
en particulier, $E'(h_0)$ est $C^\infty$ \`a support compact. Par contre $h_0$ n'est pas 
normalis\'ee, et m\^eme pire que \c{c}a puisque $\beta h_0 \notin L^2$. 

On va commencer par normaliser localement $h_0$, c'est-\`a-dire que 
l'on cherche une $1$-forme $\eta$ telle que $\beta(h_0 - \delta^* \eta)$ soit nul pr\`es du 
lieu singulier. Si on prend $\eta$ de la forme $f(r) e^r$, alors la fonction $f$ v\'erifie 
l'\'equation diff\'erentielle (ordinaire) 
$$-f''(r) - (\frac{1}{\tanh(r)} + (n-2) \tanh(r))f'(r) + (\frac{1}{\tanh(r)^2} + (n-2) 
\tanh(r)^2+n-1)f(r) = \frac{2}{\tanh(r)}.$$
Cette \'equation diff\'erentielle est \`a singularit\'e r\'eguli\`ere (selon la terminologie de 
\cite{Wasow}) ; elle admet au voisinage de $r=0$ une solution de la forme $f(r) = - r \ln(r) + 
r^3(f_1(r) + \ln(r) f_2(r))$, o\`u $f_1$ et $f_2$ sont d\'eveloppables en s\'eries enti\`eres. On 
d\'efinit alors $\chi(r)$ comme \'etant un fonction lisse \'egale \`a $- r \ln(r) + 
r^3(f_1(r) + \ln(r) f_2(r))$ pr\`es de $\Sigma_k$, et nulle en dehors d'un certain voisinage de 
$\Sigma_k$.

On note maintenant $h_1 = h_0 - \delta^*(\chi(r) e^r)$ ; pr\`es du lieu singulier, on a encore 
$E'(h_1) =0$, et en plus $\beta h_1=0$. Au premier ordre, $h_1 = (1 + \ln(r))(dr^2 + \sinh(r)^2 
d\theta^2) + O(r^2\ln(r))$, c'est-\`a-dire que $h_1$ ressemble asymptotiquement \`a une d\'eformation 
conforme de la m\'etrique du c\^one, et $h_1$ appartient \`a $L^2$. Par contre $\n h_1$ n'est pas 
$L^2$.

\medskip

On va maintenant normaliser globalement $h_1$, c'est-\`a-dire que l'on va r\'esoudre l'\'equation 
$ \beta \delta^* \eta = \beta h_1$ sur tout $M$ et non plus seulement sur un voisinage du lieu 
singulier. Comme $\beta h_1$ est $C^\infty$ \`a support compact et que les angles coniques sont 
inf\'erieurs \`a $\pi$, d'apr\`es le lemme \ref{nnul2} l'\'equation admet une unique solution $\eta$ 
dans $L^{2,2}(T^*M)$ (et son comportement au voisinage du lieu singulier est assez bien 
compris).

On pose alors $h_2 = h_1 - \delta^* \eta$. Cette d\'eformation est dans $L^2$, mais $\n h_2$ n'est 
pas dans $L^2$ (puisque $\n \delta^* \eta \in L^2$ et $\n h_1 \notin L^2$). On a aussi $E'(h_2) 
\in C^\infty_0(S^2M)$ et $\beta h_2 = 0$.

On peut maintenant appliquer le corollaire. On construit ainsi une d\'eformation Einstein 
infinit\'esimale $h_{\dot{\alpha}}$, normalis\'ee, et qui ne diff\`ere de $h_2$ que par un \'el\'ement de 
$D(P) \subset L^{1,2}$. Cette d\'eformation $h_{\dot{\alpha}}$ modifie donc les angles coniques de 
la m\^eme fa\c{c}on que $h_2$ et que $h_0$ ; elle induit donc bien la variation voulue des angles 
coniques.
\end{proof}

\bigskip

Comme on l'a mentionn\'e dans l'introduction, ces deux r\'esultats montrent que dans un certain sens, l'espace tangent \`a une c\^one-vari\'et\'e hyperbolique parmi les structures de c\^ones-vari\'et\'es Einstein est de dimension finie, param\'etr\'e par les variations du $p$-uplet des angles coniques. On voudrait maintenant en savoir plus sur le comportement des d\'eformations Einstein infinit\'esimales que l'on vient de construire. En fait celui-ci est relativement facile \`a comprendre ; cet \'etude est l'objet de la section suivante.

\section{Comportement des d\'eformations}\label{compdef}

On a vu (th\'eor\`eme \ref{consdef}) qu'\`a toute variation du $p$-uplet des angles coniques correspondait une unique d\'eformation Einstein infinit\'esimale normalis\'ee. Cette d\'eformation est de la forme $h_1 - \delta^* \eta - h$. On conna\^it le terme $h_1$ : il ressemble 
asymptotiquement \`a un certain changement conforme de la m\'etrique du c\^one. Le terme en $\delta^* 
\eta$ est une d\'eformation triviale correspondant \`a la 
normalisation. On sait que $L\eta = 0$ pr\`es du lieu singulier ; les r\'esultats de la section \ref{vois3} d\'etaillent son comportement. Le terme $h$ est moins bien connu. On sait cependant qu'il v\'erifie la condition de jauge de Bianchi ainsi que l' \'equation $Ph=0$ au voisinage du lieu singulier. 

On peut alors proc\'eder avec l'op\'erateur $P=\na\n -2 \RO$, agissant sur les $2$-tenseurs sym\'etriques, comme on l'a fait avec l'op\'erateur $L=\na \n + (n-1)Id$, agissant sur les $1$-formes. La connaissance du comportement des solutions de l'\'equation $Pu=0$ pr\`es du lieu singulier nous permettra d'obtenir des r\'esultats sur la r\'egularit\'e des d\'eformations Einstein infinit\'esimales ; on montrera en particulier \`a la section \ref{reg} que la d\'eformation induite de la m\'etrique du lieu singulier est $C^\infty$.

\subsection{Etude de l'op\'erateur lin\'earis\'e}\label{secLap2}

Soit $u$ une d\'eformation infinit\'esimale Einstein. Si elle est normalis\'ee (c'est-\`a-dire si elle 
v\'erifie la condition de jauge de Bianchi $\beta(u)=0$), elle v\'erifie l'\'equation
$$\na \n u - 2 \RO u =0.$$
L'\'etude des d\'eformations Einstein est donc fortement reli\'ee \`a l'\'etude de l'op\'erateur
$P = \na \n - 2 \mathring{R}$, le lin\'earis\'e de l'\'equation d'Einstein pour une d\'eformation 
normalis\'ee. Dans le cas g\'en\'eral, son \'etude peut \^etre assez compliqu\'ee. Cependant si la forme 
quadratique $\<\RO u, u\>$ n'est pas trop positive, l'op\'erateur $P$ est coercif, ce qui permet 
d'obtenir des r\'esultats int\'eressants. Et les propri\'et\'es voulues de $\RO$ peuvent se d\'eduire de 
certaines hypoth\`eses de courbure sur la m\'etrique (voir \cite{Besse} \S 12.67 et 12.71).
Dans le cas qui nous int\'eresse ici, la m\'etrique est hyperbolique, et on a l'expression plus 
simple 
$$P u = \na \n u - 2 u + 2 (\tr u) g.$$

L'op\'erateur $P$ pr\'esente de nombreuses similarit\'es avec l'op\'erateur $L = \na \n + (n-1)Id$ 
agissant sur les $1$-formes, dont l'\'etude approfondie \'etait l'objet de la section \ref{secLap}. 
Il y a donc de nombreuses ressemblances dans le plan, les r\'esultats et les formulations entre 
les deux sections.

\bigskip
  
\subsubsection{Expression de l'op\'erateur en coordonn\'ees cylindriques}\label{cocyl}

On voudrait maintenant avoir plus de pr\'ecisions sur les solutions de l'\'equation $Pu=0$ au 
voisinage du lieu singulier. Pour cela, on va se placer en coordonn\'ees cylindriques au voisinage d'une 
composante connexe du lieu singulier.

\medskip

On utilisera les m\^emes notations que dans la section \ref{vois1}. Pour rappel,
les notations $e_1,\ldots e_{n-2}$ d\'esignent des champs de vecteurs locaux tels que
$(e_r,e_\theta,e_1,\ldots e_{n-2})$ forme un rep\`ere mobile orthonorm\'e (local),
v\'erifiant $\n_{e_r} e_k = \n_{e_\theta} e_k = 0$
pour tout $k$ dans $1\ldots n-2$. On d\'efinit de m\^eme des 1-formes locales
$e^1,\ldots e^{n-2}$ telles que $(e^r,e^\theta,e^1,\ldots e^{n-2})$ soit le rep\`ere mobile dual 
du pr\'ec\'edent. La notation $N$ d\'esigne le (sous-)fibr\'e vectoriel au-dessus de $U_a$, dont la 
fibre au-dessus de $x \in U_a$ est le sous-espace vectoriel de $T^*_x M$
orthogonal \`a $e^\theta$ et $e^r$, et $N^*$ d\'esigne le (sous-)fibr\'e vectoriel au-dessus
de $U_a$, dont la fibre au-dessus de $x \in U_a$ est le sous-espace vectoriel
de $T_x M$ orthogonal \`a $e_\theta$ et $e_r$. Les sections $(e_1,\ldots,
e_{n-2})$ forment localement une base de $N^*$, de m\^eme pour $(e^1,\ldots,
e^{n-2})$ et $N$. Si $s$ est une section de $N^*$,
et $t$ une section de $N$ ou de $N^*$, on note $\n_{\Sigma\,s}t$, ou de fa\c{c}on plus lisible 
$(\n_\Sigma)(s,t)$, la projection orthogonale sur $N$ ou sur $N^*$ de $\n_s t$.

\medskip

On introduit en plus le 
sous-fibr\'e $S^2N$, engendr\'e par les $e_i \odot e_j$, $1 \leq i \leq n-2$, $0 \leq j \leq n-2$.
(Si $a$ et $b$ sont deux $1$-formes, on note $a\odot b$, ou 
plus simplement s'il n'y a pas de risque de confusion,
$a.b$ ou $ab$ pour $\frac{1}{2}(a \otimes b + b \otimes a)$. En particulier,
$x\odot x = x \otimes x$.)
Si $k$ est une section de $S^2N$ et $s$ est une section de $N^*$, on d\'efinit de m\^eme $\n_\Sigma 
(s,k)$ comme \'etant la 
projection orthogonale sur $S^2N$ de la d\'eriv\'ee covariante $\n_s k$ prise dans $S^2(M)$.

Si $\omega$ est une section de $N$, on d\'efinit $\delta^*_\Sigma \omega$, section de $S^2N$, par
$$\delta^*_\Sigma \omega = \sum_{i = 1}^{n-2} e^i \odot \n_\Sigma (e_i,\omega).$$

Enfin, si $k$ est une section de $S^2M$, on d\'efinit $\delta_\Sigma k$, section de $N$, et 
$\trs k$ par
$$\delta_\Sigma k = -\cosh(r)^2 \sum_{i = 1}^{n-2} \n_\Sigma (e_i,k) (e_i)$$
et
$$\trs k = \cosh(r)^2 \sum_{i = 1}^{n-2} k(e_i,e_i).$$

\bigskip

Si $u$ est une section de $S^2M$, on peut la d\'ecomposer orthogonalement au-dessus de $U_a$ en
$$u = f e^r.e^r + g e^\theta.e^\theta + h e^r.e^\theta + \sigma.e^r +
\eta.e^\theta + k$$
o\`u $\sigma$ et $\eta$ sont des sections de $N$, et o\`u $k$ est une section de $S^2N$. 
Quelques calculs \'el\'ementaires permettent d'appliquer la m\^eme d\'ecomposition \`a $\na \n u - 2 \mathring{R}u = \na \n u - 2
u + 2 (\tr u) g$. En regroupant les termes on obtient alors, pour la
composante suivant $e^r.e^r$ :
\begin{multline*} \Delta f + 2\left(\frac{1}{\tanh(r)^2} +
(n-2)\tanh(r)^2 \right)f  + \left(2 - \frac{2}{ \tanh(r)^2}\right)g + \frac{2 e_\theta.h }{
\tanh(r)} - \frac{2 \tanh(r) }{ \cosh(r)^2} (\delta_\Sigma \sigma) 
+ \frac{2 - 2 \tanh(r) }{ \cosh(r)^2} ({\rm tr}_\Sigma\ k)\end{multline*}

pour la composante suivant $e^\theta.e^\theta$ :
$$ \left(2- \frac{2}{ \tanh(r)^2}\right)f + \Delta g + \frac{2 }{ \tanh(r)^2} g 
- \frac{2 e_\theta.h }{ \tanh(r)} + \frac{2
  }{ \cosh(r)^2} {\rm tr}_\Sigma\ k $$

pour la composante suivant $e^r.e^\theta$ :
$$ - 4 \frac{e_\theta.f }{ \tanh(r)} + 4 e_\theta.g \frac{1 }{ \tanh(r)} +
\Delta h + \left(\frac{4 }{ \tanh(r)^2} + (n - 2)\tanh(r)^2 -2 \right) h
- \frac{2 \tanh(r) }{ \cosh(r)^2} (\delta_\Sigma \eta)
$$

pour la composante incluse dans $N.e^r$ :
\begin{multline*} 
 - 4 \tanh(r) d_\Sigma f -\n_{e_r}\n_{e_r} \sigma - 
\n_{e_\theta}\n_{e_\theta} \sigma - \left(\frac{1 }{ \tanh(r)} + 
  (n-2)\tanh(r)\right) \n_{e_r} \sigma + \frac{1 }{ \cosh(r)^2} (\na
  \n)_\Sigma \sigma  \\
+ \left(\frac{1 }{ \tanh(r)^2} + (n + 1) \tanh(r)^2 - 2\right) \sigma + \frac{2 }{
  \tanh(r)} (\n_{e_\theta} \eta) - \frac{4 \tanh(r) }{ \cosh(r)^2} \delta_\Sigma
k
\end{multline*} 

pour la composante incluse dans $N.e^\theta$ :
\begin{multline*} 
 - 2 \tanh(r) d_\Sigma h -\n_{e_r}\n_{e_r} \eta - \n_{e_\theta}\n_{e_\theta}
 \eta - \left(\frac{1 }{ \tanh(r)} + (n-2)\tanh(r)\right) \n_{e_r} \eta + \frac{1
   }{ \cosh(r)^2} (\na \n)_\Sigma \eta  \\
+ \left(\frac{1 }{ \tanh(r)^2} + \tanh(r)^2 -2 \right) \eta - \frac{2 }{ \tanh(r)}
(\n_{e_\theta} \sigma)
\end{multline*} 

et enfin pour la composante incluse dans $S^2N$ :
\begin{multline*} 
\left(2 - 2 \tanh(r)^2\right)f\,\cosh(r)^2 g_\Sigma  +2 g - 2 \tanh(r) \delta^*_\Sigma \sigma  -\n_{e_r}
\n_{e_r} k - \n_{e_\theta} \n_{e_\theta} k - (\frac{1 }{ \tanh(r)} + (n-2)
\tanh(r)) \n_{e_r} k \\
+ \frac{1 }{ \cosh(r)^2} (\na \n)_\Sigma k  + \left(2 \tanh(r)^2 -2\right)k
+ 2 {\rm tr}_\Sigma\ k\, g_\Sigma
\end{multline*}

\bigskip

Pour pouvoir manipuler cette expression, on va effectuer dans la suite une sorte de 
d\'e\-com\-po\-si\-tion en s\'eries de Fourier g\'en\'eralis\'ees, c'est-\`a-dire une d\'ecomposition sur des vecteurs 
propres d'op\'erateurs elliptiques du second degr\'e. Mais il faut une d\'ecomposition suffisamment 
astucieuse pour qu'elle se comporte bien avec les op\'erateurs $d_\Sigma$, $\delta^*_\Sigma$, etc. 
qui apparaissent dans les expressions ci-dessus.

\subsubsection{D\'ecomposition en s\'erie de Fourier g\'en\'eralis\'ee}

Lors de l'\'etude de l'op\'erateur $\na \n + (n-1) Id$ agissant sur les $1$-formes, on va d\'ej\`a
eu besoin d'une d\'ecomposition bien choisie des espace $L^2(M)$ et $L^2(N)$. Ce qu'il nous faut 
maintenant est une base hilbertienne de $L^2(S^2N)$, dans laquelle le comportement des 
op\'erateurs $\delta^*_\Sigma$, $\trs$ etc. se comprennent bien.

On va utiliser la proposition suivante, qui compl\`ete la proposition \ref{bashib} (voir la 
section \ref{secbashib} pour certaines notations).

\medskip

\begin{Prop}
Il existe une base hilbertienne $(\psi_j)_{j\in \mathbb{N}}$ du complexifi\'e de
$L^2(\Sigma_a)$, telle que pour tout indice $j$, il existe un r\'eel $\lambda_j
\geq 0$ et un entier relatif $p_j$, pour lesquels
$$\begin{cases} \Delta_\Sigma \psi_j = \lambda_j \psi_j \\
e_\theta.\psi_j = \frac{ip_j\gamma}{\sinh(a)}\psi_j.
\end{cases}$$

\medskip

Soit $J$ l'ensemble des $j$ pour lesquels $\lambda_j > 0$. Il existe une base
hilbertienne $(\phi_j)_{j\in J } \cup (\varphi_j)_{j\in \mathbb{N}}$ du
complexifi\'e de $L^2(N)$, telle que :
\begin{itemize}
\item pour tout indice $j$ appartenant \`a $J$, $\phi_j = \frac{\cosh(a) }{
    (\lambda_j)^{1/2}} d_\Sigma \psi_j$, et donc
$$\begin{cases} (\na \n)_\Sigma \phi_j = (\lambda_j + n - 3) \phi_j \\
\n_{e_\theta}\phi_j = \frac{ip_j\gamma}{\sinh(a)}\phi_j\\
\delta_\Sigma \phi_j = \cosh(a) (\lambda_j)^{1/2} \psi_j;
\end{cases}$$
\item pour tout indice $j \in \mathbb{N}$, il existe un r\'eel $\mu_j$ et un
  entier relatif $p'_j$, pour lesquels
$$\begin{cases} (\na \n)_\Sigma \varphi_j = \mu_j \varphi_j \\
\n_{e_\theta}\varphi_j = \frac{ip'_j\gamma}{\sinh(a)}\varphi_j,
\end{cases}$$
et on a de plus $\delta_\Sigma \varphi_j = 0$.
\end{itemize}

\medskip

Il existe une base hilbertienne $(a_j)_{j\in \mathbb{N}} \cup (b_j)_{j\in J}
\cup (c_j)_{j\in \mathbb{N}} \cup (d_j)_{j\in \mathbb{N}}$ du
complexifi\'e de $L^2(S^2 N)$, telle que :
\begin{itemize}
\item pour tout indice $j \in \mathbb{N}$, $a_j = \frac{\psi_j }{ \sqrt{n-2}}
   \cosh(a)^2 g_\Sigma$, et donc 
   $$\begin{cases} (\na \n)_\Sigma a_j = \lambda_j a_j \\
   \n_{e_\theta} a_j = \frac{ip_j\gamma}{\sinh(a)} a_j\\
   {\rm tr}_\Sigma\ a_j = \sqrt{n-2} \cosh(a)^2 \psi_j \\ 
   \delta_\Sigma a_j = - \frac{\cosh(a)^2 }{ \sqrt{n-2}} d_\Sigma \psi_j
 \ ( = - \cosh(a) (\frac{\lambda_j }{ n-2})^{1/2} \phi_j {\rm \ si\ } \lambda_j
   \neq 0 );
   \end{cases}$$
\item pour tout indice $j \in J$, $b_j = \frac{\cosh(a) }{ \sqrt{n-3}}
   (\frac{\lambda_j }{ n-2} + 1)^{-1/2} \left(\delta^*_\Sigma \phi_j + \frac{1 }{
   n-2} (\delta_\Sigma \phi_j) g_\Sigma \right)$, et donc 
   $$\begin{cases} (\na \n)_\Sigma b_j = (\lambda_j + 2(n-2)) b_j \\
   \n_{e_\theta} b_j = \frac{ip_j\gamma}{\sinh(a)} b_j\\
   {\rm tr}_\Sigma\ b_j = 0 \\ 
   \delta_\Sigma b_j = \cosh(a) \sqrt{n-3} (\frac{\lambda_j }{ n-2} + 1)^{1/2}
   \phi_j; 
   \end{cases}$$
\item pour tout indice $j \in \mathbb{N}$, $c_j = \cosh(a) (\frac{\mu_j + n-3 }{
   2})^{-1/2}\, \delta^*_\Sigma \varphi_j$, et donc 
   $$\begin{cases} (\na \n)_\Sigma c_j = (\mu_j + n - 1) c_j \\
   \n_{e_\theta} c_j = \frac{ip'_j\gamma}{\sinh(a)} c_j\\
   {\rm tr}_\Sigma\ c_j = 0 \\ 
   \delta_\Sigma c_j = \cosh(a) (\frac{\mu_j + n-3 }{ 2})^{1/2} \varphi_j;
   \end{cases}$$
\item pour tout indice $j \in \mathbb{N}$, il existe un r\'eel $\nu_j \geq 0$ et
   un entier relatif $p''_j$, pour lesquels
   $$\begin{cases} (\na \n)_\Sigma d_j = \nu_j d_j \\
   \n_{e_\theta}\varphi_j = \frac{ip''_j\gamma}{\sinh(a)} d_j,
   \end{cases}$$
   et on a de plus 
   $$\begin{cases} \delta_\Sigma d_j = 0\\
   {\rm tr}_\Sigma\ d_j = 0.
   \end{cases}$$
\end{itemize}
\end{Prop}

\medskip

Cette proposition, comme la proposition \ref{bashib}, se d\'emontre \`a l'aide des r\'esultats sur la th\'eorie spectrale des submersions riemanniennes \`a fibres totalement g\'eod\'esiques, voir \cite{BessBo} et \cite{BouBeBer}. Les liens entre les fonctions, les $1$-formes et les $2$-tenseurs sym\'etriques sont donn\'es, en plus de la relation \eqref{0030}, par les relations de commutations suivantes, qui utilisent le fait que la m\'etrique sur le lieu singulier $\Sigma$ est hyperbolique : pour toute section $\omega$ de $N$,
$$\na \n_\Sigma (\delta^*_\Sigma \omega) = (n-1)\delta^*_\Sigma \omega + 2 (\delta_\Sigma 
\omega) g_\Sigma + \delta^*_\Sigma (\na \n_\Sigma \omega)$$
et
\begin{align*}\na \n_\Sigma ((\delta_\Sigma \omega) g_\Sigma) &= \Delta_\Sigma (\delta_\Sigma 
\omega) g_\Sigma\\
&= (\delta_\Sigma (\na \n_\Sigma \omega) - (n-3) (\delta_\Sigma \omega)) g_\Sigma.\end{align*}
La premi\`ere est la relation \eqref{0001} vue p.\pageref{0001}, appliqu\'ee \`a la m\'etrique 
$g_\Sigma$. Elle d\'ecoule directement de l'invariance par diff\'eomorphisme du tenseur de 
Ricci. Dans la deuxi\`eme relation, la premi\`ere \'egalit\'e est \'el\'ementaire, et la deuxi\`eme utilise le 
fait que $\Delta_\Sigma$ et $\delta_\Sigma$ commutent ainsi que la formule de Weitzenb\"ock 
\eqref{0015}. En combinant les deux relations de commutations on trouve une derni\`ere relation :
\begin{eqnarray*}
\na \n_\Sigma (\delta^*_\Sigma \omega + \frac{1 }{ n-2} (\delta_\Sigma \omega) g_\Sigma) & = &
(n-1)\delta^*_\Sigma \omega + 2 (\delta_\Sigma \omega) g_\Sigma +
\delta^*_\Sigma (\na \n_\Sigma \omega) \\ & & + \frac{1 }{ n-2}(\delta_\Sigma (\na \n_\Sigma
\omega) - (n-3) (\delta_\Sigma \omega)) g_\Sigma\\
& = & (n-1)(\delta^*_\Sigma \omega + \frac{1 }{ n-2} (\delta_\Sigma \omega) g_\Sigma) +  
\delta^*_\Sigma (\na \n_\Sigma \omega) \\ & &+ \frac{1 }{ n-2}(\delta_\Sigma (\na \n_\Sigma 
\omega)) g_\Sigma.
\end{eqnarray*}

\bigskip

Maintenant, comme dans la section \ref{secbashib}, on prolonge ces sections \`a tout $U_a$ par 
transport parall\`ele le long des g\'eod\'esiques int\'egrales du champ de vecteurs $e_r$.
Par simple changement d'\'echelle, on observe que les tenseurs prolong\'es se comportent de la fa\c{c}on 
suivante sur le voisinage $U_a$ de $\Sigma$ :

 $$\begin{cases}\displaystyle \frac{\d}{\d r}\psi_j = 0\\
 \displaystyle \frac{\d}{\d \theta}\psi_j = ip_j\gamma\,\psi_j\\
 \displaystyle \Delta_\Sigma \psi_j = \lambda_j \psi_j\\
 \displaystyle d_\Sigma \psi_j = \frac{(\lambda_j)^{1/2}}{\cosh(r)}
 \phi_j {\rm \ (ou\ } 0 {\rm \ si\ } \lambda_j = 0 {\rm )}\\
 \displaystyle \psi_j g_\Sigma = \frac{(n-2)^{1/2} }{ \cosh(r)^2} a_j,
 \end{cases}$$
$$\begin{cases}
 \displaystyle \n_{\frac{\d}{\d r}} \phi_j = 0\\
 \displaystyle \n_{\frac{\d}{\d \theta}} \phi_j = ip_j\gamma\,\phi_j \\
 \displaystyle (\na \n)_\Sigma \phi_j =  (\lambda_j + n-3) \phi_j \\
 \displaystyle \delta_\Sigma \phi_j = \cosh(r) (\lambda_j)^{1/2} \psi_j\\
 \displaystyle \delta^*_\Sigma \phi_j = \cosh(r)^{-1} (n-3)^{1/2} (\frac{\lambda_j
 }{ n-2} +1)^{1/2} b_j - \cosh(r)^{-1}(\frac{\lambda_j }{ n-2})^{1/2}a_j,
\end{cases}$$
$$\begin{cases} \displaystyle \n_{\frac{\d}{\d r}} \varphi_j = 0\\
 \displaystyle \n_{\frac{\d}{\d \theta}} \varphi_j = ip'_j\gamma\,\varphi_j\\
 \displaystyle (\na \n)_\Sigma \varphi_j = \mu_j \varphi_j\\
\displaystyle  \delta_\Sigma \varphi_j = 0\\
 \displaystyle \delta^*_\Sigma \varphi_j = \cosh(r)^{-1} (\frac{\mu_j + n-3 }{
 2})^{1/2}\, c_j, 
\end{cases}$$
$$\begin{cases} \displaystyle \n_{\frac{\d}{\d r}} a_j = 0\\
 \displaystyle \n_{\frac{\d}{\d \theta}} a_j = ip_j\gamma\,a_j\\
 \displaystyle (\na \n)_\Sigma a_j = \lambda_j a_j\\
 \displaystyle  \delta_\Sigma a_j =  - \cosh(r) (\frac{\lambda_j }{ n-2})^{1/2}
 \phi_j {\rm \ (ou\ } 0 {\rm \ si\ } \lambda_j = 0 {\rm )}\\ 
 \displaystyle {\rm tr}_\Sigma\ a_j = \sqrt{n-2} \cosh(r)^2 \psi_j,
\end{cases}$$
$$\begin{cases} \displaystyle \n_{\frac{\d}{\d r}} b_j = 0\\
 \displaystyle \n_{\frac{\d}{\d \theta}} b_j = ip_j\gamma\,b_j\\
 \displaystyle (\na \n)_\Sigma b_j = (\lambda_j + 2(n-2)) b_j \\
  \displaystyle  \delta_\Sigma b_j = \cosh(r) \sqrt{n-3} (\frac{\lambda_j }{ n-2}
 + 1)^{1/2} \phi_j\\
 \displaystyle {\rm tr}_\Sigma\ b_j = 0,
\end{cases}$$
$$\begin{cases} \displaystyle \n_{\frac{\d}{\d r}} c_j = 0\\
 \displaystyle \n_{\frac{\d}{\d \theta}} c_j = ip'_j\gamma\,c_j\\
 \displaystyle (\na \n)_\Sigma c_j = (\mu_j + n - 1) c_j \\
 \displaystyle  \delta_\Sigma c_j = \cosh(r) (\frac{\mu_j + n-3 }{ 2})^{1/2}
 \varphi_j\\
 \displaystyle {\rm tr}_\Sigma\ c_j = 0,
\end{cases}$$
et enfin
$$\begin{cases} \displaystyle \n_{\frac{\d}{\d r}} d_j = 0\\
 \displaystyle \n_{\frac{\d}{\d \theta}} d_j = ip''_j\gamma\,d_j\\
 \displaystyle (\na \n)_\Sigma d_j = \nu_j d_j\\
 \displaystyle  \delta_\Sigma d_j = 0\\
 \displaystyle {\rm tr}_\Sigma\ d_j = 0
\end{cases}.$$

L'existence de ces bases hilbertiennes va permettre de r\'eduire l'\'equation aux d\'eriv\'ees 
partielles $Pu=0$ en une infinit\'e d'\'equations diff\'erentielles ordinaires.
On rappelle la d\'ecomposition orthogonale d'un $2$-tenseur sym\'etrique $u$ au voisinage d'une 
composante connexe du lieu singulier :
$$u = f e^r.e^r + g e^\theta.e^\theta + h e^r.e^\theta + \sigma.e^r +
\eta.e^\theta + k$$
o\`u $\sigma$ et $\eta$ sont des sections de $N$, et o\`u $k$ est une section de $S^2N$. En utilisant les r\'esultats de la section 
pr\'ec\'edente, on peut \'ecrire :
\begin{eqnarray*}
u & = & \sum_{j\in \mathbb{N}} f_j(r)\psi_j e^r.e^r + \sum_{j\in \mathbb{N}}
g_j(r)\psi_j e^\theta.e^\theta + \sum_{j\in \mathbb{N}} h_j(r)\psi_j
e^r.e^\theta \\
& + & \sum_{j\in J} \sigma_j(r)\phi_j.e^r + \sum_{j\in \mathbb{N}}
\overline{\sigma}_j(r) \varphi_j.e^r \\
& + & \sum_{j\in J} \eta_j(r)\phi_j.e^\theta + \sum_{j\in \mathbb{N}}
\overline{\eta}_j(r) \varphi_j.e^\theta \\
& + & \sum_{j\in \mathbb{N}} k^1_j(r) a_j + \sum_{j\in J} k^2_j(r) b_j
+ \sum_{j\in \mathbb{N}} k^3_j(r) c_j + \sum_{j\in \mathbb{N}} k^4_j(r) d_j.
\end{eqnarray*}

Il est plus judicieux de regrouper les termes de cette d\'ecomposition de la
fa\c{c}on suivante, faisant appara\^itre des ``blocs \'el\'ementaires'' de m\^eme
fr\'equence :
\begin{eqnarray}
u & = & \sum_{j\in J}\left(f_j(r)\psi_j e^r.e^r + g_j(r)\psi_j e^\theta.e^\theta + h_j(r)\psi_j 
e^r.e^\theta + \sigma_j(r)\phi_j.e^r +  \eta_j(r)\phi_j.e^\theta \right. \nonumber \\ 
& & \left. \mbox{} + k^1_j(r) a_j +  k^2_j(r) b_j \right) \nonumber \\
& + & \sum_{j\in \mathbb{N}\setminus J} \left(f_j(r)\psi_j e^r.e^r +
  g_j(r)\psi_j e^\theta.e^\theta + h_j(r)\psi_j e^r.e^\theta + k^1_j(r) a_j
  \right)\nonumber \\ 
& + & \sum_{j\in \mathbb{N}} \left( \overline{\sigma}_j(r) \varphi_j.e^r + 
\overline{\eta}_j(r) \varphi_j.e^\theta +  k^3_j(r) c_j \right)\nonumber \\ 
& + & \sum_{j\in \mathbb{N}} k^4_j(r) d_j. \label{0036}
\end{eqnarray}

Maintenant, on fait la m\^eme d\'ecomposition pour $\na \n u - 2 u + 2 (\tr u) g$, dont l'expression 
en coordonn\'ees cylindriques est donn\'ee \`a la section \ref{cocyl}.
On obtient alors, pour la composante en $\psi_j e^r.e^r$, si $j\in J$ :
\begin{multline} 
- f_j'' - \left(\frac{1 }{ \tanh(r)} + (n-2) \tanh(r)\right) f_j' + \left( \frac{2
 }{ \tanh(r)^2} + 2(n-2) \tanh(r)^2 + \frac{p_j^2 \gamma^2 }{ \sinh(r)^2} +
\frac{\lambda_j }{ \cosh(r)^2}\right) f_j \\
+ (2 - \frac{2 }{ \tanh(r)^2}) g_j  + \frac{2 i p_j \gamma }{ \sinh(r) \tanh(r)}
h_j -\frac{2 \tanh(r) (\lambda_j)^{1/2} }{ \cosh(r)} \sigma_j 
+ (2 - 2 \tanh(r)^2) (n-2)^{1/2} k^1_j, \label{0034}
\end{multline} 

si $j \notin J$ :
\begin{multline} 
- f_j'' - \left(\frac{1 }{ \tanh(r)} + (n-2) \tanh(r)\right) f_j' + \left( {2
 \over \tanh(r)^2} + 2(n-2) \tanh(r)^2 + {p_j^2 \gamma^2 \over \sinh(r)^2}
\right) f_j  \\
+ (2 - {2 \over \tanh(r)^2}) g_j  + {2 i p_j \gamma \over \sinh(r) \tanh(r)}
h_j + (2 - 2 \tanh(r)^2) (n-2)^{1/2} k^1_j, 
\end{multline} 

pour la composante en $\psi_j e^\theta.e^\theta$ :
\begin{multline} 
- g_j'' - \left({1 \over \tanh(r)} + (n-2) \tanh(r)\right) g_j' + \left( {2
 \over \tanh(r)^2} + {p_j^2 \gamma^2 \over \sinh(r)^2} +
{\lambda_j \over \cosh(r)^2}\right) g_j \\
+ (2 - {2 \over \tanh(r)^2}) f_j  - {2 i p_j \gamma \over \sinh(r) \tanh(r)}
h_j + 2 (n-2)^{1/2} k^1_j, 
\end{multline} 

pour la composante en $\psi_j e^r.e^\theta$, si $j\in J$ :
\begin{multline} 
- h_j'' - \left({1 \over \tanh(r)} + (n-2) \tanh(r)\right) h_j' + \left( {4
 \over \tanh(r)^2} + (n-2) \tanh(r)^2 + {p_j^2 \gamma^2 \over \sinh(r)^2} +
{\lambda_j \over \cosh(r)^2} - 2 \right) h_j \\
- {4 i p_j \gamma \over \sinh(r) \tanh(r)} f_j + {4 i p_j \gamma \over \sinh(r)
  \tanh(r)} g_j - {2 \tanh(r) (\lambda_j)^{1/2} \over \cosh(r)} \eta_j, 
\end{multline} 

si $j \notin J$ :
\begin{multline} 
- h_j'' - \left({1 \over \tanh(r)} + (n-2) \tanh(r)\right) h_j' + \left( {4
 \over \tanh(r)^2} + (n-2) \tanh(r)^2 + {p_j^2 \gamma^2 \over \sinh(r)^2} -2
\right) h_j  \\
- {4 i p_j \gamma \over \sinh(r) \tanh(r)} f_j + {4 i p_j \gamma \over \sinh(r)
  \tanh(r)} g_j, 
\end{multline} 

pour la composante en $\phi_j.e^r$ :
\begin{multline} 
- \sigma_j'' - \left({1 \over \tanh(r)} + (n-2) \tanh(r)\right) \sigma_j' +
\left( 
  {1 \over \tanh(r)^2} + (n+1) \tanh(r)^2 + {p_j^2 \gamma^2 \over \sinh(r)^2}
+ {\lambda_j + n-3 \over \cosh(r)^2} - 2 \right) \sigma_j  \\
- {4 \tanh(r) (\lambda_j)^{1/2} \over \cosh(r)} f_j 
+ {2 i p_j \gamma \over \sinh(r) \tanh(r)} \eta_j 
+ {4 \tanh(r) \over \cosh(r)} ({\lambda_j \over n-2})^{1/2} k^1_j 
- {4 \tanh(r) (n-3)^{1/2} \over \cosh(r)} ({\lambda_j \over n-2} + 1)^{1/2}
k^2_j,
\end{multline} 

pour la composante en $\phi_j.e^\theta$ :
\begin{multline} 
- \eta_j'' - \left({1 \over \tanh(r)} + (n-2) \tanh(r)\right) \eta_j' +
\left( 
  {1 \over \tanh(r)^2} + \tanh(r)^2 + {p_j^2 \gamma^2 \over \sinh(r)^2}
+ {\lambda_j + n-3 \over \cosh(r)^2} - 2 \right) \eta_j  \\
- {2 \tanh(r) (\lambda_j)^{1/2} \over \cosh(r)} h_j 
- {2 i p_j \gamma \over \sinh(r) \tanh(r)} \sigma_j,
\end{multline} 

pour la composante en $\varphi_j.e^r$ :
\begin{multline} 
- \overline{\sigma}_j'' - \left({1 \over \tanh(r)} + (n-2) \tanh(r)\right)
\overline{\sigma}_j' + \left({1 \over \tanh(r)^2} + (n+1) \tanh(r)^2 +
  {{p'_j}^2 \gamma^2 \over \sinh(r)^2} + {\mu_j \over \cosh(r)^2} - 2 \right) 
\overline{\sigma}_j \\
+ {2 i p'_j \gamma \over \sinh(r) \tanh(r)} \overline{\eta}_j 
- {4 \tanh(r) \over \cosh(r)} ({\mu_j + n-3 \over 2})^{1/2} k^3_j,
\end{multline} 

pour la composante en $\varphi_j.e^\theta$ :
\begin{multline} 
- \overline{\eta}_j'' - \left({1 \over \tanh(r)} + (n-2) \tanh(r)\right)
\overline{\eta}_j' + \left({1 \over \tanh(r)^2} + \tanh(r)^2 + {{p'_j}^2
 \gamma^2 \over \sinh(r)^2} + {\mu_j \over \cosh(r)^2} - 2 \right)
\overline{\eta}_j \\
- {2 i p'_j \gamma \over \sinh(r) \tanh(r)} \overline{\sigma}_j,
\end{multline} 

pour la composante en $a_j$, si $j \in J$ :
\begin{multline} 
- {k^1_j}'' - \left({1 \over \tanh(r)} + (n-2) \tanh(r)\right) {k^1_j}' +
\left(2 \tanh(r)^2 + {p_j^2 \gamma^2 \over \sinh(r)^2} +
{\lambda_j \over \cosh(r)^2} - 2 + 2(n-2)\right) k^1_j \\
+ 2 (n-2)^{1/2}(1 - \tanh(r)^2) f_j + 2 (n-2)^{1/2} g_j + {2 \tanh(r) \over
  \cosh(r)}({\lambda_j \over n-2})^{1/2} \sigma_j,
\end{multline} 

si $j \notin J$ :
\begin{multline} 
- {k^1_j}'' - \left({1 \over \tanh(r)} + (n-2) \tanh(r)\right) {k^1_j}' +
\left(2 \tanh(r)^2 + {p_j^2 \gamma^2 \over \sinh(r)^2} - 2 + 2(n-2)\right) k^1_j \\
+ 2 (n-2)^{1/2}(1 - \tanh(r)^2) f_j + 2 (n-2)^{1/2} g_j,
\end{multline} 

pour la composante en $b_j$ :
\begin{multline} 
- {k^2_j}'' - \left({1 \over \tanh(r)} + (n-2) \tanh(r)\right) {k^2_j}' +
\left(2 \tanh(r)^2 + {p_j^2 \gamma^2 \over \sinh(r)^2} +
{\lambda_j + 2(n-2) \over \cosh(r)^2} - 2 \right) k^2_j \\
- {2 \tanh(r) \over \cosh(r)}(n-3)^{1/2}({\lambda_j \over n-2}+1)^{1/2}
\sigma_j, 
\end{multline} 

pour la composante en $c_j$ :
\begin{multline} 
- {k^3_j}'' - \left({1 \over \tanh(r)} + (n-2) \tanh(r)\right) {k^3_j}' +
\left(2 \tanh(r)^2 + {{p'_j}^2 \gamma^2 \over \sinh(r)^2} +
{\mu_j + n-1 \over \cosh(r)^2} - 2 \right) k^3_j \\
- {2 \tanh(r) \over \cosh(r)}({\mu_j + n - 3 \over 2})^{1/2}
\overline{\sigma}_j,
\end{multline} 

et pour la composante en $d_j$ :
\begin{equation}
- {k^4_j}'' - \left({1 \over \tanh(r)} + (n-2) \tanh(r)\right) {k^4_j}' +
\left(2 \tanh(r)^2 + {{p''_j}^2 \gamma^2 \over \sinh(r)^2} +
{\nu_j \over \cosh(r)^2} - 2 \right) k^4_j.\label{0035}
\end{equation}

\bigskip

\subsubsection{Comportement des solutions de l'\'equation homog\`ene}\label{eqhom2}

Comme \`a la section \ref{vois3}, on constate que pour chaque indice $j$, l'\'equation (ou plut\^ot le 
syst\`eme) que l'on obtient pr\'esente une 
singularit\'e ``r\'eguli\`ere'' en $r=0$. Les solutions de ces
\'equations sont donc des combinaisons lin\'eaires de fonctions de la forme
$r^\kappa f(r)$ avec $f$ une fonction analytique, o\`u les exposants $\kappa$ s'obtiennent
comme racines de l'\'equation indicielle (en cas de racines multiples ou
s\'epar\'ees par des entiers, il faut \'eventuellement
rajouter des termes en $\ln r$ dans l'expression des solutions).

On pose donc, pour un indice $j$ donn\'e,
\begin{align*}
f_j(r) & = r^\kappa(f_0 + f_1 r + f_2 r^2 + \cdots), \\
g_j(r) & = r^\kappa(g_0 + g_1 r + g_2 r^2 + \cdots), \\
{\rm \ etc.} &
\end{align*}

A partir des expressions \eqref{0034} \`a \eqref{0035}, on aboutit aux syst\`emes d'\'equations indicielles suivants : si $j \in J$,
$$\left\{ 
\begin{array}{rcrcrcl}
(-\kappa^2 + 2 + p_j^2\gamma^2) f_0 & - & 2 g_0 & + & 2ip_j\gamma h_0 & = & 0 \\
 -2 f_0 & + & (-\kappa^2 + 2 + p_j^2\gamma^2) g_0 & - & 2ip_j\gamma h_0 & = & 0 \\
 -4 ip_j\gamma f_0 & + & 4 ip_j\gamma g_0 & + & (-\kappa^2 + 4 + p_j^2\gamma^2) h_0 & =
 & 0 \\ 
& & (-\kappa^2 + 1 + p_j^2\gamma^2) \sigma_0 & + & 2 ip_j\gamma \eta_0 & = & 0 \\
& & -2ip_j\gamma \sigma_0 & + & (-\kappa^2 + 1 + p_j^2\gamma^2) \eta_0 & = & 0 \\
& & & & (-\kappa^2 + p_j^2\gamma^2) k^1_0 & = & 0 \\
& & & & (-\kappa^2 + p_j^2\gamma^2) k^2_0 & = & 0,
\end{array}\right.$$
si $j \notin J$,
$$\left\{ 
\begin{array}{rcrcrcl}
(-\kappa^2 + 2 + p_j^2\gamma^2) f_0 & - & 2 g_0 & + & 2ip_j\gamma h_0 & = & 0 \\
 -2 f_0 & + & (-\kappa^2 + 2 + p_j^2\gamma^2) g_0 & - & 2ip_j\gamma h_0 & = & 0 \\
 -4 ip_j\gamma f_0 & + & 4 ip_j\gamma g_0 & + & (-\kappa^2 + 4 + p_j^2\gamma^2) h_0 & =
 & 0 \\ 
& & & & (-\kappa^2 + p_j^2\gamma^2) k^1_0 & = & 0,
\end{array}\right.$$
et 
$$\left\{ 
\begin{array}{rcrcl}
(-\kappa^2 + 1 + {p_j'}^2\gamma^2) \overline{\sigma}_0 & + & 2 ip_j'\gamma
\overline{\eta}_0 & = & 0 \\
-2ip_j'\gamma  \overline{\sigma}_0 & + & (-\kappa^2 + 1 + {p_j'}^2\gamma^2)
\overline{\eta}_0 & = & 0 \\ 
 & & (-\kappa^2 + {p_j'}^2\gamma^2) k^3_0 & = & 0,
\end{array}\right.$$
et enfin,
$$
(-\kappa^2 + {p_j''}^2\gamma^2) k^4_0 = 0.$$

Un simple calcul de d\'eterminant donne maintenant les racines indicielles, c'est-\`a-dire les valeurs de $\kappa$ pour lesquelles les syst\`emes ci-dessus admettent des solutions non triviales.
On trouve que les exposants dominants sont de la forme $\pm p_j\gamma$, $\pm p_j\gamma \pm 1$ et $\pm p_j\gamma \pm 2$, $\pm p_j' \gamma$ et $\pm p_j'\gamma \pm 1$, et $\pm p_j''\gamma$.
La seule racine multiple posant probl\`eme est en $0$ (uniquement) si $p_j$, $p_j'$ ou $p_j''$ est nul (ce qui arrive toujours), ou si $p_j\gamma \in \{ -2,-1,1,2\}$ ou $p_j'\gamma \in \{-1,1\}$ (ce qui n'arrivera jamais avec nos
conditions d'angles). Il faut dans ces cas rajouter une solution logarithmique.

Plus pr\'ecis\'ement : \\
si $\kappa = \pm (p\gamma + 2)$, alors $(f_0,g_0,h_0)$ est multiple de
$(-1,1,2i)$ ;\\
si $\kappa = \pm (p\gamma - 2)$, alors $(f_0,g_0,h_0)$ est multiple de
$(1,-1,2i)$ ; \\ 
si $\kappa = \pm p\gamma$, alors $(f_0,g_0,h_0,k^1_0,k^2_0)$ est combinaison
lin\'eaire de $(1,1,0,0,0)$, $(0,0,0,1,0)$ et $(0,0,0,0,1)$ ; 
si $\kappa = \pm (p\gamma + 1)$, alors $(\sigma_0,\eta_0)$ est multiple de
$(1,-i)$ ; \\
si $\kappa = \pm (p\gamma - 1)$, alors $(\sigma_0,\eta_0)$ est multiple de
$(1,i)$.

\noindent Si $\kappa = \pm (p'\gamma + 1)$, alors
$(\overline{\sigma}_0,\overline{\eta}_0,k^3_0)$ est multiple de $(1,-i,0)$ ;\\
si $\kappa = \pm (p'\gamma - 1)$, alors
$(\overline{\sigma}_0,\overline{\eta}_0,k^3_0)$ est multiple de $(1,i,0)$ ;\\
si $\kappa = \pm p'\gamma$, alors
$(\overline{\sigma}_0,\overline{\eta}_0,k^3_0)$ est multiple de $(0,0,1)$.

\noindent Enfin, l'\'equation pour $k^4$ donne deux solutions d'exposants dominants $\kappa = \pm 
p''\gamma$.

\medskip

On obtient ainsi toute une famille de solutions \'el\'ementaires, que l'on pourrait regrouper dans un \'enonc\'e inutilement long, similaire \`a celui de la proposition \ref{dvpt}. Rappelons juste que si $u$ est solution de l'\'equation $Pu=0$ au voisinage d'une composante connexe du lieu singulier, alors chacun des termes de la d\'ecomposition en s\'erie \eqref{0036}
\begin{eqnarray*}
u & = & \sum_{j\in J}\left(f_j(r)\psi_j e^r.e^r + g_j(r)\psi_j e^\theta.e^\theta + h_j(r)\psi_j 
e^r.e^\theta + \sigma_j(r)\phi_j.e^r +  \eta_j(r)\phi_j.e^\theta \right. \\ 
& & \left. \mbox{} + k^1_j(r) a_j +  k^2_j(r) b_j \right)\\
& + & \sum_{j\in \mathbb{N}\setminus J} \left(f_j(r)\psi_j e^r.e^r +
  g_j(r)\psi_j e^\theta.e^\theta + h_j(r)\psi_j e^r.e^\theta + k^1_j(r) a_j
  \right)\\ 
& + & \sum_{j\in \mathbb{N}} \left( \overline{\sigma}_j(r) \varphi_j.e^r + 
\overline{\eta}_j(r) \varphi_j.e^\theta +  k^3_j(r) c_j \right)\\ 
& + & \sum_{j\in \mathbb{N}} k^4_j(r) d_j
\end{eqnarray*}
est une solution de l'\'equation $Pu=0$, et est une combinaison lin\'eaire des solutions \'el\'ementaires correspondantes.

\subsection{R\'egularit\'e des d\'eformations induites}\label{reg}

La connaissance des exposants dominants permet de comprendre ce qu'il se passe au voisinage du lieu singulier, quand $r$ tend vers $0$. En particulier, on va d\'emontrer la proposition suivante :

\medskip

\begin{Prop}\label{induit} Soit $u$ un $2$-tenseur sym\'etrique, appartenant \`a $L^{1,2}$ et solution de l'\'equation ${Pu=0}$ au voisinage d'une composante connexe $\Sigma$ du lieu singulier, dont l'angle conique $\alpha$ n'est pas un multiple de $\pi$. Alors $u$ induit sur $\Sigma$ un $2$-tenseur sym\'etrique $u_\Sigma$, qui est $C^\infty$.
\end{Prop}

\begin{proof} Le fait que soit $u$ soit dans $L^{1,2}$ impose des restrictions sur les solutions \'el\'ementaires apparaissant dans la d\'ecomposition de $u$. Comme dans la d\'emonstration du lemme \ref{nnd}, on trouve que les composants de $u$ sont combinaisons lin\'eaires des seules solutions \'el\'ementaires d'exposant dominant $\kappa \geq 0$, sans terme logarithmique. 

Par cons\'equent, tous les termes de la d\'ecomposition en s\'erie de $u$ sont born\'es, et les seuls termes ne tendant pas vers $0$ quand $r$ tend vers $0$ sont ceux d'exposant dominant $\kappa$ nul. On a vu que $\kappa$ \'etait de la forme $\pm p_j\gamma$, $\pm p_j\gamma \pm 1$, $\pm p_j\gamma \pm 2$, $\pm p_j' \gamma$, $\pm p_j'\gamma \pm 1$, ou $\pm p_j''\gamma$. Si l'angle conique $\alpha =  \frac{2\pi}{\gamma}$ n'est pas un multiple de $\pi$, les exposants dominants nuls sont donc ceux de la forme $\pm p_j \gamma$, $\pm p_j'\gamma$ et $\pm p_j''\gamma$ pour $p_j=0$, $p_j'=0$ et $p_j''=0$.
Il s'agit donc uniquement de termes ne d\'ependant pas de la variable d'angle $\theta$, ce qui permet de regarder leur limite quand $r$ tend vers $0$.

On peut alors consid\'erer le tenseur induit sur $\Sigma$, qui s'\'ecrit formellement comme une s\'erie 
\begin{equation}
u_\Sigma = \sum_{\setlength{\extrarowheight}{-12pt}\begin{array}{c} \scriptstyle j \in \N \\ \scriptstyle p_j =0\end{array}} \left(k^1_j(0) a_j + k^2_j(0) b_j\right) 
+ \sum_{\setlength{\extrarowheight}{-12pt}\begin{array}{c}  \scriptstyle j \in \N \\  \scriptstyle p_j' =0\end{array}} k^3_j(0) c_j
+ \sum_{\setlength{\extrarowheight}{-12pt}\begin{array}{c}  \scriptstyle j \in \N \\  \scriptstyle p_j'' =0\end{array}} k^4_j(0) d_j \label{0037}
\end{equation}
Notons que comme les sections de $S^2N$ $a_j$, $b_j$, $c_j$ et $d_j$ qui apparaissent ici ne d\'ependent pas la variable d'angle $\theta$, elles s'identifient naturellement avec des sections de $S^2\Sigma$, et forment une base hilbertienne de $L^2(S^2\Sigma)$.

On va maintenant montrer que la s\'erie \eqref{0037} d\'efinit bien un $2$-tenseur sym\'etrique $C^\infty$ sur $\Sigma$. De m\^eme que dans la d\'emonstration des lemmes \ref{nnd} et \ref{nnul2}, cela peut se faire de deux mani\`eres : soit en utilisant les r\'esultats de \cite{Mazzeo} sur les op\'erateurs d'ar\^etes elliptiques, soit en majorant les termes de la s\'erie, \`a la fa\c{c}on de \cite{These}. 

Cette deuxi\`eme m\'ethode est ici assez simple. Rappelons que pour $r$ assez petit, on note $\Sigma_r$ le bord d'un $r$-voisinage (tubulaire) de $\Sigma$. Il h\'erite naturellement d'une structure de fibr\'e en cercle sur $\Sigma$, dont la projection canonique sur $\Sigma$ est (\`a une homoth\'etie de facteur $\cosh(r)$ pr\`es) une submersion riemannienne \`a fibres totalement g\'eod\'esiques. La m\'etrique sur $\Sigma_r$, induite de celle de $M$, est localement celle d'un produit, et le fibr\'e tensoriel $S^2N$ est trivial le long des fibres de la submersion. Cela permet de consid\'erer l'int\'egrale le long des fibres de la submersion d'une section de $S^2N$ au-dessus de $\Sigma_r$ ; le r\'esultat est une section de $S^2\Sigma$.

Si on int\`egre ainsi, le long des fibres de la submersion, la composante selon $S^2N$ de la restriction \`a $\Sigma_r$ du tenseur $u$, on obtient une section $u_\Sigma(r)$ de $S^2\Sigma$, qui s'exprime sous la forme d'une s\'erie
$$u_\Sigma(r) = \alpha\, \sinh(r) \cosh(r)^2 \Bigg( \sum_{\setlength{\extrarowheight}{-12pt}\begin{array}{c} \scriptstyle j \in \N \\ \scriptstyle p_j =0\end{array}} \left(k^1_j(r) a_j + k^2_j(r) b_j\right) 
+ \sum_{\setlength{\extrarowheight}{-12pt}\begin{array}{c}  \scriptstyle j \in \N \\  \scriptstyle p_j' =0\end{array}} k^3_j(r) c_j
+ \sum_{\setlength{\extrarowheight}{-12pt}\begin{array}{c}  \scriptstyle j \in \N \\  \scriptstyle p_j'' =0\end{array}} k^4_j(r) d_j \Bigg),$$
les termes pour lesquels $p_j$ (ou $p_j'$, ou $p_j''$) est diff\'erent de $0$ disparaissant lors de l'int\'egration.

Or comme l'op\'erateur $P$ est elliptique, le tenseur $u$ est $C^\infty$, ainsi que la composante selon $S^2N$ de sa restriction \`a $\Sigma_r$. Le tenseur $u_\Sigma(r)$ est donc lui aussi $C^\infty$, ce qui implique en particulier une d\'ecroissance rapide de ces coefficients : pour tout polyn\^ome $P$,
les quantit\'es $|P(\lambda_j) k^1_j(r)|$ et $|P(\lambda_j) k^2_j(r)|$ tendent vers $0$ quand $j$ parcourt l'ensemble $\{ j \in \N \ | \ p_j=0\}$, et il en est de m\^eme pour $|P(\mu_j) k^3_j(r)|$ et $|P(\nu_j) k^4_j(r)|$ en rempla\c{c}ant $p_j$ par respectivement $p_j'$ et $p_j''$.

On utilise ensuite le fait que les solutions \'el\'ementaires d'exposant dominant nul, c'est-\`a-dire  pr\'ecis\'ement celles apparaissant dans l'expression en s\'erie de $u_\Sigma$ et $u_\Sigma(r)$, sont croissantes  en module sur un voisinage de $0$, sauf \'eventuellement pour un nombre fini d'indice $j$ ; les techniques \`a employer pour la d\'emonstration sont exactement celles de \cite{These}. Pour tout polyn\^ome $P$, on a donc $|P(\lambda_j) k^1_j(0)| \leq |P(\lambda_j) k^1_j(r)|$, sauf \'eventuellement pour un nombre fini d'indice $j$, ind\'ependent de $P$, et il en est de m\^eme avec $|P(\lambda_j) k^2_j(0)|$ etc.

On en d\'eduit que les coefficients de $u_\Sigma$ d\'ecroissent plus rapidement que tout polyn\^ome, ce qui d\'emontre que la s\'erie en question d\'efinit bien une section $C^\infty$ de $S^2\Sigma$.
\end{proof}

\bigskip

Cette proposition permet de d\'emontrer imm\'ediatement le th\'eor\`eme suivant :

\medskip

\begin{Thm} \label{regdef} Soit $M$ une c\^one-vari\'et\'e hyperbolique dont tous les angles coniques $\alpha_1,\ldots \alpha_p$ sont strictement inf\'erieurs \`a $\pi$. Soit $\dot{\alpha} = (\dot{\alpha_1}, \ldots 
\dot{\alpha_p})$ une variation donn\'ee du $p$-uplet des angles coniques, et soit $h_{\dot{\alpha}}$ la 
d\'eformation Einstein infinit\'esimale normalis\'ee correspondante. Alors la d\'eformation infinit\'esimale $h_\Sigma$ de la m\'etrique du lieu singulier, induite par $h_{\dot{\alpha}}$, est $C^\infty$.
\end{Thm}

\begin{proof} La d\'eformation Einstein infinit\'esimale $h_{\dot{\alpha}}$, telle qu'on l'a construite au th\'eor\`eme \ref{consdef}, est de la forme $h_1 - \delta^* \eta - h$. Le tenseur $h_1$ est asymptotique \`a une d\'eformation conforme de la m\'etrique du c\^one ; c'est lui qui r\'ealise la variation des angles coniques au premier ordre, mais il ne contribue pas \`a la d\'eformation induite $h_\Sigma$.

L'autre terme $h+\delta^*\eta$ v\'erifie, par construction, $P(h+\delta^*\eta)=0$ au voisinage du lieu singulier. Donc d'apr\`es la proposition \ref{induit} pr\'ec\'edente, il induit une d\'eformation $C^\infty$ de la m\'etrique du lieu singulier, dont la partie provenant de $h$ est a priori non triviale.
\end{proof}

\medskip

Les d\'eformations Einstein infinit\'esimales sont donc $C^\infty$ sur le lieu singulier. 
Notons n\'eanmoins qu'\`a cause du terme $h_1$, de la forme $(1 + \ln(r))(dr^2 + \sinh(r)^2 d\theta^2) + O(r^2\ln(r))$, certaines composantes divergent pr\`es du lieu singulier. Mais cela est en partie d\^u au fait que la d\'eformation est normalis\'ee : on peut trouver des d\'eformations non normalis\'ees, dans la m\^eme classe, qui reste born\'ees pr\`es de $\Sigma$ (il suffit de supprimer le terme en $\delta^*(\chi(r) e^r)$ qui appara\^it dans la construction). 
Notons aussi que, comme on pouvait s'y 
attendre, la d\'eformation $h_{\dot{\alpha}}$ ne se contente pas de changer les angles 
coniques : le fait de rester Einstein impose aussi de modifier de fa\c{c}on non triviale la m\'etrique du lieu singulier.

\vskip 1cm

Pour conclure cet article, je voudrais dire un mot des possibilit\'es d'int\'egration des d\'eformations Einstein infinit\'esimales modifiant les angles coniques. Le fait que, en un sens, l'espace tangent \`a une c\^one-vari\'et\'e hyperbolique parmi les structures de c\^ones-vari\'et\'es Einstein est aussi simple, param\'etr\'e par les variations du $p$-uplet des angles coniques, est encourageant, ainsi que le comportement r\'egulier des d\'eformations infinit\'esimales. Un autre point positif est donn\'e dans \cite{These} (th\'eor\`eme 3.2.4) : l'op\'erateur Einstein lin\'earis\'e pour les d\'eformations normalis\'ees, $P=\na\n -2 \RO$, est un isomorphisme de $L^{2,2}(S^2M)$ dans $L^2(S^2M)$ si tous les angles coniques sont inf\'erieurs \`a $2 \pi/ 3$. Tous ces r\'esultats vont dans le sens d'un th\'eor\`eme d'inversion locale.

 A cela s'oppose deux difficult\'es.
Premi\`erement, il y a de bonnes raisons de penser que le cadre des espaces de Sobolev (qui est celui de cet article) n'est pas le plus appropri\'e pour travailler avec des d\'eformations qui ne sont plus infinit\'esimales. Les espaces de H\"older sont mieux adapt\'es \`a ce genre de probl\`emes, mais les estim\'ees n\'ecessaires y sont plus difficiles \`a \'etablir.
Deuxi\`emement, les espaces fonctionnels consid\'er\'es (que ce soit de H\"older ou de Sobolev) d\'ependent de la valeur des angles coniques ; travailler avec des angles variables implique donc de travailler avec des espaces fonctionnels variables. Ces deux points contribuent \`a rendre le probl\`eme d'analyse g\'eom\'etrique plus compliqu\'e, mais c'est justement cela qui le rend plus int\'eressant.

\end{document}